\documentclass[smallextended]{svjour3}       
\smartqed  
\usepackage[pdftex,dvipsnames]{xcolor}  
\usepackage{graphicx}
\usepackage{subfigure}
\usepackage{amsmath,amssymb,stmaryrd}
\usepackage{tabu}
\usepackage[numbers]{natbib}

\usepackage{xargs}    
\usepackage[colorinlistoftodos,prependcaption,textsize=tiny]{todonotes}
\newcommandx{\unsure}[2][1=]{\todo[linecolor=blue,backgroundcolor=blue!25,bordercolor=blue,#1]{#2}}
\newcommandx{\change}[2][1=]{\todo[linecolor=red,backgroundcolor=red!25,bordercolor=red,#1]{#2}}
\newcommandx{\info}[2][1=]{\todo[linecolor=OliveGreen,backgroundcolor=OliveGreen!25,bordercolor=OliveGreen,#1]{#2}}
\newcommandx{\improvement}[2][1=]{\todo[linecolor=Plum,backgroundcolor=Plum!25,bordercolor=Plum,#1]{#2}}
\newcommandx{\thiswillnotshow}[2][1=]{\todo[disable,#1]{#2}}

\newcommand{\jump}[1]{\left\llbracket#1\right\rrbracket}   
\newcommand{\average}[1]{\left\{\!\!\left\{#1\right\}\!\!\right\}}   
\newcommand\interiorfaces{{\mathrm{interior}\atop\mathrm{faces}}}
\newcommand\boundaryfaces{{\mathrm{boundary}\atop\mathrm{faces}}}
\newcommand*\diff{\mathop{}\!\mathrm{d}}
\begin{document}

\title{A Free--Energy Stable Nodal Discontinuous Galerkin Approximation with Summation--By--Parts Property for the Cahn--Hilliard Equation}


\titlerunning{A Free--Energy Stable Nodal DG approximation with SBP Property for the CHE}        

\author{Juan Manzanero        \and Gonzalo Rubio \and David A. Kopriva \and Esteban Ferrer \and 
Eusebio Valero 
}


\institute{Juan Manzanero (\email{juan.manzanero@upm.es}), Gonzalo Rubio, Esteban Ferrer, Eusebio Valero  \at
	ETSIAE-UPM - School of Aeronautics, Universidad Polit\'ecnica de Madrid. Plaza Cardenal Cisneros 3, E-28040 Madrid, Spain. //
	Center for Computational Simulation, Universidad Polit\'ecnica de Madrid, Campus de Montegancedo, Boadilla del Monte, 28660, Madrid, Spain. \\
              David A. Kopriva \at Department of Mathematics, Florida State University and Computational Science Research Center, San Diego State University. 
}

\date{Received: date / Accepted: date}

\maketitle

\begin{abstract}
We present a nodal Discontinuous Galerkin (DG) scheme for the Cahn--Hilliard equation that satisfies the summation--by--parts simultaneous--approximation--term (SBP--SAT) property. 
The latter permits us to show that the discrete free--energy is bounded, and as a result, the scheme is provably stable.
 The scheme and the stability proof are presented for general curvilinear three--dimensional hexahedral meshes. We use the 
 Bassi--Rebay 1 (BR1) scheme to compute interface fluxes, and an IMplicit--EXplicit (IMEX) scheme to integrate in time. 
 Lastly, we test the theoretical findings numerically and present examples for two and three--dimensional problems.
\end{abstract}

\section{Introduction}\label{sec:Introduction}

Phase field models describe the phase separation dynamics of two immiscible liquids by minimizing a chosen free--energy. For an arbitrary free--energy function, it is possible to construct different phase field models. Amongst the most popular, one can find the Cahn--Hilliard and the Allen--Cahn models. The popularity of the first one, despite being a fourth order operator in space, comes from its ability to conserve phases \cite{2016:Lee}.

In this paper we present a nodal Discontinuous Galerkin (DG) spectral element method (DGSEM) for the Cahn--Hilliard equation \cite{1958:Cahn}. In particular, this work uses the Gauss--Lobatto version of the DGSEM, which makes it possible to obtain energy--stable schemes using the summation--by--parts simultaneous--aproximation--term 
(SBP--SAT) property. Moreover, it handles arbitrary three dimensional curvilinear hexahedral meshes whilst maintaining high--order spectral accuracy and free--energy stability. An alternative, which is not considered in this work, 
but has been studied in \cite{2007:Shu}, is to use consistent integration in all 
quadratures involved in the weak--formulation. However, that alternative yields a considerably more expensive solver when compared to the approach presented here.

The Gauss--Lobatto variant of the DGSEM has seen increased popularity in recent 
years. Although it is less accurate than its Gauss counterpart (for the same number of quadrature nodes), it
satisfies the SBP--SAT property, thus allowing one to construct schemes that are provably stable  \cite{2017:Kopriva}. Precisely, different authors have presented energy-- and entropy--stable schemes using this framework for the linear advection equation \cite{Kopriva2,2017:Manzanero}, 
Burgers equation \cite{2013:Gassner,gassner2017br1}, shallow water equations \cite{2016:Gassner:Shallow}, Euler and Navier--Stokes equations \cite{2016:Gassner,gassner2017br1}, and the magneto--hydrodynamics equations \cite{2016:Winters}, among others.

Following these ideas, we present a free--energy stable approximation for the 
Cahn--Hilliard equation. The stability analysis presented here is both semi--discrete (assuming exact integration in time) and fully--discrete (i.e. considering the discrete approximation of space and time). Note that previous work on DG schemes satisfying 
the SBP--SAT property use a semi--discrete energy analysis (i.e. continuous in 
time), without considering numerical errors introduced by the time discretization. 
This is because the equations previously considered are of second order at most 
in space, and hence can be efficiently integrated in time using explicit methods. Since the 
Cahn--Hilliard equation is fourth order in space, the numerical stiffness 
of the scheme leads to impractical time step limitations when using 
explicit methods. Therefore, we study the stability of an IMplicit--EXplicit (IMEX) approximation in time, which has been previously used for the Cahn--Hilliard equation \cite{2018:Dong}, and whose efficiency is similar to that of explicit methods. This is 
possible since the scheme is designed to have a constant (in time) coefficient matrix of the linear system,  which can be solved using LU factorization and 
Gauss elimination. 

The rest of this paper is organised as follows: In Sec. \ref{sec:CahnHilliard}, 
we introduce the Cahn--Hilliard equation and derive a continuous energy 
estimate. Next, in Sec. \ref{sec:DGSEM}, we construct the DG approximation. In 
Sec. \ref{sec:DiscreteAnalysis}, we perform the energy analysis in a semi--discrete fashion in Sec. \ref{subsec:SemidiscreteAnalysis} and 
 fully--discrete in Sec. \ref{subsec:FullyDiscreteIMEX}. Lastly, we provide numerical 
experiments in Sec. 
\ref{sec:NumericalExperiments} that assess the capabilities of the method.

\section{Cahn--Hilliard equation and continuous energy estimates}\label{sec:CahnHilliard}

In this section we give a brief description of the Cahn--Hilliard equation and 
its properties. The Cahn--Hilliard equation describes the phase separation 
dynamics of binary alloys or two phase flows. The phase field 
variable, $\phi$, satisfies the evolution equation
\begin{equation}
  \phi_t = \nabla\cdot\left(M\nabla w\right),~~\text{ in }\Omega,
  \label{eq:cahn--hilliard-eqn}
\end{equation}
where $M$ is a positive parameter named mobility, $\Omega$ is the physical domain (with boundaries $\partial\Omega$), and $w$ is a scalar field called \textit{chemical potential}, that is 
designed to minimize an arbitrary free--energy functional, $\mathcal F(\phi,\nabla\phi)$, which 
depends on the phase field and its gradients,
\begin{equation}
  w :=\frac{\delta \mathcal F}{\delta \phi}.
\end{equation}
For the chemical potential, $w$, we apply an homogeneous Neumann boundary condition to guarantee mass conservation,
\begin{equation}
\nabla w \cdot\vec{n}=0~~\text{ in } \partial\Omega.
\label{eq:continuous:BC-w}
\end{equation}

The free--energy is constructed so that two opposing effects 
balance: the chemical free--energy, $\psi$, which favors phase separation, and 
the interfacial energy $\frac{1}{2}k|\nabla\phi|^2$, which favors homogenization,
\begin{equation}
  \mathcal F = \int_{\Omega}\left(\psi(\phi) + \frac{1}{2}k|\nabla\phi|^2  
  \right)\diff\vec{x} - \int_{\partial \Omega}g(\phi)\diff S = F_v(\phi) + F_s(\phi).
  \label{eq:free-energy}
\end{equation}
In \eqref{eq:free-energy} we introduced $F_v(\phi)$ and $F_s(\phi)$  as the volumetric and surface free--energies respectively, where $g(\phi)$ represents a boundary energy that will also be minimized with appropriate 
boundary conditions, and $k$ is the interfacial energy coefficient. 

To perform the minimization, one linearizes the free--energy 
\eqref{eq:free-energy} around an equilibrium solution,
\begin{equation}
  \delta \mathcal F = \int_{\Omega}\left( \psi'(\phi)\delta \phi + k\nabla\phi\cdot\nabla\left(\delta 
  \phi\right)\right)\diff\vec{x}-\int_{\partial\Omega}g'(\phi)\delta\phi \diff 
  S,
  \label{eq:free-energy-minimization-1}
\end{equation}
where $\delta\phi$ is a small perturbation. Since we will also apply Neumann boundary conditions for $\phi$, the perturbation $\delta \phi$ is not restricted to vanish at the boundaries $\partial\Omega$. We integrate the second term of the first integral in \eqref{eq:free-energy-minimization-1} by parts,
\begin{equation}
  \delta \mathcal F = \int_{\Omega}\left( \psi'(\phi) - k\nabla^2\phi\right)\delta\phi \diff\vec{x}-\int_{\partial\Omega}\left(g'(\phi)-k\nabla\phi\cdot\vec{n}\right)\delta\phi 
  \diff S,
\end{equation}
which yields both the chemical potential definition
\begin{equation}
  w = \psi'(\phi) - k\nabla^2 \phi,
  \label{eq:chemical-potential-def}
\end{equation}
and the appropriate Neumann boundary conditions prescription

\begin{equation}
  k\nabla\phi\cdot\vec{n}\biggr|_{\partial\Omega} = g'(\phi).
  \label{eq:continuous:BC-phi}
\end{equation}

In this work we use the polynomial double--well function for the chemical 
free--energy \cite{1958:Cahn},
\begin{equation}
  \psi(\phi) = \frac{1}{4}(1-\phi)^2(1+\phi)^2,
  \label{eq:continuous:chemical-free-energy}
\end{equation}
and a linear function for the boundary energy,
\begin{equation}
g(\phi) = \beta \phi,
\end{equation}
but other choices that are not covered here exist (e.g. logarithmic chemical free--energy 
\cite{1995:Debussche}).

\subsection{Continuous free--energy stability bound}\label{subsec:continuousEnergyBound}

We first show that the free--energy is bounded by that computed with a given initial condition (i.e. the problem free--energy is well--posed). 
To do so, we follow \cite{2007:Shu} and transform the fourth order equation into a system of four first order equations. As a result, we construct four weak forms
\begin{subequations}\label{eq:continuous:weakform}
\begin{align}
\left\langle \phi_t, \varphi_\phi \right\rangle &= \left\langle \nabla\cdot \left(M\vec{f}\right),\varphi_\phi\right\rangle, \label{eq:continuous:weakform-phi}\\
\left\langle \vec{f}, \vec{\varphi}_f\right\rangle &= \left\langle \nabla w,\vec{\varphi}_f\right\rangle\label{eq:continuous:weakform-f},\\
\left\langle w, \varphi_w\right\rangle &= \left\langle \frac{\diff \psi}{\diff \phi}, \varphi_w\right\rangle - k\left\langle \nabla\cdot\vec{q}, \varphi_w\right\rangle, \label{eq:continuous:weakform-w}\\
\left\langle \vec{q}, \vec{\varphi}_{q} \right\rangle &= \left\langle \nabla\phi, \vec{\varphi}_q\right\rangle,\label{eq:continuous:weakform-q}
\end{align}
\end{subequations}
where we introduced the auxiliary variables $\vec{q}=\nabla\phi$ and $\vec{f}=\nabla w$, two arbitrary $L^2$ scalar test functions $\varphi_\phi$ and $\varphi_w$, and their two vectorial counterparts $\vec{\varphi}_f$ and $\vec{\varphi}_q$. The operator $\left\langle f,g \right\rangle$ is the $L^2$ inner product
\begin{equation}
\left\langle f, g\right\rangle = \int_{\Omega} f g \diff\vec{x},
\end{equation}
which induces the $L^2$ norm,
\begin{equation}
\left\langle f,f\right\rangle = \int_{\Omega} f^2 \diff\vec{x} = ||f||^2.
\end{equation}
We integrate \eqref{eq:continuous:weakform-phi} and \eqref{eq:continuous:weakform-w} by parts,
\begin{equation}
\begin{split}
\left\langle \phi_t, \varphi_\phi \right\rangle &= \int_{\partial \Omega} M\vec{f}\cdot\vec{n}\varphi_\phi \diff S - \left\langle M\vec{f}, \nabla\varphi_\phi \right\rangle, \\
\left\langle w, \varphi_w\right\rangle &= \left\langle \frac{\diff \psi}{\diff \phi}, \varphi_w\right\rangle - k\int_{\partial\Omega}\vec{q}\cdot\vec{n}\varphi_w \diff S + k\left\langle \vec{q}, \nabla \varphi_w\right\rangle,
\end{split}
\label{eq:continuous:weak-form-integrated}
\end{equation}
and apply the boundary conditions \eqref{eq:continuous:BC-phi} and \eqref{eq:continuous:BC-w} to \eqref{eq:continuous:weak-form-integrated},
\begin{subequations}\label{eq:continuous:weak-form-withBC}
\begin{align}
\left\langle \phi_t, \varphi_\phi \right\rangle &= - \left\langle M\vec{f}, \nabla\varphi_\phi\right\rangle, \label{eq:continuous:weak-form-withBC-phi}\\
\left\langle w, \varphi_w\right\rangle &= \left\langle \frac{\diff \psi}{\diff \phi},\varphi_w\right\rangle - k\int_{\partial\Omega}g'(\phi)\varphi_w \diff S + k\left\langle \vec{q}, \nabla \varphi_w\right\rangle.\label{eq:continuous:weak-form-withBC-w}
\end{align}
\end{subequations}
Next, we set $\varphi_\phi = w$ in \eqref{eq:continuous:weak-form-withBC-phi},
\begin{equation}
\left\langle \phi_t, w \right\rangle =- \left\langle M\vec{f}, \nabla w \right\rangle =- \left\langle M\vec{f}, \vec{f} \right\rangle \leqslant 0 ,
\label{eq:continuous:well-posedness-phi}
\end{equation}
where we used \eqref{eq:continuous:weakform-f} in the last equality. Similarly, we let $\varphi_w=\diff \phi / \diff t = \phi_t$ in \eqref{eq:continuous:weak-form-withBC-w},
\begin{equation}
\left\langle w, \phi_t\right\rangle = \left\langle \frac{\diff \psi}{\diff \phi}, \phi_t\right\rangle + k\left\langle \vec{q}, \nabla\phi_t\right\rangle - \int_{\partial\Omega}g'(\phi)\phi_t \diff S.
\label{eq:continuous:w-phit-equation}
\end{equation}
Lastly, we use the chain rule in time for each of the three terms in \eqref{eq:continuous:w-phit-equation},
\begin{equation}
\begin{split}
\left\langle\frac{\diff \psi}{\diff \phi},\phi_t\right\rangle &=\int_{\Omega}\frac{\diff \psi}{\diff \phi}\phi_t \diff\Omega =  \frac{\diff }{\diff t}\int_{\Omega}\psi \diff\Omega, \\
k\left\langle\vec{q},\nabla\phi_t\right\rangle &= k\left\langle \vec{q},\vec{q}_t\right\rangle = k\int_{\Omega}\vec{q}\cdot\vec{q}_t \diff\Omega = \frac{k}{2}\frac{\diff }{\diff t}\int_{\Omega}|\vec{q}|^2 \diff\Omega, \\
\int_{\partial\Omega}\frac{\diff g}{\diff \phi}\phi_t \diff S &= \frac{\diff }{\diff t}\int_{\partial\Omega} g(\phi)\diff S,
\end{split}
\end{equation}
which replaced in \eqref{eq:continuous:w-phit-equation} yields the time derivative of the free--energy $\mathcal F$,
\begin{equation}
\left\langle w,\phi_t\right\rangle = \frac{\diff }{\diff t}\int_{\Omega}\left(\psi + \frac{k}{2}\left|\vec{q}\right|^2\right)\diff\Omega - \frac{\diff }{\diff t}\int_{\partial\Omega}g(\phi)\diff S = \frac{\diff \mathcal 
F}{\diff t}.
\label{eq:continuous:well-posedness-w}
\end{equation}
Therefore, by substituting $\left\langle w,\phi_t\right\rangle$ from \eqref{eq:continuous:well-posedness-w} into \eqref{eq:continuous:well-posedness-phi} we find the bound for  the free--energy time derivative,

\begin{equation}
\frac{\diff  \mathcal F}{\diff t} = -\left\langle M\vec{f}, \vec{f} \right\rangle \leqslant 
0,
\label{eq:continuous:well-posedness-timeder}
\end{equation}
which we can integrate in time to show that
\begin{equation}
\mathcal F(T) = \mathcal F(0)-\int_{0}^T \left\langle M\vec{f}, \vec{f} \right\rangle dt \leqslant \mathcal F(0).
\label{eq:continuous:well-posedness}
\end{equation}
As a result, the problem \eqref{eq:cahn--hilliard-eqn} with chemical potential \eqref{eq:chemical-potential-def} and Neumann boundary conditions  \eqref{eq:continuous:BC-w} and \eqref{eq:continuous:BC-phi}  is well--posed in the sense that the free--energy $\mathcal F$, defined in \eqref{eq:free-energy}, is bounded in time, and is the property to be mimicked by the subsequent approximation.

\section{The nodal discontinuous Galerkin spectral element method}\label{sec:DGSEM}

In this section we describe the construction of the nodal Discontinuous Galerkin 
Spectral Element Method (DGSEM). From all the variants, we restrict 
ourselves to the tensor product DGSEM with Gauss--Lobatto (GL) points (DGSEM--GL), since it satisfies the 
summation--by--parts simultaneous--approximation--term (SBP--SAT) property. The latter is used to 
prove the scheme's stability without relying on
exact integration.

The computational domain $\Omega$ is tessellated with non--overlapping hexahedral 
elements, which are then geometrically transformed to a reference element $E=[-1,1]^3$ 
by means of a polynomial mapping that relates physical ($\vec{x}=(x^1, x^2, x^3)=(x,y,z)$) and local ($\vec{\xi}=(\xi^1, \xi^2, \xi^3)=(\xi,\eta,\zeta)$) 
coordinates through
\begin{equation}
  \vec{x} = \vec{X}(\vec{\xi})=\vec{X}(\xi,\eta,\zeta).
\end{equation}

Let $\vec{S}_{L/R}(\eta,\zeta)$, $\vec{S}_{F/Ba}(\xi,\zeta)$, and $\vec{S}_{T/Bo}(\xi,\eta)$ 
be order $N$ polynomial approximations of the (curvilinear) left, right, front, back, 
top, and bottom faces respectively of an element. We construct the transfinite mapping with a linear interpolation between those,
\begin{equation}
\begin{split}
  \vec{X} =&\phantom{{}+{}} \frac{1}{2}\vec{S}_{L}(\eta,\zeta)(1-\xi)+\frac{1}{2}\vec{S}_{R}(\eta,\zeta)(\xi+1)+\frac{1}{2}\vec{S}_{F}(\xi,\zeta)(1-\eta) + \frac{1}{2}\vec{S}_{Ba}(\xi,\zeta)(\eta+1) \\
  &+\frac{1}{2}\vec{S}_{Bo}(\xi,\eta)(1-\zeta) + \frac{1}{2}\vec{S}_{T}(\xi,\eta)(\zeta+1)-\frac{1}{4}\vec{S}_{Bo}(\xi,-1)(1-\eta)(1-\zeta) \\
  & -\frac{1}{4}\vec{S}_{Bo}(1,\eta)(1+\xi)(1-\zeta) - \frac{1}{4}\vec{S}_{Bo}(\xi,1)(\eta+1)(1-\zeta)-\frac{1}{4}\vec{S}_{Bo}(-1,\eta)(1-\xi)(1-\zeta)\\
  &-\frac{1}{4}\vec{S}_{T}(\xi,-1)(1-\eta)(1+\zeta)  -\frac{1}{4}\vec{S}_{T}(1,\eta)(1+\xi)(1+\zeta) - \frac{1}{4}\vec{S}_{T}(\xi,1)(\eta+1)(1+\zeta) \\
  &-\frac{1}{4}\vec{S}_{T}(-1,\eta)(1-\xi)(1+\zeta)-\frac{1}{4}\vec{S}_{L}(-1,\zeta)(1-\xi)(1-\eta)-\frac{1}{4}\vec{S}_{L}(1,\zeta)(1-\xi)(1+\eta)\\
  &-\frac{1}{4}\vec{S}_{R}(-1,\zeta)(1+\xi)(1-\eta)-\frac{1}{4}\vec{S}_{R}(1,\zeta)(1+\xi)(1+\eta)\\
  &+\frac{1}{8}\vec{S}_{Bo}(-1,-1)(1-\xi)(1-\eta)(1-\zeta)+\frac{1}{8}\vec{S}_{Bo}(1,-1)(1+\xi)(1-\eta)(1-\zeta)\\
  &+\frac{1}{8}\vec{S}_{Bo}(-1,1)(1-\xi)(1+\eta)(1-\zeta)+\frac{1}{8}\vec{S}_{Bo}(1,1)(1+\xi)(1+\eta)(1-\zeta)\\  
  &+\frac{1}{8}\vec{S}_{T}(-1,-1)(1-\xi)(1-\eta)(1+\zeta)+\frac{1}{8}\vec{S}_{T}(1,-1)(1+\xi)(1-\eta)(1+\zeta)\\
&+\frac{1}{8}\vec{S}_{T}(-1,1)(1-\xi)(1+\eta)(1+\zeta)+\frac{1}{8}\vec{S}_{T}(1,1)(1+\xi)(1+\eta)(1+\zeta).\\    
  \end{split}
\end{equation}

We use polynomial interpolants of order $N$ to approximate the solutions inside an element $E$. 
These polynomials are written as tensor products of the Lagrange interpolating polynomials, $l_j(\xi)$,
\begin{equation}
l_j(\xi) = \prod_{\substack{i = 0 \\ i \neq j }}^N \frac{\xi - \xi_i}{\xi_j - \xi_i},~~~j=0,...,N,
\end{equation}
whose nodes are a set of Gauss--Lobatto points $\{\xi_i\}_{i=0}^N$, $\{\eta_i\}_{i=0}^N$, and $\{\zeta_i\}_{i=0}^N$ in the reference element $E$. The polynomial interpolant of a function $u$ in $E$ is therefore
\begin{equation}
\mathcal I^{N}\left[u(x,y,z,t)\right]_{E} = U(\xi,\eta,\zeta,t) = \sum_{i,j,k=0}^N 
 U_{ijk}(t)l_i(\xi)l_j(\eta)l_k(\zeta).
 \label{eq:interpolant-definition}
\end{equation}
In \eqref{eq:interpolant-definition}, $U_{ijk}(t)$ represents the (time dependent) 
nodal values of an arbitrary function $u$. Note that we use lower cases for functions, whilst upper cases represent a polynomial interpolant. 

The Lagrange polynomials satisfy by construction the cardinal property
\begin{equation}
l_j(\xi_i) = \delta_{ij},
\end{equation}
where $\delta_{ij}$ is the Kronecker delta. To approximate integrals, we use a quadrature rule using GL nodes and weights $\{w_j\}_{j=0}^N$ which provides a precision of order $2N-1$ (see \cite{2009:Kopriva}),

\begin{equation}
\int_{-1}^1 F G \diff\xi \approx \int_{E,N} F G \diff\xi \equiv \sum_{m=0}^{N} w_m F_m  G_m.
\label{eq:quadrature-definition}
\end{equation}
The associated Lagrange polynomials are discretely orthogonal
\begin{equation}
\int_{E,N}l_i(\xi)l_j(\xi)\diff\xi = w_i \delta_{ij},
\end{equation}
and the definition of the quadrature weights $w_i$ for $i=j$ are

\begin{equation}
w_i =\int_{E,N}l_i(\xi)^2\diff\xi = \int_{E,N}l_i(\xi)\diff\xi =  \int_{-1}^{1} l_i(\xi)\diff\xi.
\label{eq:quadrature-weights-definition}
\end{equation}

Furthermore, we chose GL nodes since the quadrature rule \eqref{eq:quadrature-definition} with weights \eqref{eq:quadrature-weights-definition} satisfies the discrete summation--by--parts simultaneous--approximation--term (SBP--SAT) property, that is, the integration by parts rule holds discretely, which in one dimension is
\begin{equation}
\int_{E,N} \frac{\diff U}{\diff\xi} V\diff\xi = U_{N}V_{N} - U_{0}V_{0} - \int_{E,N} U \frac{\diff V}{\diff\xi} \diff\xi.
\end{equation}
The SBP--SAT property is used to follow the continuous analysis steps that prove the boundedness of the free energy $\mathcal F$ in \eqref{eq:continuous:well-posedness} discretely.

To construct the integrals that define the discrete weak--formulation of \eqref{eq:continuous:weakform} in a general curvilinear 3D configuration, we construct the covariant
\begin{equation}
\vec{a}_j = \frac{\partial\vec{X}}{\partial \xi^j},\quad j = 1,2,3,
\end{equation}
and contravariant
\begin{equation}
\vec{a}^j = \nabla \xi^j,\quad j = 1,2,3,
\end{equation}
vector bases. The covariant and (volume weighted) contravariant bases are related by
\begin{equation}
J\vec{a}^i = \vec{a}_j \times \vec{a}_k,\quad (i,j,k)\; \text{cyclic},
\label{eq:metrics:contravariant-covariant-rel}
\end{equation}
where the Jacobian of the transformation is
\begin{equation}
J = \vec{a}_1 \cdot\left(\vec{a}_2 \times \vec{a}_3\right).
\label{eq:metrics:Jacobian}
\end{equation}
However, the continuous metric identities
\begin{equation}
\sum_{i=1}^3 \frac{\partial Ja_n^i}{\partial \xi^i} = 0,\quad n=1,2,3
\label{eq:metrics:metric-identities}
\end{equation}
do not hold discretely for  \eqref{eq:metrics:contravariant-covariant-rel}, since the product $\vec{a}_i \times \vec{a}_j$ is a polynomial of order $2N$. Therefore, we construct the discrete contravariant basis in curl form
\begin{equation}
Ja^i_n = -\hat{\vec{x}}^i\cdot\nabla_\xi \times \mathcal I^{N}\left(X_l \nabla_\xi X_m\right),~~~i,n=1,2,3,~~~(n,m,l) \text{ 
cyclic},
\label{eq:metrics:curl-form}
\end{equation}
so that they satisfy the metric identities discretely \cite{2006:Kopriva}.
In \eqref{eq:metrics:curl-form}, $Ja_n^i$ is the $n-$th Cartesian component of the contravariant vector $J\vec{a}^i$, $\hat{\vec{x}}^i$ is the $i-$th Cartesian unit vector, and $\nabla_\xi =(\partial/\partial \xi, \partial/\partial \eta, \partial/\partial \zeta)$. 

We use the contravariant basis to transform differential operators from physical to computational space. The divergence of a vector is
\begin{equation}
\nabla\cdot\vec{F} = \frac{1}{J}\nabla_\xi\cdot \left(\boldsymbol{\mathcal M}^T\vec{F}\right) = \frac{1}{J}\nabla_\xi\cdot\tilde{\boldsymbol{F}},
\label{eq:metrics:transform-div}
\end{equation}
where,
\begin{equation}
\left[\tilde{F}^{1} ~\tilde{F}^{2} ~\tilde{F}^{3}\right]^T= \tilde{\boldsymbol{F}} = \boldsymbol{\mathcal M}^T\vec{F},
 \label{eq:metrics:contravariant-flux}
 \end{equation}
is the contravariant flux, and
\begin{equation}
\boldsymbol{\mathcal M} =[J\vec{a}^1 ~~ J\vec{a}^2 ~~ J\vec{a}^3]= [J\vec{a}^\xi ~~ J\vec{a}^\eta ~~ J\vec{a}^\zeta],
\end{equation}
is the Jacobian matrix of the transformation. The gradient of a scalar is
\begin{equation}
\nabla U = \frac{1}{J}\boldsymbol{\mathcal M}\nabla_\xi U.
\label{eq:metrics:transform-grad}
\end{equation}

Lastly, using the transformation Jacobian \eqref{eq:metrics:Jacobian}, we approximate the three dimensional integrals in an element using Gauss-Lobatto 
quadrature. Let $\mathcal J$ be the polynomial approximation of the mapping 
Jacobian \eqref{eq:metrics:Jacobian} (using \eqref{eq:interpolant-definition}), 
then the quadrature is,
\begin{equation}
\int_{e}F~ G\diff e \approx \int_{E,N}\mathcal J F~ G\diff E = \left\langle \mathcal J F, G\right\rangle_{E,N} \equiv \sum_{i,j,k=0}^N w_{ijk} \mathcal J_{ijk} F_{ijk}G_{ijk},
\label{eq:DGSEM:InnerProduct}
\end{equation}
with associated norm $||F||_{J,N}^2 = \left\langle \mathcal J F,F\right\rangle_{E,N}$. Eq. \eqref{eq:DGSEM:InnerProduct} allows us to write the discrete summation--by--parts property as in \cite{2017:Kopriva}
\begin{equation}
\left\langle \nabla_\xi U,\tilde{\boldsymbol{F}}\right\rangle_{E,N} = \int_{\partial E,N} U\tilde{\boldsymbol{F}}\cdot\hat{\boldsymbol{n}}\diff S_{\xi} - \left\langle \nabla_\xi\cdot \tilde{\boldsymbol{F}},U\right\rangle_{E,N}. 
\label{eq:DGSEM:SBP-property}
\end{equation}
In \eqref{eq:DGSEM:SBP-property}, $\hat{\vec{n}}$ is the reference space unit outward  normal vector at the element faces, $\diff S_\xi$ is the surface local integration variables ($\diff S^{i}_{\xi}=\pm \diff\xi^j\diff\xi^k$ for $i$--oriented faces).  
To compute the surface integral, we write two dimensional quadratures in each of 
the six faces that define the element,
\begin{equation}
\begin{split}
\int_{\partial E,N} U\tilde{\boldsymbol{F}}\cdot\hat{\boldsymbol{n}}\diff S_{\xi} = &\int_{f,N}U \tilde{F}^{\xi}\diff\eta \diff\zeta\biggr|_{\xi=-1}^{\xi=1} +  \int_{f,N}U \tilde{F}^{\eta}\diff\xi \diff\zeta\biggr|_{\eta=-1}^{\eta=1} + \\
& \int_{f,N}U \tilde{F}^{\zeta}\diff\xi \diff\eta\biggr|_{\zeta=-1}^{\zeta=1}.
\end{split}
\end{equation}

Moreover, we can write surface integrals in either physical or computational space. The relation to the physical surface integration variables is
\begin{equation}
\diff S^i = \left| J\vec{a}^i\right|\diff\xi^{j}\diff\xi^{k} = \mathcal J_f^i \diff S_\xi^i,
\end{equation}
where we defined the face Jacobian $\mathcal J_{f}^i = \left|\mathcal J \vec{a}^i\right|$. We can relate the surface flux in both reference element, $\tilde{\boldsymbol{F}}\cdot\hat{\boldsymbol{n}}$, and physical, $\vec{F}\cdot\vec{n}$, variables through
\begin{equation}
\tilde{\vec{F}}\cdot\hat{\vec{n}}^{i}\diff S_\xi = \left(\boldsymbol{\mathcal M}^T\vec{F}\right)\cdot\hat{\vec{n}}^{i}\diff S_\xi = \vec{F}\cdot\left(\boldsymbol{\mathcal M}\hat{\vec{n}}^{i}\right)\diff S_\xi = \vec{F}\cdot\vec{n}\left|J\vec{a}^i\right|\diff S_\xi = \vec{F}\cdot\vec{n}^{i}\diff S.
\label{eq:DG:local-physical-fluxes-relationship}
\end{equation}
Therefore, quadratures can be represented both in physical and computational 
spaces,
\begin{equation}
  \int_{\partial E,N}\tilde{\vec{F}}\cdot\hat{\vec{n}}\diff S_\xi = \int_{\partial 
  e,N}\vec{F}\cdot\vec{n}\diff S,
\end{equation}
and we will use one or the other depending on whether we are studying an 
isolated element (computational space) or the whole combination of elements in 
the mesh (physical space).

\subsection{Discontinuous Galerkin spectral element approximation of the Cahn--Hilliard equation}\label{subsec:DG-CH}

We now construct the discrete version of \eqref{eq:continuous:weakform}. 
We first transform \eqref{eq:continuous:weakform} to the local coordinate system as described in \eqref{eq:metrics:transform-div} and \eqref{eq:metrics:transform-grad}, and construct four weak forms inside the reference element $E$. We do so by multiplying each equation by test functions $\varphi_\Phi$, $\varphi_{W}$ (scalar), $\vec{\varphi}_{\vec{F}}$, $\vec{\varphi}_{\vec{Q}}$ (vectorial), which are restricted to the order $N$ polynomial space,
\begin{subequations}\label{eq:DG:weakforms}
\begin{align}
\left\langle J\phi_t,\varphi_\Phi\right\rangle_{E} &= \left\langle\nabla_\xi\cdot\left(M\tilde{\boldsymbol{f}}\right), \varphi_\Phi\right\rangle_{E}, \\
\left\langle  J\vec{f},\vec{\varphi}_F \right\rangle_{E} &= \left\langle\boldsymbol{\mathcal M}\nabla_\xi w,\vec{\varphi}_F\right\rangle_{E} = \left\langle\nabla_\xi w,\boldsymbol{\mathcal M}^T\vec{\varphi}_F\right\rangle_{E}= \left\langle\nabla_\xi w,\tilde{\boldsymbol{\varphi}}_F\right\rangle_{E},\\
\left\langle  J w,\varphi_W\right\rangle_{E} &= \left\langle  J \frac{\diff \psi}{\diff \phi} ,\varphi_W\right\rangle_{E} - k\left\langle\nabla_\xi\cdot\tilde{\boldsymbol{q}},\varphi_W\right\rangle_{E}, \\
\left\langle J\vec{q},\vec{\varphi}_Q\right\rangle_{E} &= \left\langle\boldsymbol{\mathcal M}\nabla_\xi\phi,\vec{\varphi}_Q\right\rangle_{E} = \left\langle\nabla_\xi\phi,\boldsymbol{\mathcal M}^T\vec{\varphi}_Q\right\rangle_{E}=\left\langle\nabla_\xi\phi,\tilde{\boldsymbol{\varphi}}_Q\right\rangle_{E}.
\end{align}
\end{subequations}
Next, we integrate the right hand side terms that contain a $\nabla_\xi$ operator by parts, replace the continuous functions by their polynomial approximations, and replace exact integrals by quadratures, 
\begin{subequations}\label{eq:DG:weakforms-with-surface}
\begin{align}
\left\langle \mathcal J\Phi_t,\varphi_\Phi\right\rangle_{E,N} &= \int_{\partial E,N}\varphi_\Phi\left(M\tilde{\boldsymbol{F}}\right)^\star\cdot\hat{\boldsymbol{n}} \diff S_\xi - \left\langle M\tilde{\boldsymbol{F}}, \nabla_\xi\varphi_\Phi\right\rangle_{E,N}, \label{eq:DG:weakforms-with-surface-phi}\\
\left\langle\mathcal J\vec{F},\vec{\varphi}_F \right\rangle_{E,N} &=  \int_{\partial E,N}W^\star \tilde{\boldsymbol{\varphi}}_F\cdot\hat{\boldsymbol{n}}\diff S_\xi - \left\langle W,\nabla_\xi\cdot\tilde{\boldsymbol{\varphi}}_F\right\rangle_{E,N},\label{eq:DG:weakforms-with-surface-f}\\
\left\langle \mathcal J W,\varphi_W\right\rangle_{E,N} &= \left\langle \mathcal J \frac{\diff \Psi}{\diff \Phi} ,\varphi_W\right\rangle_{E,N} -k\int_{\partial E,N}\varphi_W\tilde{\boldsymbol{Q}}^\star\cdot\hat{\boldsymbol{n}}\diff S_\xi + k\left\langle\tilde{\boldsymbol{Q}},\nabla_\xi\varphi_W\right\rangle_{E,N}, \label{eq:DG:weakforms-with-surface-w}\\
\left\langle \mathcal J\vec{Q},\vec{\varphi}_Q\right\rangle_{E,N} &= \int_{\partial E,N}\Phi^\star \tilde{\boldsymbol{\varphi}}_Q\cdot\hat{\boldsymbol{n}}\diff S_\xi -  \left\langle\Phi,\nabla_\xi\cdot\tilde{\boldsymbol{\varphi}}_Q\right\rangle_{E,N}.\label{eq:DG:weakforms-with-surface-q}
\end{align}
\end{subequations}

The terms with star superscript in \eqref{eq:DG:weakforms-with-surface} are the \textit{numerical fluxes}, which make the flux uniquely defined at the boundaries. The DG variant implemented depends on the choice of the numerical fluxes. In this work, we use the Bassi--Rebay 1 scheme (BR1), 
\begin{equation}
W^\star = \average{W}, \Phi^\star = \average{\Phi},
\label{eq:DG:BR1+IP-solution-riemann}
\end{equation}
where
\begin{equation}
\average{U} = \frac{U^{\partial e^+} + U^{\partial e^-}}{2}, ~~\average{\vec{F}} = \frac{\vec{F}^{\partial e^+} + \vec{F}^{\partial e^-}}{2},
\end{equation}
is the average operator. For the divergence weak forms, we use the average for $\vec{Q}^\star$ (i.e. the BR1 method). For $\left(M\vec{F}\right)^\star$ we propose to use the average with additional interface dissipation,
\begin{equation}
\left(M\vec{F}\right)^\star = \average{M\vec{F}}-\sigma M \jump{W}, \vec{Q}^\star = 
\average{\vec{Q}},
\label{eq:DG:BR1-div-riemann}
\end{equation}
where $\sigma$ is a positive penalty parameter ($\sigma\geqslant 0$). The jump operators (with built--in normal vectors) for a scalar $U$ and a vector $\vec{F}$ are defined as
\begin{equation}
\llbracket U \rrbracket = U\vec{n}\bigr|^{\partial e^+} + U\vec{n}\bigr|^{\partial e^-},~~ \llbracket \vec{F} \rrbracket = \vec{F}\cdot\vec{n}\bigr|^{\partial e^+} + \vec{F}\cdot\vec{n}\bigr|^{\partial e^-}.
\end{equation}
They satisfy the algebraic identity,
\begin{equation}
\llbracket U\vec{F}\rrbracket = \average{U}\llbracket \vec{F}\rrbracket + \llbracket U\rrbracket\cdot\average{\vec{F}}
\label{eq:dg:average-jump-rel}.
\end{equation}

For the penalty parameter $\sigma$, we propose to use the estimate introduced for the DGSEM--GL variant in \cite{2018:Manzanero},
\begin{equation}
\sigma=\kappa_{\sigma}\frac{N(N+1)}{2}|\mathcal J_{f}| \average{\mathcal J^{-1}},
\label{eq:PenaltyParameter}
\end{equation}
which has only one dimensionless free parameter, $\kappa_{\sigma}$ (also positive), and the dependency of $\sigma$ on the mesh and the polynomial order are taken into account automatically. In \eqref{eq:PenaltyParameter}, $N$ is the polynomial order, $|\mathcal J_{f}|$ is the surface Jacobian of the face, and $\average{\mathcal J^{-1}}$ is the average of the inverse of the Jacobians of the elements that share the face evaluated at the face.
It will be shown that the scheme is stable without interface stabilization. Nevertheless, the advantage of adding interface stabilization will be demonstrated with a numerical test. We have not introduced dissipation in $\vec{Q}^\star$ since it generates non--physical terms in the discrete 
free--energy, see Appendix \ref{app:StabilizationQ}.
%

For Neumann boundary conditions, we use the adjacent element interior value to 
compute gradients in \eqref{eq:DG:weakforms-with-surface-f} and \eqref{eq:DG:weakforms-with-surface-q},
\begin{equation}
  W^\star = W \bigr|_{\partial e},  ~~ \Phi^\star = \Phi\bigr|_{\partial e}, 
  \label{eq:DG:Neumann-bcs-wphi}
\end{equation}
and  directly impose Neumann boundary values for divergence weak forms \eqref{eq:DG:weakforms-with-surface-phi} and \eqref{eq:DG:weakforms-with-surface-w},
\begin{equation}
\left(M\vec{F}\right)^\star\cdot\vec{n} = 0,~~ \vec{Q}^\star\cdot\vec{n} = G'\left(\Phi\bigr|_{\partial 
e}\right)=\beta.
\label{eq:DG:Neumann-bcs-fq}
\end{equation}
Note that we have written the physical interface fluxes, rather than the contravariant fluxes, so that we can more easily relate the interface values shared by two elements (recall that the same physical flux yields different contravariant flux values, as it depends on each element geometry).

\section{Stability analysis}\label{sec:DiscreteAnalysis}

In this section, we will follow the steps in Sec. \ref{subsec:continuousEnergyBound} for the continuous
analysis to show the stability of the numerical scheme \eqref{eq:DG:weakforms-with-surface}. 
We first show in Sec. \ref{subsec:SemidiscreteAnalysis} the semi--discrete free--energy estimate 
assuming exact time integration. Then, in Sec. \ref{subsec:FullyDiscreteIMEX} we do the same for a 
fully discrete approximation in space and time. We cannot follow exactly the 
continuous steps in the semi--discrete analysis. To obtain the estimate we change the 
numerical scheme to one whose stability can be verified. Since the spatial approximations have to be 
different, we require two separate analyses to perform the semi--discrete and the 
fully--discrete stability proofs.

\subsection{Semi--discrete stability analysis}\label{subsec:SemidiscreteAnalysis}

In this section we follow the steps used to derive the free--energy bound 
 \eqref{eq:continuous:well-posedness} to derive an equivalent discrete bound for a modification of \eqref{eq:continuous:weakform-q} assuming exact time integration.  
To prove semi--discrete stability, we need to modify the original set of equations \eqref{eq:continuous:weakform}, replacing \eqref{eq:continuous:weakform-q} by its time derivative,
\begin{equation}
\left\langle\vec{q}_t,\vec{\varphi}_{\vec{Q}}\right\rangle_{E} = \left\langle\nabla\phi_t,\vec{\varphi}_{\vec{Q}}\right\rangle_{E},
\end{equation} 
whose associated discrete weak formulation is obtained following the same steps used for 
\eqref{eq:DG:weakforms-with-surface-q},
 \begin{equation}
\left\langle \mathcal J\vec{Q}_t,\vec{\varphi}_Q\right\rangle_{E,N} = \int_{\partial E,N}\Phi^\star_t \tilde{\boldsymbol{\varphi}}_Q\cdot\tilde{\boldsymbol{n}}\diff S_\xi -  \left\langle\Phi_t,\nabla_\xi\cdot\tilde{\boldsymbol{\varphi}}_Q\right\rangle_{E,N}.\label{eq:SemiDiscrete:weakforms-with-surface-q}
\end{equation}
 %
 
Now we follow the steps performed in the continuous energy analysis  \eqref{eq:continuous:well-posedness}. First, we need to sum \eqref{eq:SemiDiscrete:weakforms-with-surface-q} again by parts,
\begin{equation}
\left\langle \mathcal J\vec{Q}_t,\vec{\varphi}_Q\right\rangle_{E,N} = \int_{\partial E,N}\left(\Phi^\star_t - \Phi_t\right) \tilde{\boldsymbol{\varphi}}_Q\cdot\tilde{\boldsymbol{n}}\diff S_\xi + \left\langle\nabla_\xi\Phi_t,\tilde{\boldsymbol{\varphi}}_Q\right\rangle_{E,N},
\end{equation}
and replace $\vec{\varphi}_Q = \vec{Q}$ (and $\tilde{\vec{\varphi}}_{\vec{Q}} =\boldsymbol{\mathcal M}^{T}\vec{\varphi}_{\vec{Q}}$ as shown in \eqref{eq:DG:weakforms}), to get the time derivative,
\begin{equation}
\left\langle \mathcal J\vec{Q}_t,\vec{Q}\right\rangle_{E,N} = \frac{1}{2}\frac{\diff }{\diff t}||\vec{Q}||_{J,N}^2= \int_{\partial E,N}\left(\Phi^\star_t - \Phi_t\right) \tilde{\boldsymbol{Q}}\cdot\tilde{\boldsymbol{n}}\diff S_\xi + \left\langle\nabla_\xi\Phi_t,\tilde{\boldsymbol{Q}}\right\rangle_{E,N}.
\label{eq:SemiDiscrete:NormQ}
\end{equation}
Next, we set $\varphi_W = \Phi_t$ in 
\eqref{eq:DG:weakforms-with-surface-w},
\begin{equation}
\left\langle \mathcal J W,\Phi_t\right\rangle_{E,N} = \left\langle \mathcal J \frac{\diff \Psi}{\diff \Phi} ,\Phi_t\right\rangle_{E,N} -k\int_{\partial E,N}\Phi_t\tilde{\boldsymbol{Q}}^\star\cdot\hat{\boldsymbol{n}}\diff S_\xi  + k\left\langle\tilde{\boldsymbol{Q}},\nabla_\xi\Phi_t\right\rangle_{E,N}.
\label{eq:SemiDiscrete:W-equation-wQ}
\end{equation}
Since the time derivative is exact, we can use the chain rule for the chemical free--energy potential derivative,
\begin{equation}
\left\langle \mathcal J \frac{\diff \Psi}{\diff \Phi} ,\Phi_t\right\rangle_{E,N} = \frac{\diff }{\diff t}\left\langle \mathcal J \Psi, 1\right\rangle_{E,N}.
\label{eq:SemiDiscrete:ChainRuleChemPot}
\end{equation}

We subtract \eqref{eq:SemiDiscrete:NormQ} multiplied by $k$ from \eqref{eq:SemiDiscrete:W-equation-wQ}, and use the result in \eqref{eq:SemiDiscrete:ChainRuleChemPot} to obtain,
\begin{equation}
\left\langle \mathcal J W,\Phi_t\right\rangle_{E,N} = \frac{\diff }{\diff t}\left\langle \mathcal J \Psi, 1\right\rangle_{E,N} + \frac{k}{2}\frac{\diff }{\diff t}||Q||_{J,N}^2 -k\int_{\partial E,N}\left(\Phi_t\tilde{\boldsymbol{Q}}^\star\cdot\hat{\boldsymbol{n}}+\Phi^\star_t \tilde{\boldsymbol{Q}}\cdot\hat{\boldsymbol{n}}-\Phi_t \tilde{\boldsymbol{Q}}\cdot\hat{\boldsymbol{n}}\right)\diff S_\xi.
\label{eq:SemiDiscrete:W-equation}
\end{equation}
Next, we sum \eqref{eq:DG:weakforms-with-surface-f} by parts,
\begin{equation}
\left\langle\mathcal J\vec{F},\vec{\varphi}_F \right\rangle_{E,N} =  \int_{\partial E,N}\left(W^\star-W\right) \tilde{\boldsymbol{\varphi}}_F\cdot\hat{\boldsymbol{n}}\diff S_\xi + \left\langle \nabla_\xi 
W,\tilde{\boldsymbol{\varphi}}_F\right\rangle_{E,N},
\end{equation}
and we set $\vec{\varphi}_F = M\vec{F}$, so that
\begin{equation}
\left\langle\mathcal J\vec{F},M\vec{F} \right\rangle_{E,N} =\left|\left|\sqrt{M}\vec{F}\right|\right|^2_{J,N} = \int_{\partial E,N}\left(W^\star-W\right) M\tilde{\boldsymbol{F}}\cdot\hat{\boldsymbol{n}}\diff S_\xi + \left\langle \nabla_\xi W,M\tilde{\boldsymbol{F}}\right\rangle_{E,N}.
\label{eq:SemiDiscrete:NormF}
\end{equation}
Lastly, we set $\varphi_\phi = W$ in \eqref{eq:DG:weakforms-with-surface-phi},
\begin{equation}
\left\langle \mathcal J\Phi_t,W\right\rangle_{E,N} = \int_{\partial E,N}W\left(M\tilde{\boldsymbol{F}}\right)^\star\cdot\hat{\boldsymbol{n}} \diff S_\xi - \left\langle M\tilde{\boldsymbol{F}}, \nabla_\xi W\right\rangle_{E,N},
\label{eq:SemiDiscrete:W-eqn-weak}
\end{equation}
and we sum \eqref{eq:SemiDiscrete:NormF} and \eqref{eq:SemiDiscrete:W-eqn-weak} to find,
\begin{equation}
\left\langle \mathcal J\Phi_t,W\right\rangle_{E,N} = \int_{\partial E,N}\left(W\left(M\tilde{\boldsymbol{F}}\right)^\star\cdot\hat{\boldsymbol{n}}+W^\star\left(M\tilde{\boldsymbol{F}}\right)\cdot\hat{\boldsymbol{n}}-W\left(M\tilde{\boldsymbol{F}}\right)\cdot\hat{\boldsymbol{n}}\right) \diff S_\xi - \left|\left|\sqrt{M}\vec{F}\right|\right|^2_{J,N}.
\label{eq:SemiDiscrete:Phi-equation}
\end{equation}

Now, \eqref{eq:SemiDiscrete:W-equation} and \eqref{eq:SemiDiscrete:Phi-equation} contain the term $\left\langle \mathcal J W, \Phi_t\right\rangle_{E,N}$ on the left hand side, so we can equate both right hand sides to obtain
\begin{equation}
\begin{split}
  \int_{\partial E,N}\left(W\left(M\tilde{\boldsymbol{F}}\right)^\star\cdot\hat{\boldsymbol{n}}+W^\star\left(M\tilde{\boldsymbol{F}}\right)\cdot\hat{\boldsymbol{n}}-W\left(M\tilde{\boldsymbol{F}}\right)\cdot\hat{\boldsymbol{n}}\right) \diff S_\xi - \left|\left|\sqrt{M}\vec{F}\right|\right|^2_{J,N}    \\
= \frac{\diff }{\diff t}\left\langle \mathcal J \Psi, 1\right\rangle_{E,N} + \frac{k}{2}\frac{\diff }{\diff t}||Q||_{J,N}^2 -k\int_{\partial E,N}\left(\Phi_t\tilde{\boldsymbol{Q}}^\star\cdot\hat{\boldsymbol{n}}+\Phi^\star_t \tilde{\boldsymbol{Q}}\cdot\hat{\boldsymbol{n}}-\Phi_t \tilde{\boldsymbol{Q}}\cdot\hat{\boldsymbol{n}}\right)\diff S_\xi.
\end{split}
\label{eq:UnarrangedEnergyEqn}
\end{equation}
Rearranging \eqref{eq:UnarrangedEnergyEqn} to move time derivatives to the left hand side of the equation, we get
\begin{equation}
\frac{\diff }{\diff t}\left(\left\langle \mathcal J \Psi, 1\right\rangle_{E,N} + \frac{k}{2}\left|\left| \vec{Q}\right|\right|_{J,N}^2 \right) = - \left|\left|\sqrt{M}\vec{F}\right|\right|^2_{J,N} + k\text{BT}_{E,N}(\Phi_t,\vec{Q}) + \text{BT}_{E,N}(W,M\vec{F}).
\label{eq:SemiDiscrete:Stability-one-element}
\end{equation}
In \eqref{eq:SemiDiscrete:Stability-one-element} we used the boundary operator defined as (in both physical and computational spaces according to \eqref{eq:DG:local-physical-fluxes-relationship}),
\begin{equation}
\text{BT}_{E,N}\left(\vartheta,\vec{\tau}\right) = \int_{\partial E, N}\left(\vartheta \tilde{\boldsymbol{{\tau}}}^\star+\left(\vartheta^\star-\vartheta\right) \tilde{\boldsymbol{\tau}} \right)\cdot\hat{\boldsymbol{n}}\diff S_\xi = \int_{\partial e, N}\left(\vartheta \vec{\tau}^\star+\left(\vartheta^\star-\vartheta\right) \vec{\tau} \right)\cdot\vec{n}\diff S.
\label{eq:BoundaryOperatorDefinition}
\end{equation}
Additionally, we identify the volumetric 
discrete free--energy of the element,
\begin{equation}
\frac{\diff }{\diff t}\left(\left\langle \mathcal J \Psi, 1\right\rangle_{E,N} + \frac{k}{2}\left|\left| \vec{Q}\right|\right|_{J,N}^2 \right)   
= \frac{\diff }{\diff t}\int_{E,N}\mathcal J\left(\Psi + \frac{k}{2}\vec{Q}\cdot\vec{Q}  
\right)dE = \frac{\diff \mathcal F_{v}^{E,N}}{\diff t},
\end{equation}
and thus, we simplify \eqref{eq:SemiDiscrete:Stability-one-element} to 
\begin{equation}
\frac{\diff \mathcal F_{v}^{E,N}}{\diff t} = - \left|\left|\sqrt{M}\vec{F}\right|\right|^2_{J,N}
+ k\text{BT}_{E,N}(\Phi_t,\vec{Q}) + \text{BT}_{E,N}(W,M\vec{F}).
\label{eq:SemiDiscrete:Stability-one-element-free-en}
\end{equation}
In \eqref{eq:SemiDiscrete:Stability-one-element-free-en}, we find that the 
volumetric free--energy is dissipated in the element interior by the 
chemical potential flux (similarly to the continuous counterpart 
\eqref{eq:continuous:well-posedness-timeder}), and exchanged with other elements through 
the boundary terms $\text{BT}_{E,N}$. 

To obtain an energy estimate similar to that in 
\eqref{eq:continuous:well-posedness}, we sum all element contributions, getting
  \begin{equation}
\sum_{e}\frac{\diff \mathcal F_{v}^{E,N}}{\diff t} = \frac{\diff \mathcal F_{v}^{N}}{\diff t} = -\sum_{e} \left|\left|\sqrt{M}\vec{F}\right|\right|^2_{J,N} + 
\sum_{e}\left(k\text{BT}_{E,N}(\Phi_t,\vec{Q}) + \text{BT}_{E,N}(W,M\vec{F})\right).
\label{eq:SemiDiscrete:free-energy-balance-all-elements}
\end{equation}
We then split the boundary quadratures $\sum_{e}\text{BT}_{E,N} = \text{IBT}_{N}+\text{PBT}_{N}$  into the combination of interior ($\text{IBT}_N$) and 
physical boundary ($\text{PBT}_N$) sums. We first transform the sums to physical coordinates using \eqref{eq:DG:local-physical-fluxes-relationship}, 
\begin{equation}
  \begin{split}
k\text{BT}_{E,N}&(\Phi_t,\vec{Q}) + \text{BT}_{E,N}(W,M\vec{F})  \\
= k&\int_{\partial E,N}\left(\Phi_t\tilde{\boldsymbol{Q}}^\star\cdot\hat{\boldsymbol{n}}+\Phi^\star_t \tilde{\boldsymbol{Q}}\cdot\hat{\boldsymbol{n}}-\Phi_t \tilde{\boldsymbol{Q}}\cdot\hat{\boldsymbol{n}}\right)\diff S_\xi \\
+&\int_{\partial E,N}\left(W\left(M\tilde{\boldsymbol{F}}\right)^\star\cdot\hat{\boldsymbol{n}}+W^\star\left(M\tilde{\boldsymbol{F}}\right)\cdot\hat{\boldsymbol{n}}-W\left(M\tilde{\boldsymbol{F}}\right)\cdot\hat{\boldsymbol{n}}\right) 
\diff S_\xi \\
= k&\int_{\partial e,N}\left(\Phi_t\vec{Q}^\star\cdot\vec{n}+\Phi^\star_t \vec{Q}\cdot\vec{n}-\Phi_t \vec{Q}\cdot\vec{n}\right)\diff S \\
 +&\int_{\partial e,N}\left(W\left(M\vec{F}\right)^\star\cdot\vec{n}+W^\star\left(M\vec{F}\right)\cdot\vec{n}-W\left(M\vec{F}\right)\cdot\vec{n}\right) 
\diff S.
\end{split}
\end{equation}
At interior faces there is a contribution from the left and the right sides. We account for the contribution of the two neighbouring 
elements in the following way: if $E^{+}$ and $E^{-}$ are two elements that share a face $f$, the sum of both contributions to the face is
\begin{equation}
\begin{split}
\int_{f, N}\left(\tilde{\vec{F}}\cdot\hat{\vec{n}}\biggr|^{\partial e^{+}} + \tilde{\vec{F}}\cdot\hat{\vec{n}}\biggr|^{\partial e^{-}}\right)\diff S_\xi &= \int_{f,N}\left(\vec{F}\cdot\vec{n}\biggr|^{\partial e^{+}}+\vec{F}\cdot\vec{n}\biggr|^{\partial e^{-}}\right)\diff S \\
&=\int_{f,N}\llbracket \vec{F}\rrbracket \diff S.
\end{split}
\label{eq:DG:two-surface-elements-combination}
\end{equation}

Using the definitions for the jumps at an interface between two elements, the face contributions to the sum over all elements of the boundary terms can be written as
\begin{equation}
  \begin{split}
  \text{IBT}_N = k&\sum_{\interiorfaces}\int_{f,N}\left(\jump{\Phi_t\vec{Q}^\star}+\jump{\Phi^\star_t \vec{Q}}-\jump{\Phi_t \vec{Q}}\right)\diff S \\
+&\sum_{\interiorfaces}\int_{f,N}\left(\jump{W\left(M\vec{F}\right)^\star}+\jump{W^\star M\vec{F}}-\jump{W M\vec{F}}\right) 
\diff S.
\end{split}
\end{equation}
Since interface numerical fluxes are uniquely defined at the interfaces, as noted in \eqref{eq:DG:BR1-div-riemann}, we can 
remove them from the jump operator so
\begin{equation}
  \begin{split}
  \text{IBT}_N = k&\sum_{\interiorfaces}\int_{f,N}\left(\vec{Q}^\star\cdot\jump{\Phi_t}+\Phi^\star_t \jump{\vec{Q}}-\jump{\Phi_t \vec{Q}}\right)\diff S \\
 +&\sum_{\interiorfaces}\int_{f,N}\left(\left(M\vec{F}\right)^\star\cdot\jump{W}+W^\star\jump{ M\vec{F}}-\jump{W M\vec{F}}\right) 
\diff S.
\end{split}
\label{eq:SemiDiscrete:IBT-with-numerical-fluxes}
\end{equation}
To show stability, we replace the numerical fluxes in \eqref{eq:SemiDiscrete:IBT-with-numerical-fluxes} with the averages of the BR1 method 
 (see \eqref{eq:DG:BR1+IP-solution-riemann} and \eqref{eq:DG:BR1-div-riemann}),
\begin{equation}
  \begin{split}
  \text{IBT}_N = k&\sum_{\interiorfaces}\int_{f,N}\left(\average{\vec{Q}}\cdot\jump{\Phi_t}+\average{\Phi_t} \jump{\vec{Q}}-\jump{\Phi_t \vec{Q}}\right)\diff S \\
 +&\sum_{\interiorfaces}\int_{f,N}\left(\average{M\vec{F}}\cdot\jump{W}+\average{W}\jump{ M\vec{F}}-\jump{W M\vec{F}}\right) 
\diff S \\
-&\sum_{\interiorfaces}\int_{f,N}\sigma M\jump{W}^2\diff S = -\sum_{\interiorfaces}\int_{f,N}\sigma M\jump{W}^2\diff 
S.
\end{split}
\label{eq:SemiDiscrete:IBT-finale}
\end{equation}

The two first terms in \eqref{eq:SemiDiscrete:IBT-finale} vanish by the algebraic identity  \eqref{eq:dg:average-jump-rel}, and the third is always negative. As a result, the contribution from interior faces only decreases the free--energy (or does not contribute if $\sigma=0$). 

For the physical boundary terms, 
\begin{equation}
  \begin{split}
  \text{PBT}_N = k&\sum_{\boundaryfaces}\int_{f,N}\left(\Phi_t\tilde{\boldsymbol{Q}}^\star\cdot\hat{\boldsymbol{n}}+\Phi^\star_t \tilde{\boldsymbol{Q}}\cdot\hat{\boldsymbol{n}}-\Phi_t \tilde{\boldsymbol{Q}}\cdot\hat{\boldsymbol{n}}\right)\diff S_\xi \\
+&\sum_{\boundaryfaces}\int_{f,N}\left(W\left(M\tilde{\boldsymbol{F}}\right)^\star\cdot\hat{\boldsymbol{n}}+W^\star\left(M\tilde{\boldsymbol{F}}\right)\cdot\hat{\boldsymbol{n}}-W\left(M\tilde{\boldsymbol{F}}\right)\cdot\hat{\boldsymbol{n}}\right) 
\diff S_\xi,
\end{split}
\end{equation}
we set the values specified in \eqref{eq:DG:Neumann-bcs-wphi} and 
\eqref{eq:DG:Neumann-bcs-fq}, so that the second quadrature vanishes.  
In the first quadrature we can use the chain rule in time
to see that
\begin{equation}
  \text{PBT}_N = k\sum_{\boundaryfaces} \int_{f,N}\Phi_t G'(\Phi)\diff S  = \frac{\diff }{\diff t}\left(k\sum_{\boundaryfaces} \int_{f,N} 
  G(\Phi)\diff S\right) = -\frac{\diff \mathcal F^{N}_s}{\diff t}.
  \label{eq:SemiDiscrete:PBT}
\end{equation}
Hence, we can rewrite \eqref{eq:SemiDiscrete:free-energy-balance-all-elements}  
to mimic the continuous free--energy bound \eqref{eq:continuous:well-posedness-timeder},
\begin{equation}
\begin{split}
\frac{\diff \mathcal F_{v}^{N}}{\diff t} +\frac{\diff \mathcal F_{s}^{N}}{\diff t} = \frac{\diff \mathcal F^{N}}{\diff t} &= -\sum_{e} \left|\left|\sqrt{M}\vec{F}\right|\right|^2_{J,N}-\sum_{\interiorfaces}\int_{f,N}\sigma M\jump{W}^2\diff S \\
&\leqslant -\sum_{e} \left|\left|\sqrt{M}\vec{F}\right|\right|^2_{J,N}\leqslant 0.
\end{split}
\label{eq:SemiDiscrete:well-posedness-timder}
\end{equation}
When integrated in time, we obtain the semi--discrete version of 
\eqref{eq:continuous:well-posedness},
\begin{equation}
\mathcal F^{N}(T) \leqslant  \mathcal F^{N}(0)-\int_{0}^T\left(\sum_{e} \left|\left|\sqrt{M}\vec{F}\right|\right|^2_{J,N}\right)dt \leqslant   \mathcal F^{N}(0),
\label{eq:SemiDiscrete:well-posedness}
\end{equation}
where we recall that the numerical errors incurred at the interior boundaries $\text{IBT}_{N}$
are zero for the BR1 scheme ($\sigma=0$), and dissipative otherwise ($\sigma>0$). 

Therefore, the 
semi--discrete DG scheme \eqref{eq:DG:weakforms-with-surface-phi}-\eqref{eq:DG:weakforms-with-surface-w} and \eqref{eq:SemiDiscrete:weakforms-with-surface-q} is 
stable, in the sense that the discrete free--energy,
\begin{equation}
  \mathcal F^{N} = \sum_{e}\left\{\int_{E,N}\mathcal J\left(\Psi + 
  \frac{1}{2}k\vec{Q}\cdot\vec{Q}\right)dE - \int_{\partial e\bigcap 
  \partial\Omega,N}\!\!\!\!\!\!\!\!\!\!\!\!\!\!\!\!\!\!\!\! G(\Phi)\diff S\right\},
\end{equation}
is bounded in time by its initial value.

\subsection{Fully--discrete stability analysis using an IMEX time integrator}\label{subsec:FullyDiscreteIMEX}

Because of the stiffness due to the high order derivatives in the Cahn-Hilliard equation, we use an IMplicit--EXplicit (IMEX) time integrator. We consider equally--spaced time steps $\Delta t$ 
and we use the superscript $n$ to represent state values in $t_n=n\Delta t$, 
\begin{equation}
\vartheta^n =\vartheta(t_n) = \vartheta(n\Delta t).
\end{equation}
In particular, this IMplicit--EXplicit (IMEX) solver combines a forward and 
backward Euler scheme. We split the equations into implicit and explicit parts to make the algebraic system that needs to be solved linear, despite the original equation being non--linear, 
 so that its solution can be found with a fast, efficient direct linear solver. 

We recall the spatially continuous PDE to emphasize the time discretization now
\begin{equation}
\phi_t =\nabla\cdot\left(M\nabla \mu\right)= \nabla\cdot\left(M\nabla\left(\psi(\phi) 
-k\nabla^2\phi\right)\right).
\end{equation}
Integrating over one time--step $t\in[t_n,t_{n+1}]$, we find,
\begin{equation}
\int_{t_n}^{t_{n+1}}\phi_tdt = \phi^{n+1} - \phi^{n} =\int_{t_n}^{t_{n+1}}\nabla\cdot\left(M\nabla\left(\psi(\phi) 
-k\nabla^2\phi\right)\right)dt.
\label{eq:FullyDiscrete:PDE-integrated}
\end{equation}
Depending on whether we evaluate the right hand side of \eqref{eq:FullyDiscrete:PDE-integrated} 
at $t_n$ or $t_{n+1}$ we obtain forward (explicit) and backward (implicit) Euler schemes respectively. 

On the one hand, the interface energy term $k\nabla^2 \phi$ needs to be implicit, since otherwise 
explicit time integration requires an impractical time--step restriction. On the other hand, the chemical free--energy term $\psi(\phi)$ 
is nonlinear and is more easily treated explicitly. 
For these reasons, we evaluate \eqref{eq:FullyDiscrete:PDE-integrated} in time as
\begin{equation}
\phi^{n+1} - \phi^{n} =\Delta t\nabla\cdot\left(M\nabla\left(\psi(\phi^{n}) 
-k\nabla^2\phi^{n+1}\right)\right).
\label{eq:FullyDiscrete:PDE-imex-1}
\end{equation}

The time integration scheme \eqref{eq:FullyDiscrete:PDE-imex-1} is not provably stable 
because of the non--linearities in the chemical free--energy. Following \cite{2018:Dong}, we add 
numerical stabilization proportional to solution time jumps 
$\Delta\Phi=\Phi^{n+1}-\Phi^{n}$,
\begin{equation}
\phi^{n+1} - \phi^{n} =\Delta t\nabla\cdot\left(M\nabla\left(\psi(\phi^{n}) 
+S_0\left(\phi^{n+1} -\phi^{n}\right)-k\nabla^2\phi^{n+1}\right)\right),
\label{eq:FullyDiscrete:PDE-imex}
\end{equation}
which maintains first order accuracy, and makes it possible to obtain a stable scheme, as will be seen later. 

We then introduce the time discretization \eqref{eq:FullyDiscrete:PDE-imex} into \eqref{eq:DG:weakforms-with-surface} 
to obtain the fully--discrete discontinuous Galerkin approximation,
\begin{subequations}\label{eq:FullyDiscrete:weakforms-with-surface}
\begin{align}
\left\langle \mathcal J\frac{\Phi^{n+1}-\Phi^{n}}{\Delta t},\varphi_\phi\right\rangle_{E,N} &= \int_{\partial E,N}\varphi_\phi\left(M\tilde{\boldsymbol{F}}\right)^{\star,\theta}\cdot\hat{\boldsymbol{n}} \diff S_\xi - \left\langle M\tilde{\boldsymbol{F}}^{\theta}, \nabla_\xi\varphi_\phi\right\rangle_{E,N}, \label{eq:FullyDiscrete:weakforms-with-surface-phi}\\
\left\langle\mathcal J\vec{F}^{\theta},\vec{\varphi}_F \right\rangle_{E,N} &=  \int_{\partial E,N}\left(W^{\star,\theta}-W^{\theta}\right) \tilde{\boldsymbol{\varphi}}_F\cdot\hat{\boldsymbol{n}}\diff S_\xi + \left\langle \nabla_\xi W^{\theta},\tilde{\boldsymbol{\varphi}}_F\right\rangle_{E,N},\label{eq:FullyDiscrete:weakforms-with-surface-f}\\
\left\langle \mathcal J W^{\theta},\varphi_W\right\rangle_{E,N} &= \left\langle  \left(\frac{\diff \Psi}{\diff \Phi}\right)^{n}+S_0\left(\Phi^{n+1}-\Phi^{n}\right) ,\mathcal J\varphi_W\right\rangle_{E,N} \notag \\
&-k\int_{\partial E,N}\varphi_W\tilde{\boldsymbol{Q}}^{\star,n+1}\cdot\hat{\boldsymbol{n}}\diff S_\xi + k\left\langle\tilde{\boldsymbol{Q}}^{n+1},\nabla_\xi\varphi_W\right\rangle_{E,N}, \label{eq:FullyDiscrete:weakforms-with-surface-w}\\
\left\langle \mathcal J\vec{Q}^{n+1},\vec{\varphi}_Q\right\rangle_{E,N} &= \int_{\partial E,N}\left(\Phi^{\star,n+1}-\Phi^{n+1}\right) \tilde{\boldsymbol{\varphi}}_Q\cdot\hat{\boldsymbol{n}}\diff S_\xi +  \left\langle 
\nabla_\xi\Phi^{n+1},\tilde{\boldsymbol{\varphi}}_Q\right\rangle_{E,N} \label{eq:FullyDiscrete:weakforms-with-surface-q},
\end{align}
\end{subequations}
where we use the superscript $\theta$ for variables (e.g. $\vec{F}^{\theta}$ or $W^{\theta}$) 
that are not directly evaluated at $t_n$ or $t_{n+1}$ with the IMEX strategy, but on a combination of those depending on the different terms involved in \eqref{eq:FullyDiscrete:weakforms-with-surface-w}.

To analyze the stability of the system \eqref{eq:FullyDiscrete:weakforms-with-surface}, we start by combining \eqref{eq:FullyDiscrete:weakforms-with-surface-w} and \eqref{eq:FullyDiscrete:weakforms-with-surface-q}. 
We perform the first manipulations on \eqref{eq:FullyDiscrete:weakforms-with-surface-q}, which we set for both $t_n$ and $t_{n+1}$, 
\begin{subequations} \label{eq:FullyDiscrete:Q-weak}
\begin{align}
\left\langle \mathcal J\vec{Q}^{n+1},\vec{\varphi}_Q\right\rangle_{E,N} &= \int_{\partial E,N}\left(\Phi^{\star,n+1}-\Phi^{n+1}\right) \tilde{\boldsymbol{\varphi}}_Q\cdot\hat{\boldsymbol{n}}\diff S_\xi +  \left\langle 
\nabla_\xi\Phi^{n+1},\tilde{\boldsymbol{\varphi}}_Q\right\rangle_{E,N} \label{eq:FullyDiscrete:Q-weak-n+1},\\
\left\langle \mathcal J\vec{Q}^{n},\vec{\varphi}_Q\right\rangle_{E,N} &= \int_{\partial E,N}\left(\Phi^{\star,n}-\Phi^{n}\right) \tilde{\boldsymbol{\varphi}}_Q\cdot\hat{\boldsymbol{n}}\diff S_\xi +  \left\langle 
\nabla_\xi\Phi^{n},\tilde{\boldsymbol{\varphi}}_Q\right\rangle_{E,N}. \label{eq:FullyDiscrete:Q-weak-n}
\end{align}
\end{subequations}
Then we subtract \eqref{eq:FullyDiscrete:Q-weak-n} from \eqref{eq:FullyDiscrete:Q-weak-n+1}, divide the result by $\Delta t$ (note that we have defined $\Delta \Phi = \Phi^{n+1}-\Phi^{n}$ and $\Delta \vec{Q}=\vec{Q}^{n+1}-\vec{Q}^{n}$),
\begin{equation}
  \left\langle \mathcal J\frac{\Delta\vec{Q}}{\Delta t},\vec{\varphi}_Q\right\rangle_{E,N} = \int_{\partial E,N}\left(\frac{\Delta\Phi^{\star}}{\Delta t}-\frac{\Delta\Phi}{\Delta t}\right) \tilde{\boldsymbol{\varphi}}_Q\cdot\hat{\boldsymbol{n}}\diff S_\xi +  \left\langle \frac{\nabla_\xi\left(\Delta\Phi\right)}{\Delta t},\tilde{\boldsymbol{\varphi}}_Q\right\rangle_{E,N}, \label{eq:FullyDiscrete:Q-diff-weak}
\end{equation}
and we set $\vec{\varphi}_Q = \vec{Q}^{n+1}$ in \eqref{eq:FullyDiscrete:Q-diff-weak} to obtain
\begin{equation}
  \left\langle \mathcal J\frac{\Delta\vec{Q}}{\Delta t}, \vec{Q}^{n+1}\right\rangle_{E,N} = \int_{\partial E,N}\left(\frac{\Delta\Phi^{\star}}{\Delta t}-\frac{\Delta\Phi}{\Delta t}\right)  \tilde{\boldsymbol{Q}}^{n+1}\cdot\hat{\boldsymbol{n}}\diff S_\xi +  \left\langle \frac{\nabla_\xi\left(\Delta\Phi\right)}{\Delta t}, \tilde{\boldsymbol{Q}}^{n+1}\right\rangle_{E,N} \label{eq:FullyDiscrete:Q-diff-stab}.
\end{equation}
Next, we set $\varphi_W = \Delta\Phi/\Delta t = \left(\Phi^{n+1}-\Phi^{n}\right)/\Delta t$ in 
\eqref{eq:FullyDiscrete:weakforms-with-surface-w},
\begin{equation}
  \begin{split}
\left\langle \mathcal J W^{\theta},\frac{\Delta\Phi}{\Delta t}\right\rangle_{E,N} =& \left\langle \left(\frac{\diff \Psi}{\diff \Phi}\right)^{n}+S_0\Delta\Phi ,\mathcal J \frac{\Delta\Phi}{\Delta t}\right\rangle_{E,N} \\
&-k\int_{\partial E,N}\frac{\Delta\Phi}{\Delta t}\tilde{\boldsymbol{Q}}^{\star,n+1}\cdot\hat{\boldsymbol{n}}\diff S_\xi + k\left\langle\tilde{\boldsymbol{Q}}^{n+1},\frac{\nabla_\xi\left(\Delta\Phi\right)}{\Delta t}\right\rangle_{E,N},
 \label{eq:FullyDiscrete:W-equation}
\end{split}
\end{equation}
and replace the last inner product in \eqref{eq:FullyDiscrete:W-equation} by that in \eqref{eq:FullyDiscrete:Q-diff-stab},
\begin{equation}
\begin{split}
\left\langle \mathcal J W^{\theta},\frac{\Delta\Phi}{\Delta t}\right\rangle_{E,N} =& \left\langle \left(\frac{\diff \Psi}{\diff \Phi}\right)^{n}+S_0\Delta\Phi , \mathcal J\frac{\Delta\Phi}{\Delta t}\right\rangle_{E,N}+k  \left\langle \mathcal J\frac{\Delta\vec{Q}}{\Delta t}, \vec{Q}^{n+1}\right\rangle_{E,N}  \\
&-k\int_{\partial E,N}\left(\frac{\Delta\Phi}{\Delta t}\tilde{\boldsymbol{Q}}^{\star,n+1}\cdot\hat{\boldsymbol{n}}+\left(\frac{\Delta\Phi^{\star}-\Delta\Phi}{\Delta t}\right) \tilde{\boldsymbol{Q}}^{n+1}\cdot\hat{\boldsymbol{n}}\right) \diff S_\xi.
\end{split}
\label{eq:FullyDiscrete:W+Q-equation}
\end{equation}
Using the boundary operator \eqref{eq:BoundaryOperatorDefinition},
\begin{equation}
\begin{split}
\left\langle \mathcal J W^{\theta},\frac{\Delta\Phi}{\Delta t}\right\rangle_{E,N} &= \left\langle\left(\frac{\diff \Psi}{\diff \Phi}\right)^{n}+S_0\Delta\Phi, \mathcal J \frac{\Delta\Phi}{\Delta t}\right\rangle_{E,N}+k  \left\langle \mathcal J\frac{\Delta\vec{Q}}{\Delta t}, \vec{Q}^{n+1}\right\rangle_{E,N}  \\
&-k\text{BT}_{E,N}\left(\frac{\Delta\Phi}{\Delta t}, \vec{Q}^{n+1}\right).
\end{split}
\label{eq:FullyDiscrete:W+Q-equation-simpler}
\end{equation}

Next, we combine \eqref{eq:FullyDiscrete:weakforms-with-surface-phi} and \eqref{eq:FullyDiscrete:weakforms-with-surface-f}. To do so, we set $\vec{\varphi}_{\vec{F}} = \mathcal I^{N}\left[M\vec{F}^{\theta}\right]$ in \eqref{eq:FullyDiscrete:weakforms-with-surface-f} (note we drop the $\mathcal I^{N}$ operator since the quadrature only requires nodal values, that is $\left\langle\mathcal I^{N}\left(M\vec{F}^{\theta}\right),\vartheta\right\rangle_{E,N}=\left\langle M\vec{F}^{\theta},\vartheta\right\rangle_{E,N}$),
\begin{equation}
\left\langle\mathcal J\vec{F}^{\theta},M\vec{F}^{\theta} \right\rangle_{E,N} =  \int_{\partial E,N}\left(W^{\star,\theta}-W^{\theta}\right) M\tilde{\boldsymbol{F}}^{\theta}\cdot\hat{\boldsymbol{n}}\diff S_\xi + \left\langle \nabla_\xi W^{\theta},M\tilde{\boldsymbol{F}}^{\theta}\right\rangle_{E,N},
\label{eq:FullyDiscrete:F-equation}
\end{equation}
and we set $\varphi_\phi=W^{\theta}$ in \eqref{eq:FullyDiscrete:weakforms-with-surface-phi},
\begin{equation}
\left\langle \mathcal J\frac{\Delta\Phi}{\Delta t},W^{\theta}\right\rangle_{E,N} = \int_{\partial E,N}W^{\theta}\left(M\tilde{\boldsymbol{F}}\right)^{\star,\theta}\cdot\hat{\boldsymbol{n}} \diff S_\xi - \left\langle M\tilde{\boldsymbol{F}}^{\theta}, \nabla_{\xi} 
W^{\theta}\right\rangle_{E,N}.
\label{eq:FullyDiscrete:Phi-equation}
\end{equation}
We then sum \eqref{eq:FullyDiscrete:F-equation} and \eqref{eq:FullyDiscrete:Phi-equation}, and use the boundary term operator \eqref{eq:BoundaryOperatorDefinition},
\begin{equation}
\left\langle \mathcal J\frac{\Delta\Phi}{\Delta t},W^{\theta}\right\rangle_{E,N} = - \left\langle\mathcal J\vec{F}^{\theta},M\vec{F}^{\theta} \right\rangle_{E,N} + \text{BT}_{E,N} \left( W^{\theta},M\vec{F}^{\theta} \right).
\label{eq:FullyDiscrete:Phi+F-equation}
\end{equation}
The final step is to combine \eqref{eq:FullyDiscrete:W+Q-equation-simpler} and
\eqref{eq:FullyDiscrete:Phi+F-equation}. Since they share their left hand sides, we equate both right hand sides and multiply them by the time step $\Delta t$,
\begin{equation}
\begin{split}
 &\left\langle  \left(\frac{\diff \Psi}{\diff \Phi}\right)^{n}+S_0\Delta\Phi, \mathcal J\Delta\Phi\right\rangle_{E,N}+k  \left\langle \mathcal J\Delta\vec{Q}, \vec{Q}^{n+1}\right\rangle_{E,N}  \\
&= - \Delta t\left\langle\mathcal J\vec{F}^{\theta},M\vec{F}^{\theta} \right\rangle_{E,N} +  k \text{BT}_{E,N}\left(\Delta\Phi, \vec{Q}^{n+1}\right) + \Delta t \text{BT}_{E,N} \left( W^{\theta},M\vec{F}^{\theta} 
\right).
\label{eq:FullyDiscrete:main-equation-stability}
\end{split}
\end{equation}

We then perform manipulations on the left hand side to get the free--energy $\mathcal F$. First, we perform the Taylor expansion of $\Psi(\Phi)$ centered on $\Phi^{n}$,
\begin{equation}
\Psi^{n+1} = \Psi^{n} + \left(\frac{\diff \Psi}{\diff \Phi}\right)^{n} \Delta \Phi +\frac{1}{2}\left(\frac{\diff ^2\Psi}{\diff \Phi^2}\right)^{n} \Delta \Phi^2 + \frac{1}{6}\left(\frac{\diff ^3\Psi}{\diff \Phi^3}\right)^{n} \Delta \Phi^3 +\frac{1}{24}\left(\frac{\diff ^4\Psi}{\diff \Phi^4}\right)^{n} \Delta \Phi^4 + ...
\end{equation}
For the polynomic chemical energy, \eqref{eq:continuous:chemical-free-energy}, 
\begin{equation}\left(d\Psi/d\Phi\right)^n = -\Phi^n + \left(\Phi^{n}\right)^3,\end{equation}
 \begin{equation}\left(d^2\Psi/d\Phi^2\right)^n = -1 + 3\left(\Phi^{n}\right)^2,\end{equation} 
 \begin{equation}
 \left(d^3\Psi/d\Phi^3\right)^n = 6\Phi^n,\end{equation} and 
 \begin{equation}\left(d^4\Psi/d\Phi^4\right)^n = 6,
 \end{equation} so, it follows exactly that
\begin{equation}
\Psi^{n+1} = \Psi^{n} + \left(\frac{\diff \Psi}{\diff \Phi}\right)^{n} \Delta \Phi -\frac{1}{2}\left( 1 -3 \left(\Phi^{n}\right)^2\right) \Delta \Phi^2 + \Phi^n \Delta \Phi^3 +\frac{1}{4}\Delta \Phi^4.
\label{eq:FullyDiscrete:chemical-free-energy-Taylor}
\end{equation}
We use \eqref{eq:FullyDiscrete:chemical-free-energy-Taylor} to write the first volume quadrature in \eqref{eq:FullyDiscrete:main-equation-stability} as
\begin{equation}
\left\langle  \left(\frac{\diff \Psi}{\diff \Phi}\right)^{n}+S_0\Delta\Phi ,\mathcal J\Delta\Phi\right\rangle_{E,N} = \left\langle \mathcal J \Psi^{n+1}, 1\right\rangle_{E,N} -\left\langle \mathcal J\Psi^{n}, 1\right\rangle_{E,N}+ \left\langle \mathcal J \Pi, 1\right\rangle_{E,N},
\label{eq:FullyDiscrete:chemical-free-energy-term}
\end{equation}
where $\Pi \left(\Phi^{n+1}, \Phi^{n}\right)$ is the polynomial function
\begin{equation}
\Pi=S_0 \Delta \Phi^2  +\frac{1}{2}\left( 1 - 3\left(\Phi^{n}\right)^2\right) \Delta \Phi^2 - \Phi^n \Delta \Phi^3 - \frac{1}{4}\Delta \Phi^4 = \Delta \Phi^2\Pi^\star,
\end{equation}
with
\begin{equation}
\Pi^\star = S_0+\frac{1}{2}\left( 1 - 3\left(\Phi^{n}\right)^2\right)- \Phi^n \left(\Phi^{n+1}-\Phi^{n}\right)-\frac{1}{4}\left(\Phi^{n+1}-\Phi^{n}\right)^2.
\end{equation}

The quantity $\Pi^{\star}$ is a second order concave polynomial (elliptic paraboloid), which must remain positive for the time integration to be stable, as will be shown from the stability analysis. Since the function is concave, it will always be negative for sufficiently large values of $\Phi^{n+1}$ or $\Phi^{n}$. However, the solutions of the Cahn--Hilliard equation must remain close to the range $\phi\in[-1,1]$. Thus, we can choose the value $S_0$ so that $\Pi^\star$ remains positive in a reasonably large neighborhood of $\left(\Phi^{n+1},\Phi^{n}\right)\in[-1,1]^2$. 

\begin{figure}
\centering
\subfigure[$\Pi^\star=0$ isolines]{\includegraphics[height=5.5cm]{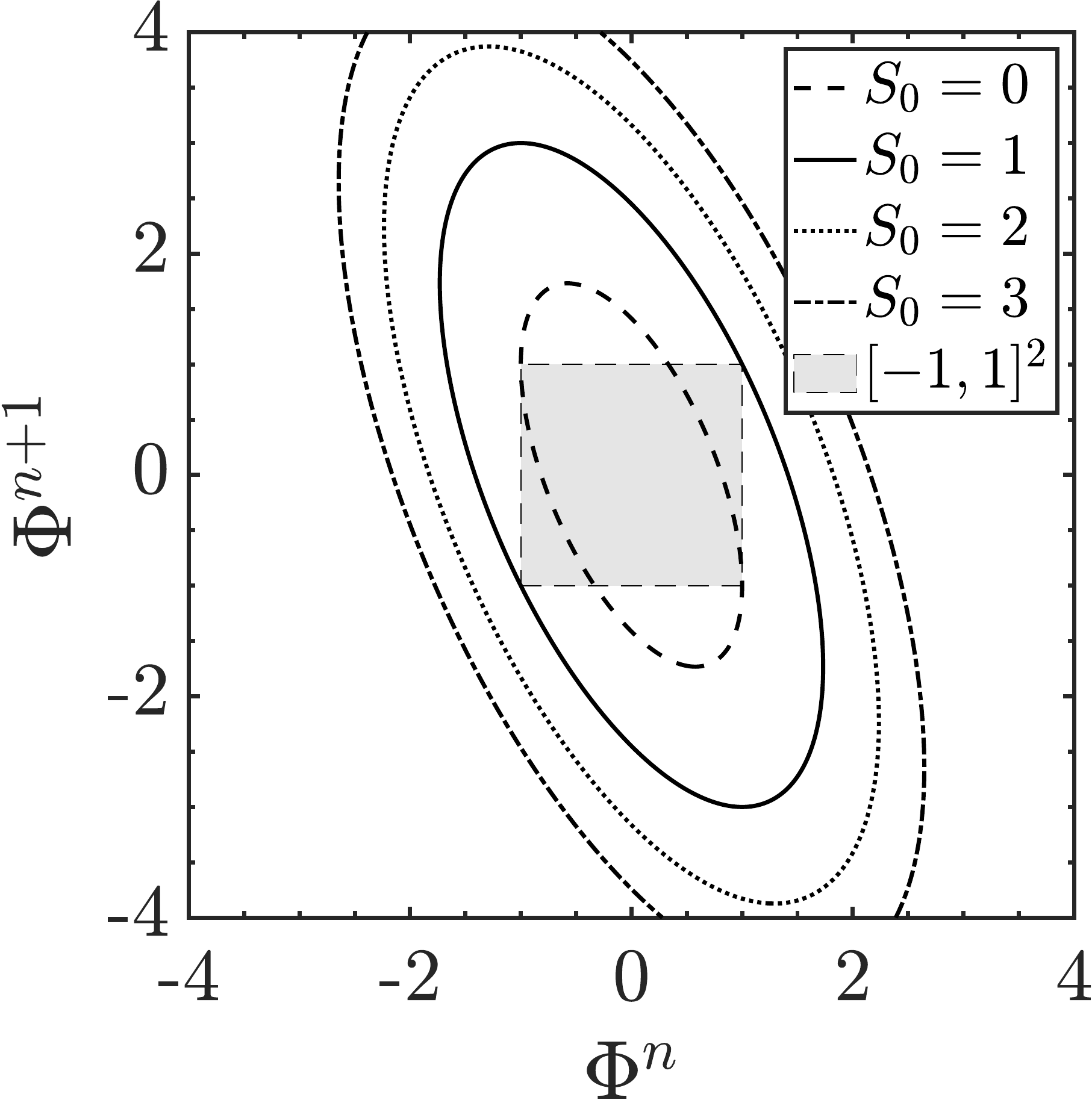}}
\subfigure[Range $\phi\in{\left[-\phi_s,\phi_s\right]^2}$ in which non--linear terms are stable]{\includegraphics[height=5.5cm]{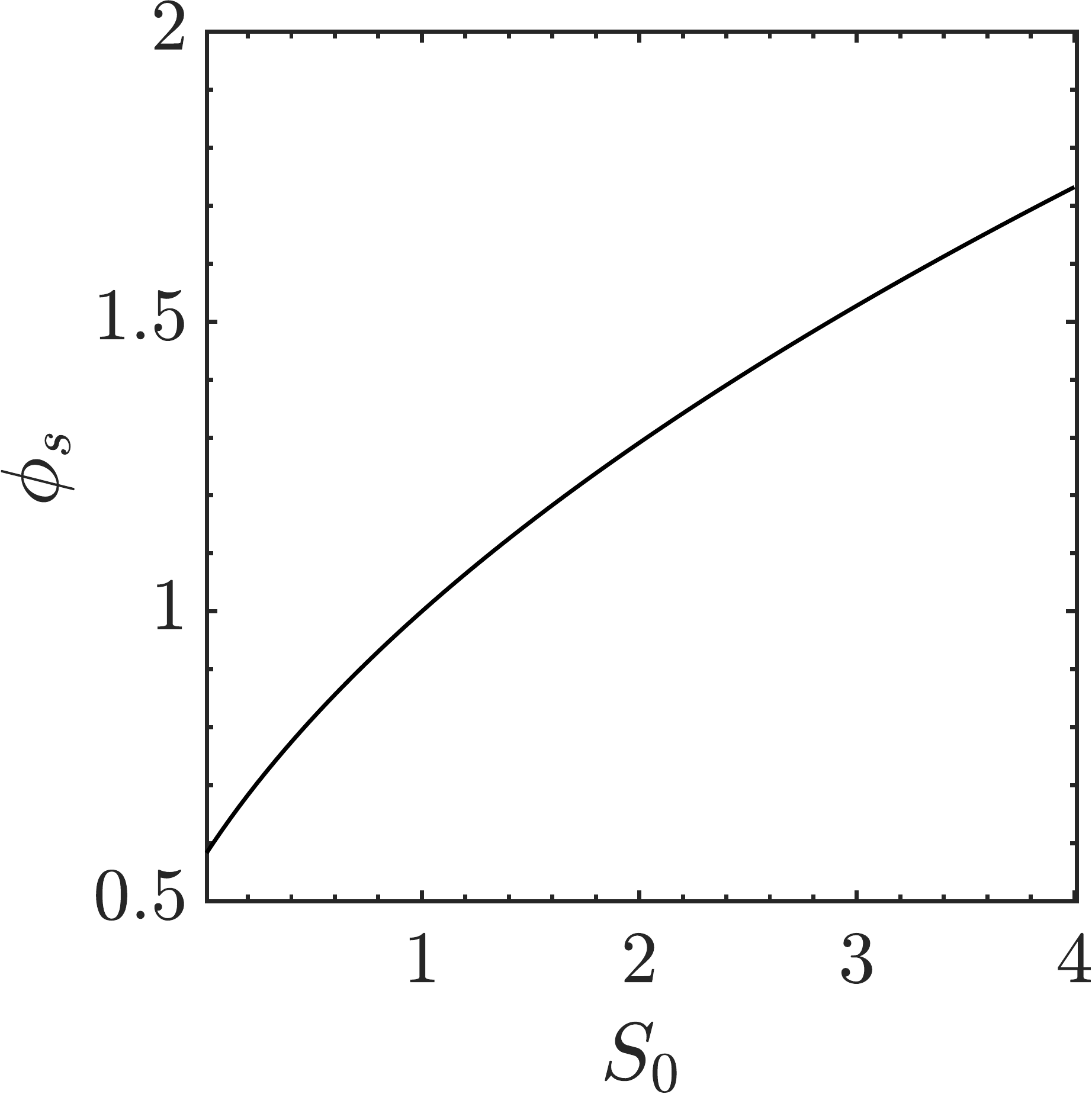}\label{fig:PsiMargin}}
\caption{Graphical representation of the elliptical isoline $\Pi^{\star} = 0$ for several $S_0$ values (indicated by the contour labels). The box $[-1,1]^2$ has also been represented. The approximation is stable if the pair $\left(\Phi^{n},\Phi^{n+1}\right)$ stays inside the ellipse for a given $S_0$ value. In \ref{fig:PsiMargin} we have represented for each $S_0$ the interval $\left[-\phi_s,\phi_s\right]^2$ in which the non--linear terms are stable, which corresponds to the formula $\phi_s = \sqrt{\frac{2S_0+1}{3}}$}
\label{fig:pi-isolines}
\end{figure}

The elliptical isocontours $\Pi^\star=0$ are shown in Fig. \ref{fig:pi-isolines} for several $S_0$ values (as labeled on each contour line). For the approximation to be stable, every pair $\left(\Phi^{n},\Phi^{n+1}\right)$ must remain inside the $\Pi^\star=0$ isoline for the given $S_0$ value. For instance, given the range $\Phi\in[-1,1]$, the scheme is stable for $S_0\geqslant 1$, and for $\Phi\in\left[-\sqrt{5/3},\sqrt{5/3}\right]\simeq[-1.29,1.29]$, the scheme remains stable for $S_0\geqslant 2$.

For the interface energy in \eqref{eq:FullyDiscrete:main-equation-stability}, we complete the square,
\begin{equation}
\begin{split}
\left\langle \mathcal J\Delta\vec{Q}, \vec{Q}^{n+1}\right\rangle_{E,N} =&\phantom{{}+{}} \frac{1}{2}\left\langle \mathcal J \vec{Q}^{n+1}, \vec{Q}^{n+1}\right\rangle_{E,N}  -\frac{1}{2}\left\langle \mathcal J \vec{Q}^{n}, \vec{Q}^{n}\right\rangle_{E,N} \\
&  +\frac{1}{2} \left\langle \mathcal J\Delta\vec{Q},\Delta\vec{Q},\right\rangle_{E,N},
\end{split}
\label{eq:FullyDiscrete:interface-energy-term}
\end{equation}
and we place \eqref{eq:FullyDiscrete:chemical-free-energy-term} and \eqref{eq:FullyDiscrete:interface-energy-term} in \eqref{eq:FullyDiscrete:main-equation-stability},
\begin{equation}
\begin{split}
 & \left\langle \mathcal J \left(\Psi^{n+1}+\frac{1}{2}k\vec{Q}^{n+1}\cdot\vec{Q}^{n+1}\right), 1\right\rangle_{E,N} - \left\langle \mathcal J \left(\Psi^{n}+\frac{1}{2}k\vec{Q}^{n}\cdot\vec{Q}^{n}\right), 1\right\rangle_{E,N}\\
&= - \Delta t\left\langle\mathcal J\vec{F}^{\theta},M\vec{F}^{\theta} \right\rangle_{E,N} +  k\text{BT}_{E,N}\left(\Delta\Phi, \vec{Q}^{n+1}\right) \\
&+ \Delta t \text{BT}_{E,N} \left( W^{\theta},M\vec{F}^{\theta} \right)- \left\langle \mathcal J \Pi, 1\right\rangle_{E,N}-\frac{1}{2}k \left\langle \mathcal J\Delta\vec{Q},\Delta\vec{Q}\right\rangle_{E,N}.
\label{eq:FullyDiscrete:main-equation-stability-1}
\end{split}
\end{equation}
We define the discrete volumetric free--energy in an element as,

\begin{equation}
\mathcal F_{v}^{n,E,N} = \left\langle \mathcal J \left(\Psi^{n}+\frac{1}{2}k\vec{Q}^{n}\cdot\vec{Q}^{n}\right), 1\right\rangle_{E,N},
\end{equation}
which simplifies \eqref{eq:FullyDiscrete:main-equation-stability-1} to
\begin{equation}
\begin{split}
 \mathcal F_{v}^{n+1,E,N}-\mathcal F_{v}^{n,E,N} =& - \Delta t\left\langle\mathcal J\vec{F}^{\theta},M\vec{F}^{\theta} \right\rangle_{E,N} +  k \text{BT}_{E,N}\left(\Delta\Phi, \vec{Q}^{n+1}\right) \\
&+ \Delta t \text{BT}_{E,N} \left( W^{\theta},M\vec{F}^{\theta} \right)- \left\langle \mathcal J \Pi, 1\right\rangle_{E,N} \\
&-\frac{1}{2}k\left\langle \mathcal J\Delta\vec{Q},\Delta\vec{Q}\right\rangle_{E,N}.
\label{eq:FullyDiscrete:main-equation-stability-finale}
\end{split}
\end{equation}
Eq. \eqref{eq:FullyDiscrete:main-equation-stability-finale} shows that free--energy changes are due to physical dissipation in the element interior by the term $ \Delta t\left\langle\mathcal J\vec{F}^{\theta},M\vec{F}^{\theta} \right\rangle_{E,N} \geqslant 0$ (the discrete counterpart of that obtained for the continuous analysis \eqref{eq:continuous:well-posedness-timeder}), numerical dissipation as a result of the IMEX scheme,
\begin{equation}
\text{diss}^{E,N}_{\text{IMEX}} = - \left\langle \mathcal J \Pi, 1\right\rangle_{E,N} -\frac{1}{2}k\left\langle \mathcal J\Delta\vec{Q},\Delta\vec{Q}\right\rangle_{E,N} \leqslant 0,
\label{eq:FullyDiscrete:IMEX-dissipation}
\end{equation}
and boundary exchanges through all $\text{BT}_{E,N}$ terms. 

The effect of boundary exchanges can only be studied from the perspective of all elements in the domain. So we sum \eqref{eq:FullyDiscrete:main-equation-stability-finale} for all mesh elements,
\begin{equation}
\begin{split}
 \mathcal F_{v}^{n+1,N}-\mathcal F_{v}^{n,N} =& - \Delta t\sum_{e}\left\langle\mathcal J\vec{F}^{\theta},M\vec{F}^{\theta} \right\rangle_{E,N} + \text{IBT}_N + \text{PBT}_N + \sum_{e} \text{diss}^{E,N}_{\text{IMEX}},
\label{eq:FullyDiscrete:stability-all-elements-1}
\end{split}
\end{equation}
where $\mathcal F_{v}^{n,N} = \sum_{e}\mathcal F_{v}^{n,E,N}$ is the sum of all element volumetric free--energies. 

On the one hand, the $\text{IBT}_{N}$ is the contribution of all interior boundary quadratures,
\begin{equation}
\begin{split}
\text{IBT}_{N} =& \phantom{{}+k{}}\sum_{\interiorfaces}\int_{f,N}\left(M\vec{F}^{\theta,\star}\cdot\jump{W^{\theta}} + W^{\theta,\star}\jump{M\vec{F}^{\theta}} - \jump{M\vec{F}^{\theta}W^{\theta}}  \right)\diff S \\
&+k\sum_{\interiorfaces}\int_{f,N}\left(\vec{Q}^{n+1,\star}\cdot\jump{\Delta \Phi} + \Delta \Phi^{\star}\jump{\vec{Q}^{n+1}}-\jump{\vec{Q}^{n+1}\Delta \Phi}\right)\diff S.
\end{split}
\end{equation}
Introducing the numerical fluxes obtained with the BR1 scheme we find that since \eqref{eq:dg:average-jump-rel} holds, we obtain the same result as in the semi--discrete analysis \eqref{eq:SemiDiscrete:IBT-finale}, 
\begin{equation}
\begin{split}
\text{IBT}_{N} =&\phantom{{}+k{}} \sum_{\interiorfaces}\int_{f,N}\left(\average{M\vec{F}^\theta}\cdot\jump{W^{\theta}} + \average{W^{\theta}}\jump{M\vec{F}^{\theta}} - \jump{M\vec{F}^{\theta}W^{\theta}}  \right)\diff S \\
&+k\sum_{\interiorfaces}\int_{f,N}\left(\average{\vec{Q}^{n+1}}\cdot\jump{\Delta \Phi} + \average{\Delta \Phi}\jump{\vec{Q}^{n+1}}-\jump{\vec{Q}^{n+1}\Delta \Phi}\right)\diff S \\
&-\phantom{{}k{}}\sum_{\interiorfaces}\int_{f,N}\sigma M\jump{W^{\theta}}^2\diff S = -\sum_{\interiorfaces}\int_{f,N}\sigma M\jump{W^{\theta}}^2\diff S.
\end{split}
\end{equation}
Therefore, we conclude that the global contribution from interior boundaries to the free--energy changes is strictly negative ($\sigma > 0$), or vanishes if $\sigma=0$. 

On the other hand, we consider the physical boundary terms, $\text{PBT}_N$,
\begin{equation}
  \begin{split}
  \text{PBT}_N =k& \sum_{\boundaryfaces} \int_{f,N}\left(\Delta \Phi\tilde{\boldsymbol{Q}}^{n+1,\star}\cdot\hat{\boldsymbol{n}}+\Delta \Phi^{\star} \tilde{\boldsymbol{Q}}^{n+1}\cdot\hat{\boldsymbol{n}}-\Delta \Phi \tilde{\boldsymbol{Q}}^{n+1}\cdot\hat{\boldsymbol{n}}\right)\diff S_\xi \\
 +&\sum_{\boundaryfaces}\int_{f,N}\left(W^{\theta}\left(M\tilde{\boldsymbol{F}}\right)^{\theta,\star}\cdot\hat{\boldsymbol{n}}+W^{\theta,\star}\left(M\tilde{\boldsymbol{F}}^{\theta}\right)\cdot\hat{\boldsymbol{n}}-W^{\theta}\left(M\tilde{\boldsymbol{F}}^{\theta}\right)\cdot\hat{\boldsymbol{n}}\right) 
\diff S_\xi.
\end{split}
\end{equation}
Introducing the numerical fluxes \eqref{eq:DG:Neumann-bcs-fq},
\begin{equation}
  \text{PBT}_N = \sum_{\boundaryfaces}\int_{f,N}\left(\Phi^{n+1}-\Phi^{n}\right)\beta \diff S =\sum_{\boundaryfaces}\int_{f,N}\left(G\left(\Phi^{n+1}\right) - G\left(\Phi^{n}\right)\right) \diff S,
\end{equation}
we transform \eqref{eq:FullyDiscrete:stability-all-elements-1} to
\begin{equation}
\begin{split}
 \mathcal F_{v}^{n+1,N}-\mathcal F_{v}^{n,N}-&\int_{\partial e\bigcup 
  \partial\Omega,N}\!\!\!\!\!\!\!\!\!\!\!\!\!\!\!\!\!\!\!\! \left(G(\Phi^{n+1})-G(\Phi^{n})\right)\diff S = \mathcal F^{n+1,N}-\mathcal F^{n,N}\\
  & =- \Delta t\sum_{e}\left\langle\mathcal J\vec{F}^{\theta},M\vec{F}^{\theta} \right\rangle_{E,N}  + \sum_{e} \text{diss}^{E,N}_{\text{IMEX}}.
\label{eq:FullyDiscrete:stability-all-elements-2}
\end{split}
\end{equation}

The last step is to sum over all time steps $n=0,...,T-1$,
\begin{equation}
\begin{split}
 \mathcal F^{T,N}-\mathcal F^{0,N} &=- \Delta t\sum_{e,n}\left\langle\mathcal J\vec{F}^{\theta},M\vec{F}^{\theta} \right\rangle_{E,N}  + \sum_{e,n} \text{diss}^{E,N}_{\text{IMEX}} -\Delta t\sum_{\interiorfaces}\int_{f,N}\sigma M\jump{W^{\theta}}^2\diff S \\
 & \leqslant - \Delta t\sum_{e,n}\left\langle\mathcal J\vec{F}^{\theta},M\vec{F}^{\theta} 
 \right\rangle_{E,N},
 \end{split}
\end{equation}
which shows that the discrete free--energy is bounded by its initial value

\begin{equation}
\mathcal F^{T,N} \leqslant \mathcal F^{0,N} - \Delta t\sum_{E,n}\left\langle\mathcal J\vec{F}^{\theta},M\vec{F}^{\theta} \right\rangle_{E,N}  \leqslant \mathcal F^{0,N} ,
\label{eq:FullyDiscrete:stability-all-timesteps}
\end{equation}
and thus, the scheme \eqref{eq:FullyDiscrete:weakforms-with-surface} is stable 
and satisfies a free--energy bound consistent with the continuous equation 
\eqref{eq:continuous:well-posedness}. Furthermore, spatial under--resolution introduces numerical dissipation proportional to interface jumps in chemical potential, $\sigma\jump{W^{\theta}}^2$, and temporal under--resolution adds dissipation proportional to $\Delta \Phi^2$ and $\Delta\vec{Q}^2$.
Note that \eqref{eq:FullyDiscrete:stability-all-timesteps} induces a time step $\Delta t $ 
restriction, which ensures the free--energy $\mathcal F^{T,N}$ remains always positive. 
In practice, we have not experienced instabilities due to large time steps (for values small enough to achieve accurate solutions).

\section{Numerical experiments}\label{sec:NumericalExperiments}

In this section we show numerical experiments to address the capabilities and robustness of the fully discrete approximation. First, we perform a convergence study using the method of manufactured solutions. The experiment performed here is two dimensional in a Cartesian mesh with straight sides. Second, we compare our results with two--dimensional results available in the literature, using a Cartesian mesh, a distorted mesh, a distorted mesh with curvilinear faces, and a fully unstructured mesh forming a ``T'' domain. Lastly, we explore the spinodal--decomposition in a three dimensional cylindrical geometry.

\subsection{Convergence study}\label{subsec:Numerical:ConvergenceAnalysis}

We now address the convergence of the fully discrete scheme. To do so, we follow \cite{2016:Kastner} and consider a two--dimensional unit square, $[-1,1]^2$, and impose the solution:
\begin{equation}
\phi_0(x,y,t) = \cos\left(\pi\alpha x\right) \cos\left(\pi\alpha y\right)\cos \left(t\right).
\label{eq:Numerical:Convergence:IC}
\end{equation}
For \eqref{eq:Numerical:Convergence:IC} to be solution of \eqref{eq:cahn--hilliard-eqn}, we add a source term to the latter,
\begin{equation}
\frac{\phi_0}{\partial t} = \nabla\cdot\left(M\nabla \left(\frac{\diff \Psi(\phi_0)}{\diff \phi} - \varepsilon^2 \nabla\phi_0^2\right)\right) + q(\vec{x},t),
\end{equation}
whose expression is,

\begin{equation}
\begin{split}
q\left(x,y,t\right)=&\phantom{{}+{}}\,2M\,\pi^4\,\alpha^4\,{\varepsilon}^2\,\cos\left(t\right)\,\cos\left(\pi \,\alpha \,x\right)\,\cos\left(\pi \,\alpha \,y\right)\\
&+3M\,\pi ^2\,\alpha ^2\,{\cos\left(t\right)}^3\,{\cos\left(\pi \,\alpha \,x\right)}^3\,{\cos\left(\pi \,\alpha \,y\right)}^3\\
&-6M\,\pi ^2\,\alpha ^2\,{\cos\left(t\right)}^3\,\cos\left(\pi \,\alpha \,x\right)\,{\cos\left(\pi \,\alpha \,y\right)}^3\,{\sin\left(\pi \,\alpha \,x\right)}^2\\
&-M\pi ^2\,\alpha ^2\,\cos\left(t\right)\,\cos\left(\pi \,\alpha \,x\right)\,\cos\left(\pi \,\alpha \,y\right)\\
&\,+2M\,\pi ^4\,\alpha ^4\,{\varepsilon}^2\,\cos\left(t\right)\,\cos\left(\pi \,\alpha \,x\right)\,\cos\left(\pi \,\alpha \,y\right)\\
&+3M\,\pi ^2\,\alpha ^2\,{\cos\left(t\right)}^3\,{\cos\left(\pi \,\alpha \,x\right)}^3\,{\cos\left(\pi \,\alpha \,y\right)}^3\\
&-6M\,\pi ^2\,\alpha ^2\,{\cos\left(t\right)}^3\,{\cos\left(\pi \,\alpha \,x\right)}^3\,\cos\left(\pi \,\alpha \,y\right)\,{\sin\left(\pi \,\alpha \,y\right)}^2\\
&-M\pi ^2\,\alpha ^2\,\cos\left(t\right)\,\cos\left(\pi \,\alpha \,x\right)\,\cos\left(\pi \,\alpha \,y\right)\\
&-\sin\left(t\right)\,\cos\left(\pi \,\alpha \,x\right)\,\cos\left(\pi \,\alpha 
\,y\right).
\end{split}
\end{equation}


We solve the Cahn--Hilliard equation to a final time $t_F=0.1$, varying the polynomial order, the element spacing, and the time step size. We check the $L^2$ norm of the error, defined as
\begin{equation}
\mathrm{error} = \Vert \phi - \phi_0 \Vert_{\mathcal J,N} = \sqrt{\sum_{e}\left\langle\mathcal J\left(\Phi - \Phi_0\right),\left(\Phi-\Phi_0\right)\right\rangle_{E,N}}.
\end{equation}
The physical parameters are set to $M=1$ and $\varepsilon = 0.1$, and we will vary $\alpha$ analyse space under--resolved and time under--resolved solutions.

\subsubsection{Polynomial order convergence study (p--refinement)}\label{subsec:Numerical:ConvergenceAnalysis:p}

\begin{figure}[ht]
	\centering
	\subfigure[Without interface stabilization ($\kappa_\sigma=0$)]{\includegraphics[scale=0.22]{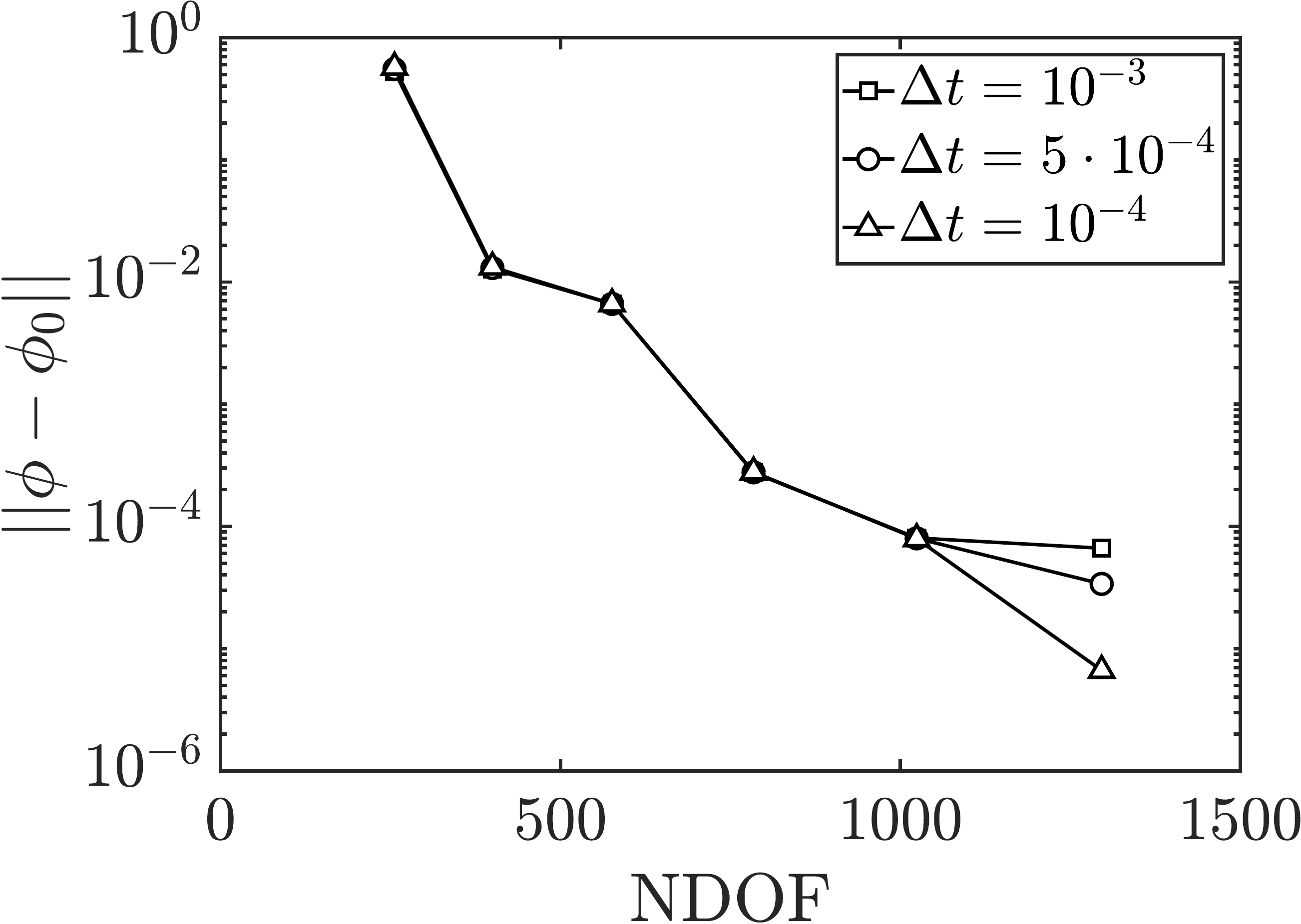}\label{fig:Conv:BR1_sigma0}}
	\subfigure[With interface stabilization ($\kappa_\sigma=3$)]{\includegraphics[scale=0.22]{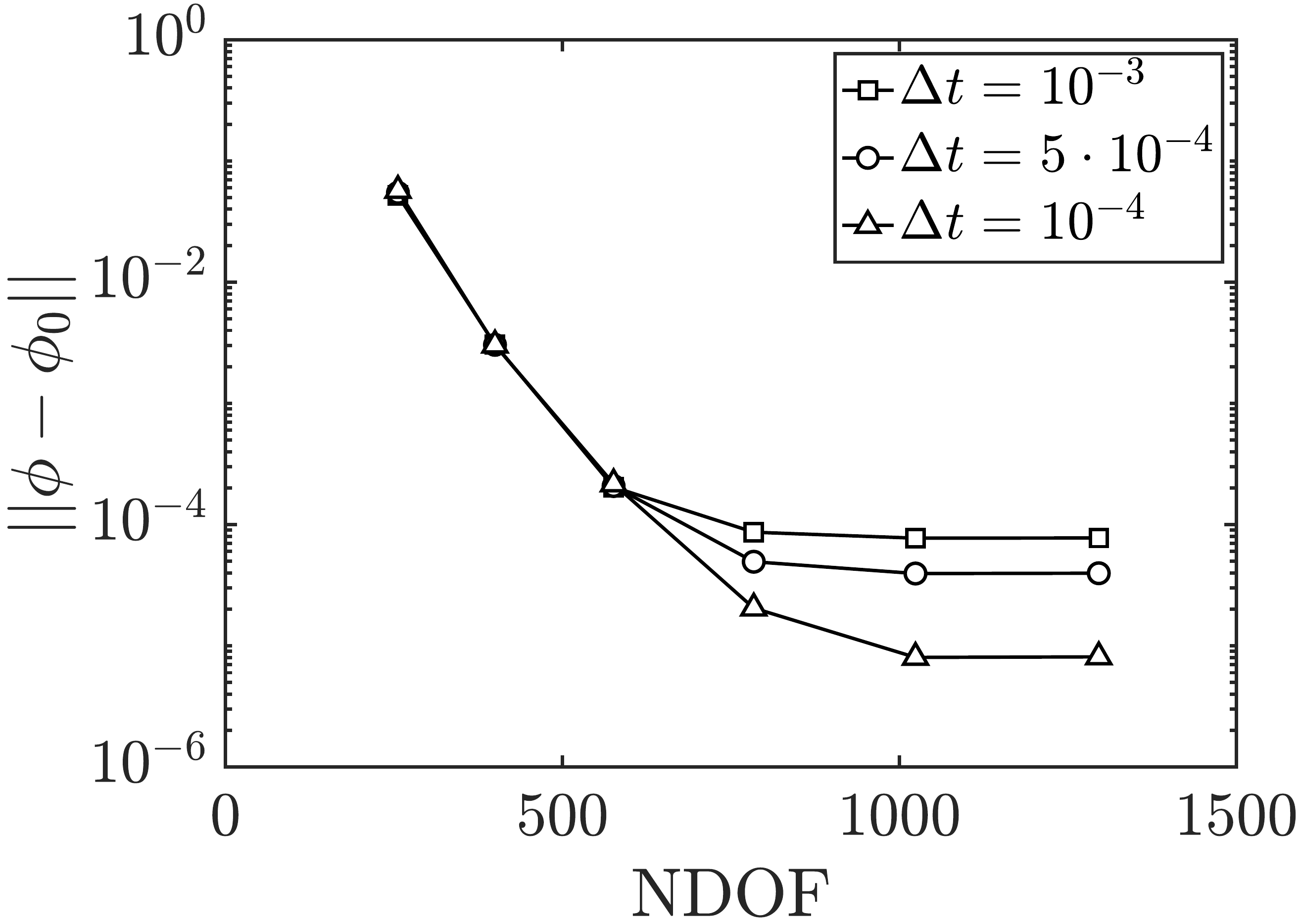}\label{fig:Conv:BR1_sigma1}}
	\caption{Polynomial order convergence study without and with interface stabilization. We show that the convergence is faster and smoother with interface stabilization, and that the convergence rate is exponential only with interface stabilization}
	\label{fig:Conv:BR1}
\end{figure}

In this test the mesh is a fixed $4\times 4$ Cartesian mesh, and  the polynomial order (represented in the $x$--axis of the figures) 
ranges from $N=3$ to $N=8$, with three time steps $\Delta t = 10^{-3}, 5\cdot 10^{-4},$ and $10^{-4}$. 
In all these tests we take $\alpha=1$, so that the solution is well--resolved in space with relatively coarse meshes. The results are shown in Fig. \ref{fig:Conv:BR1}, where we performed 
the study for two penalty parameter coefficient values $\kappa_\sigma$, which was defined in \eqref{eq:PenaltyParameter}: without stabilization ($\kappa_\sigma=0$) in Fig. \ref{fig:Conv:BR1_sigma0} 
and with stabilization  ($\kappa_\sigma=3$) in Fig. \ref{fig:Conv:BR1_sigma1}.

\begin{figure}[h]
	\centering
	\includegraphics[scale=0.3]{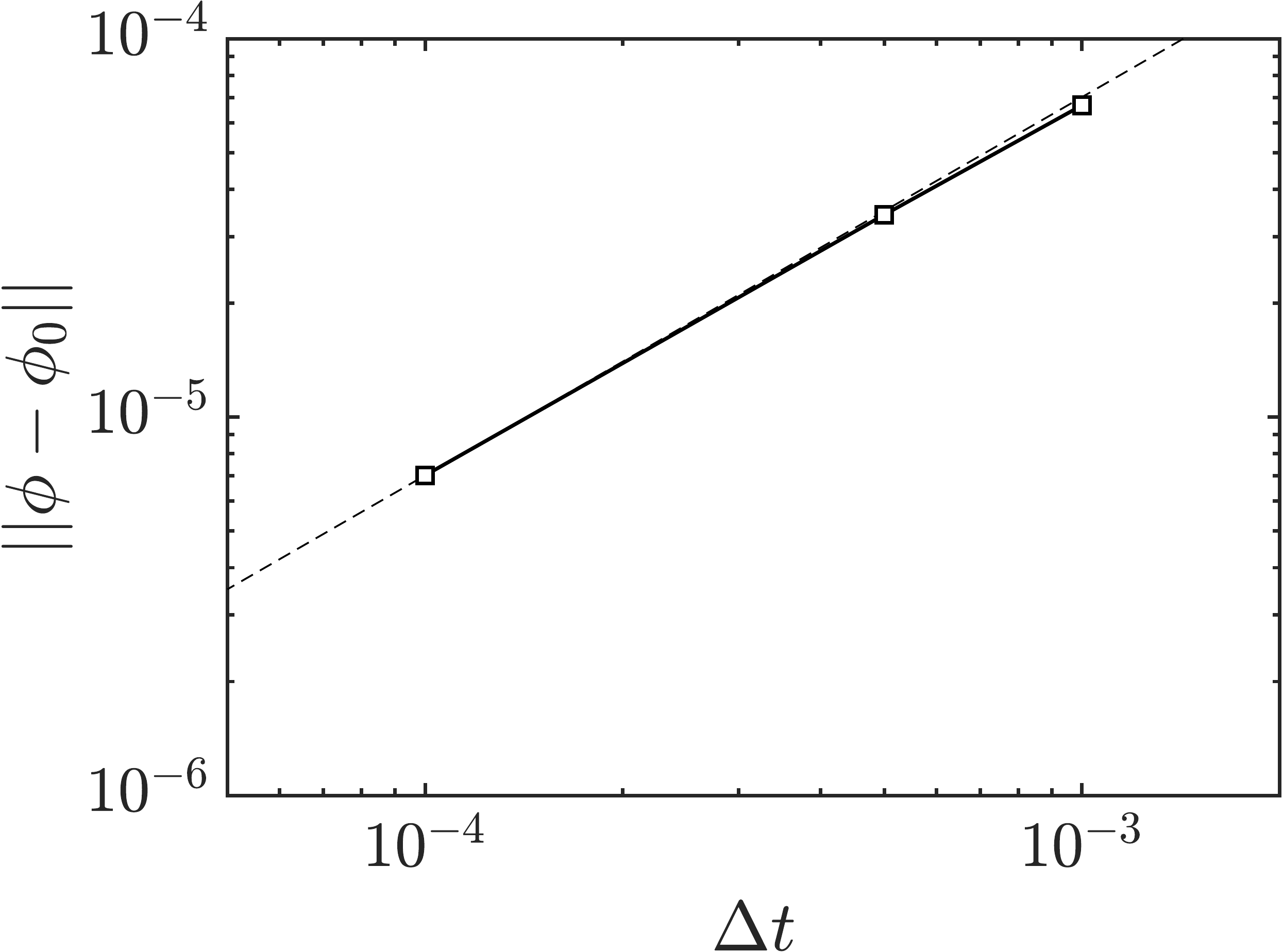}
	\caption{Temporal convergence study. The dashed line represents the theoretical linear convergence rate. We confirm that the scheme is first order accurate in time, as designed}
	\label{fig:Conv:time-convergence}
\end{figure}

We find that the convergence is slower without interface stabilization ($\kappa_\sigma=0$), and also that it is uneven. This even--odd phenomena has been also reported in \cite{2016:Gassner} in the context of the compressible Euler equations, where it was tackled by adding interface stabilization. With interface stabilization (Fig. \ref{fig:Conv:BR1_sigma1}), we find not only that the convergence is smoother, but errors are always lower. The curve shows the typical pattern of a polynomial order convergence study: on the one hand, for low polynomial orders the solution is under--resolved in space, and errors decrease exponentially with the polynomial order (linear decay in semi--logarithmic plot). In this region, errors are not affected by the time step $\Delta t$. On the other hand, for high polynomial orders, the solution is under--resolved in time, and thus it reaches a stagnation with further increase of the polynomial order. In this region, the error is controlled by the time step $\Delta t$, as the different plots in Fig. \ref{fig:Conv:BR1_sigma1} show. We have represented in Fig. \ref{fig:Conv:time-convergence} the different errors obtained once the stagnation is reached to show that the scheme is first order accurate in 
time, as designed.

\subsubsection{Mesh convergence study (h--refinement)}

\begin{figure}[h]
\centering
\includegraphics[scale=0.3]{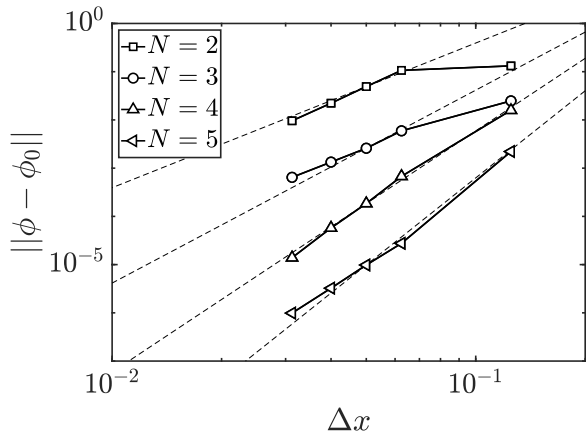}
\caption{Spatial convergence study for different polynomial orders. We vary the mesh spacing from $\Delta x = 0.125$ (16x16 mesh) to $\Delta x=0.03125$ (64x64 mesh). Dashed lines represent the theoretical convergence rate, $\Delta x^{N+1}$, for each polynomial}
\label{fig:Conv:h-convergence}
\end{figure}

In this test, we increase the manufactured solution's wavenumber to $\alpha=8$, so that we extend the region with spatial under--resolution. We vary the number of elements from a $16\times 16$ mesh to a $64\times 64$ mesh, and we consider four polynomial orders, from $N=2$ to $N=5$. The rest of the parameters remain the same as in Sec. \ref{subsec:Numerical:ConvergenceAnalysis:p}, and we use the scheme with interface stabilization. The results are represented in Fig. \ref{fig:Conv:h-convergence}, where we have drawn the theoretical convergence rates for each polynomial order
\begin{equation}
\Vert \phi-\phi_0 \Vert \propto \Delta x^{N+1}.
\end{equation}

We find that despite adding the interface stabilization, there still remains some slight even--odd effect, as the convergence is faster than the theoretical for even polynomial orders, and slower otherwise. Similar behavior was noted in \cite{2004:Ainsworth} and the reason behind this effect in this particular scheme still remains an open question.

\subsection{Two dimensional spinodal decomposition}\label{subsec:Numerical:SpinodalDecomposition}

\begin{figure}[h]
	\centering
	\subfigure[Cartesian mesh]{\includegraphics[scale=0.25]{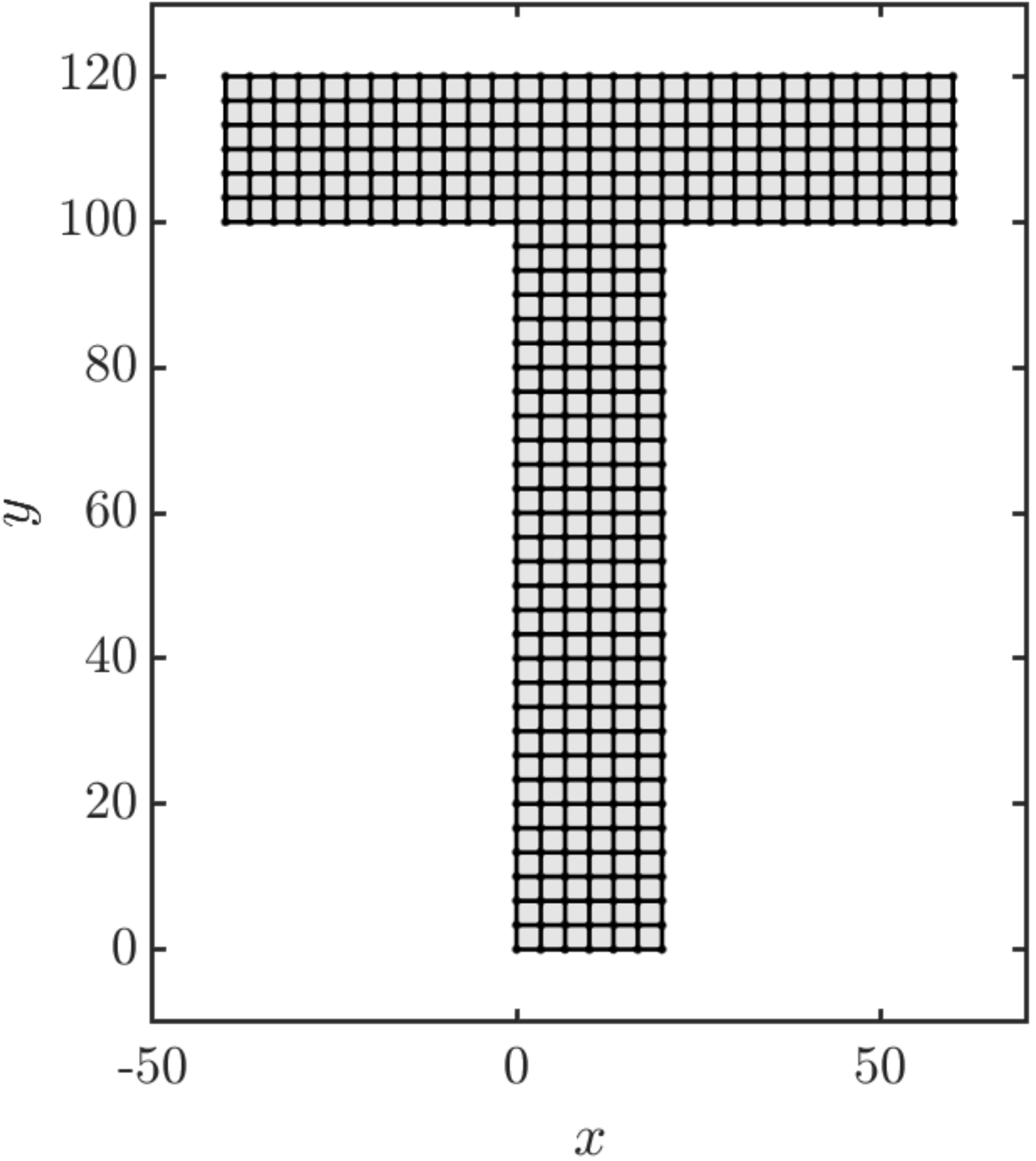}\label{fig:Numerical:Jokisaari:T-Straight}}
	\subfigure[Distorted mesh (straight edges)]{\includegraphics[scale=0.25]{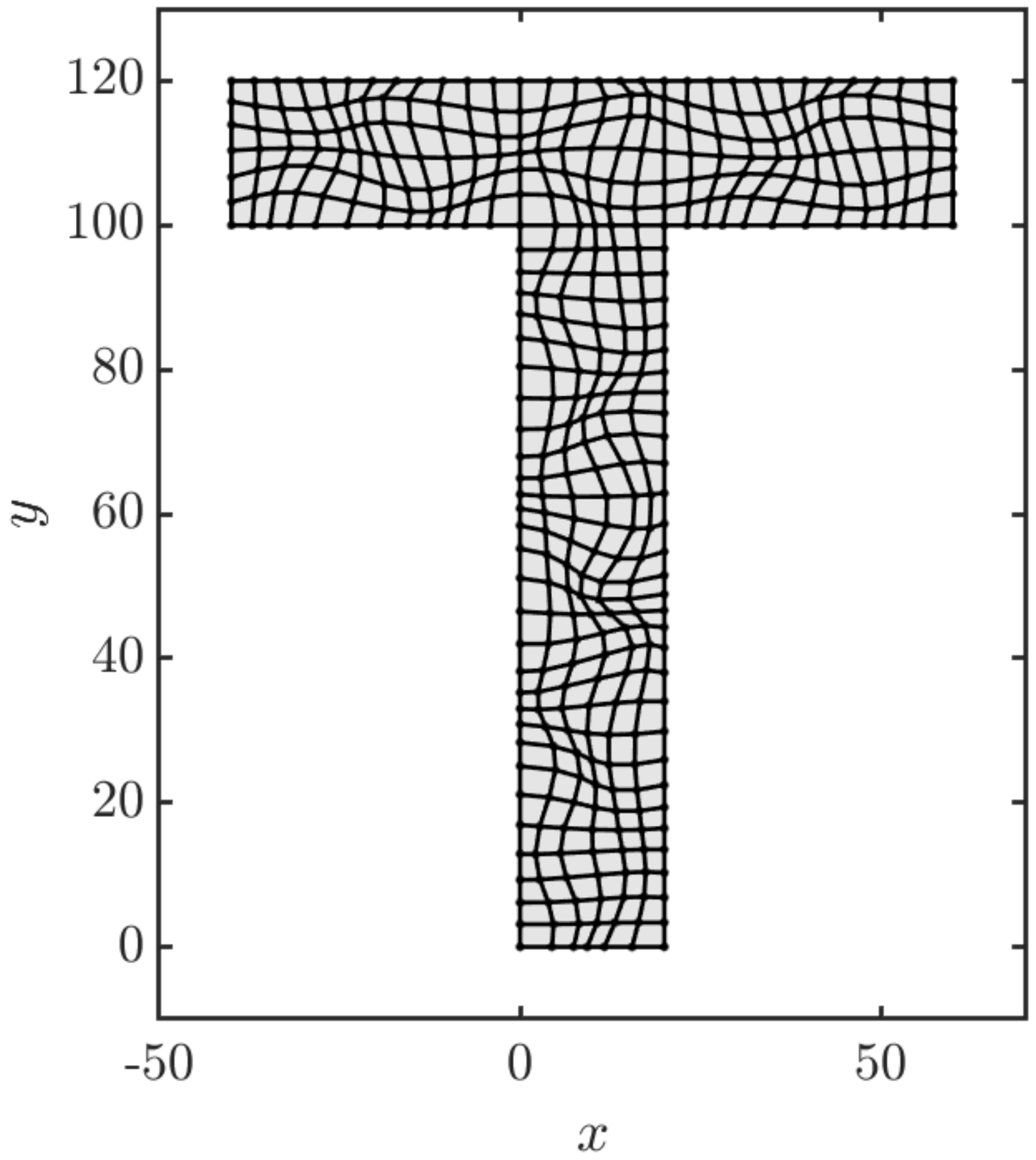}\label{fig:Numerical:Jokisaari:T-Distorted}}	
	\subfigure[Curved mesh (curvilinear edges)]{\includegraphics[scale=0.25]{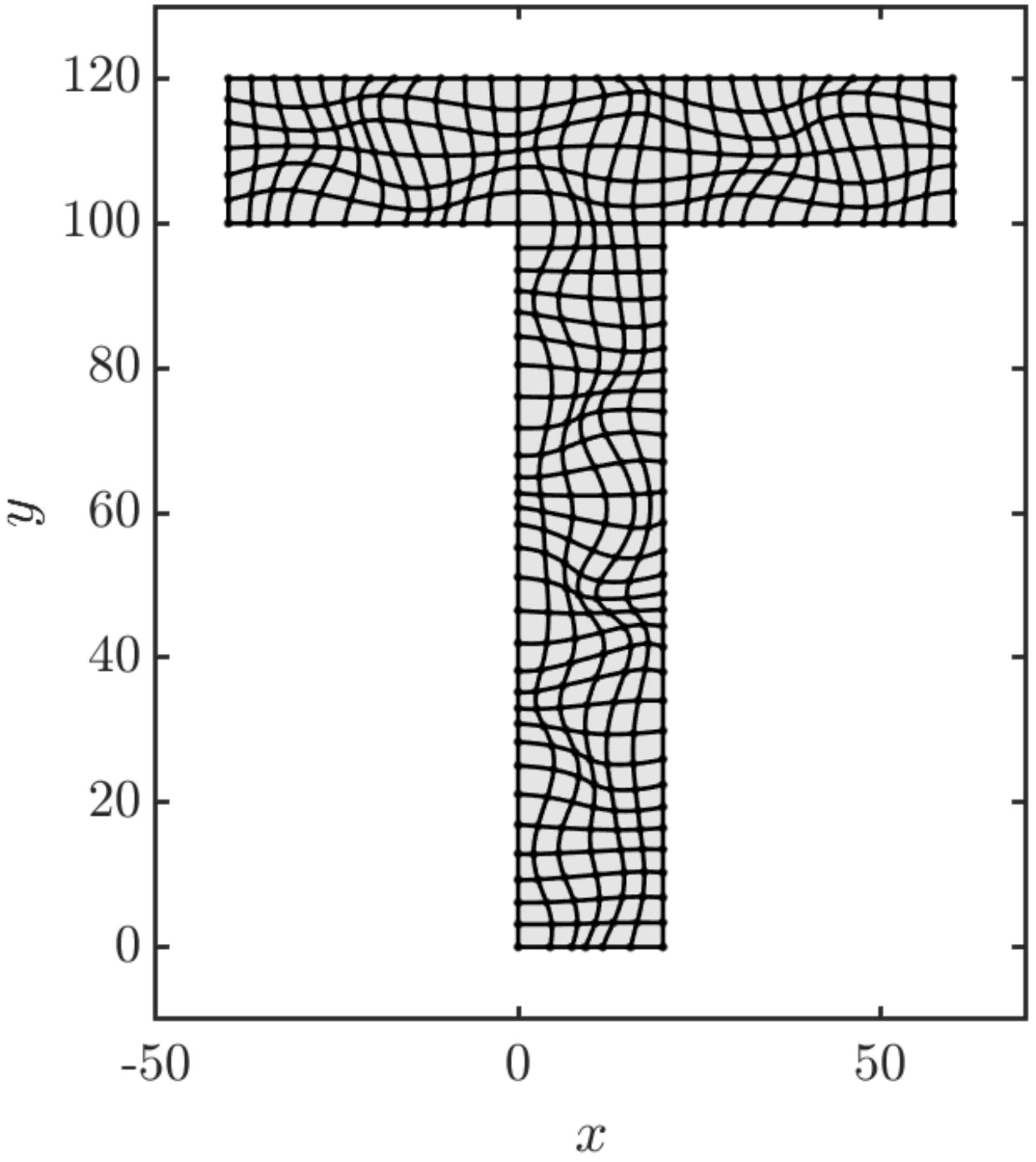}\label{fig:Numerical:Jokisaari:T-Curved}}	
	\subfigure[Unstructured quad mesh]{\includegraphics[scale=0.25]{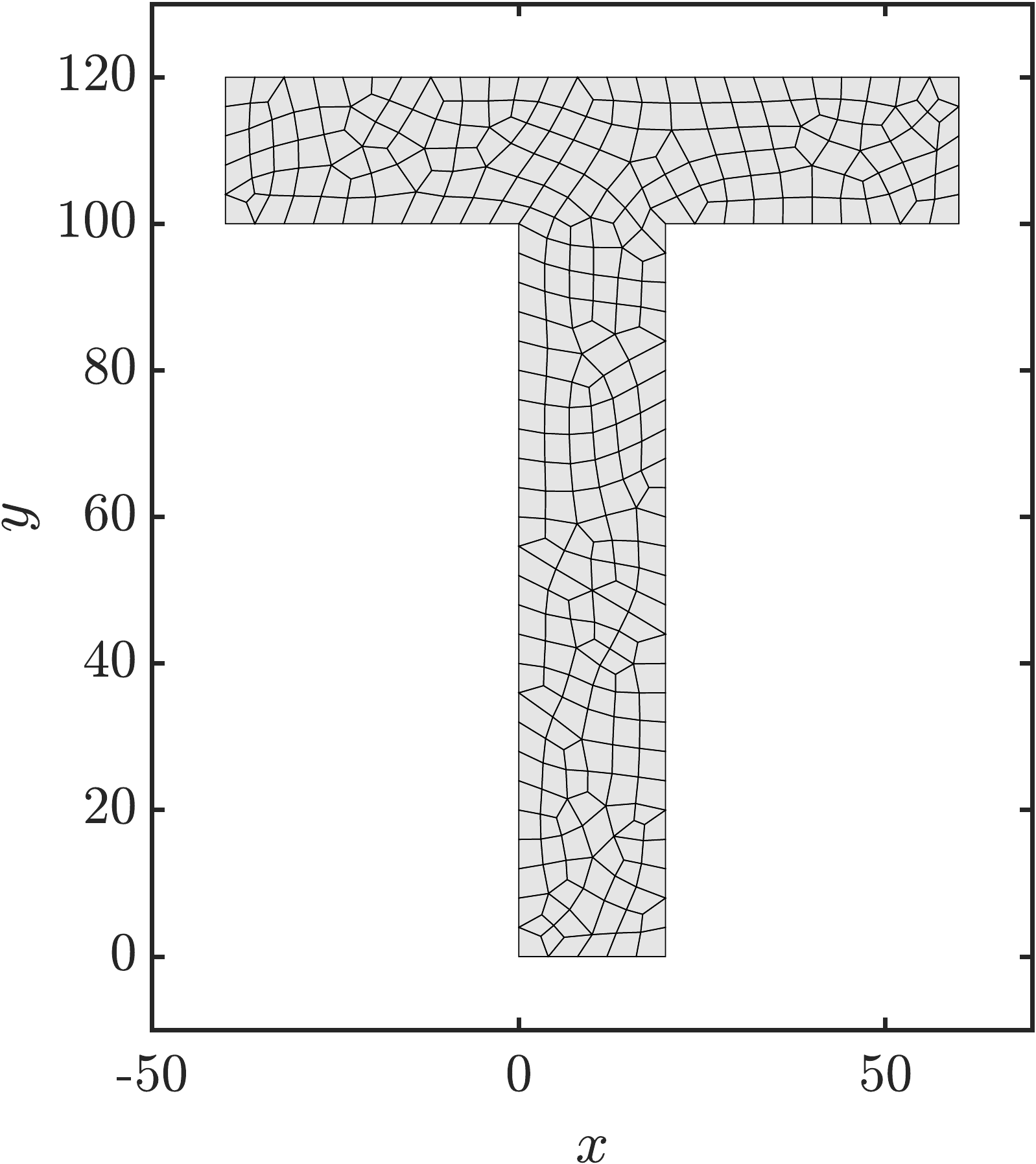}\label{fig:Numerical:Jokisaari:T-Unstructured}}	
\caption{Four meshes used for the spinodal decomposition test. The distorted mesh has straight faces, whilst the curved mesh edges are curvilinear}
\label{fig:Numerical:Jokisaari:T-meshes}	
\end{figure}

In this section, we compare the scheme when solving a two dimensional spinodal decomposition. The spinodal decomposition describes the phase separation process from an initial mixed state. To trigger the separation, we use the initial condition proposed in \cite{2016:Jokisaari},
\begin{equation}
\begin{split}
\phi_0\left(x,y\right) =& 0.05\bigl(\cos\left(0.105 x\right)\cos\left(0.11 z\right) + \left(\cos\left(0.13 x\right)\cos\left(0.087 z\right)\right)^2 \\
&+ \cos\left(0.025 x-0.15 z\right) \cos\left(0.07 x - 0.02 z\right)\bigr).
\end{split}
\end{equation}
The domain is ``T''--shaped, which we mesh using four strategies: a Cartesian mesh (Fig. \ref{fig:Numerical:Jokisaari:T-Straight}), the Cartesian mesh distorted with straight sided elements (Fig. \ref{fig:Numerical:Jokisaari:T-Distorted}), the Cartesian mesh distorted with curved elements (Fig. \ref{fig:Numerical:Jokisaari:T-Curved}), and an unstructured quad mesh (Fig. \ref{fig:Numerical:Jokisaari:T-Unstructured}).

The parameters of the Cahn--Hilliard equation used are that of \cite{2016:Jokisaari}, but adapted to $\phi$ ranging 
from $-1$ to $1$ (in \cite{2016:Jokisaari} it ranges from 0.3 to 0.7),
\begin{equation}
M = 1.0, ~~ \varepsilon = 3.1623, ~~ \Delta t = 0.1,
\end{equation}
and we use a polynomial order $N=4$. We monitor the free--energy evolution with time 
and compare with the results provided in \cite{2016:Jokisaari}.

We first study the effect of the IMEX parameters $K_0$ and 
$S_0$ in the curved (Fig. \ref{fig:Numerical:Jokisaari:T-Curved}) and the unstructured (Fig. \ref{fig:Numerical:Jokisaari:T-Unstructured}) meshes. For $K_0$ we consider two scenarios: the Crank--Nicolson scheme, $K_0=1/2$, 
and the backward Euler scheme $K_0=1$.  For $S_0$ we study 
three values: $S_0=0$, which does not guarantee stability for the non--linear 
terms, $S_0=1$, which guarantees non--linear term stability for 
$\Phi\in[-1,1]$, and $S_0=2$, which guarantees non--linear term stability for 
$\Phi\in[-1.29,1.29]$ (see Fig. \ref{fig:pi-isolines}).

\begin{figure}[ht]
\centering
\subfigure[Curved mesh]{\includegraphics[height=4.0cm]{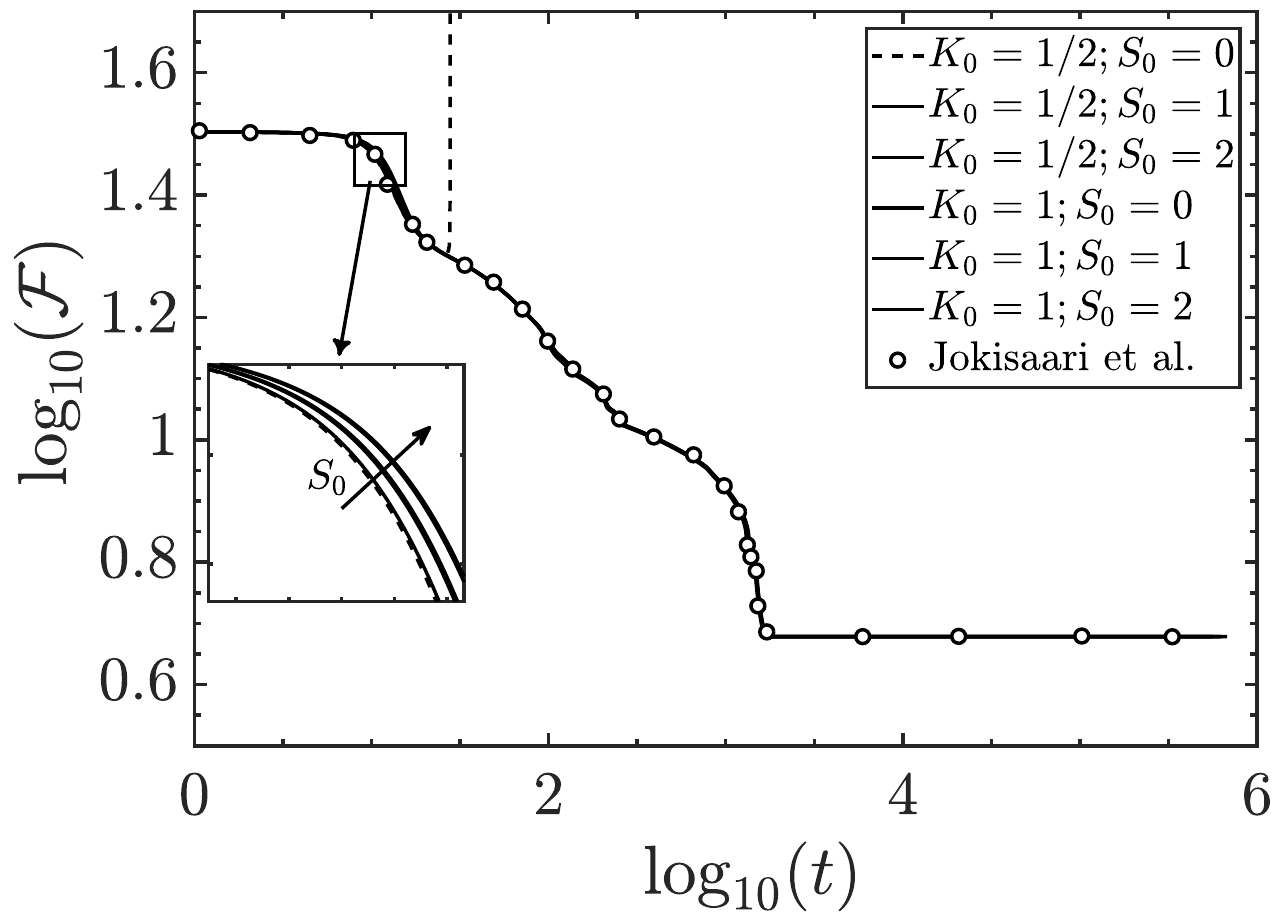}\label{fig:Numerical:Jokisaari:T-Curved-Results}}
\subfigure[Unstructured mesh]{\includegraphics[height=4.0cm]{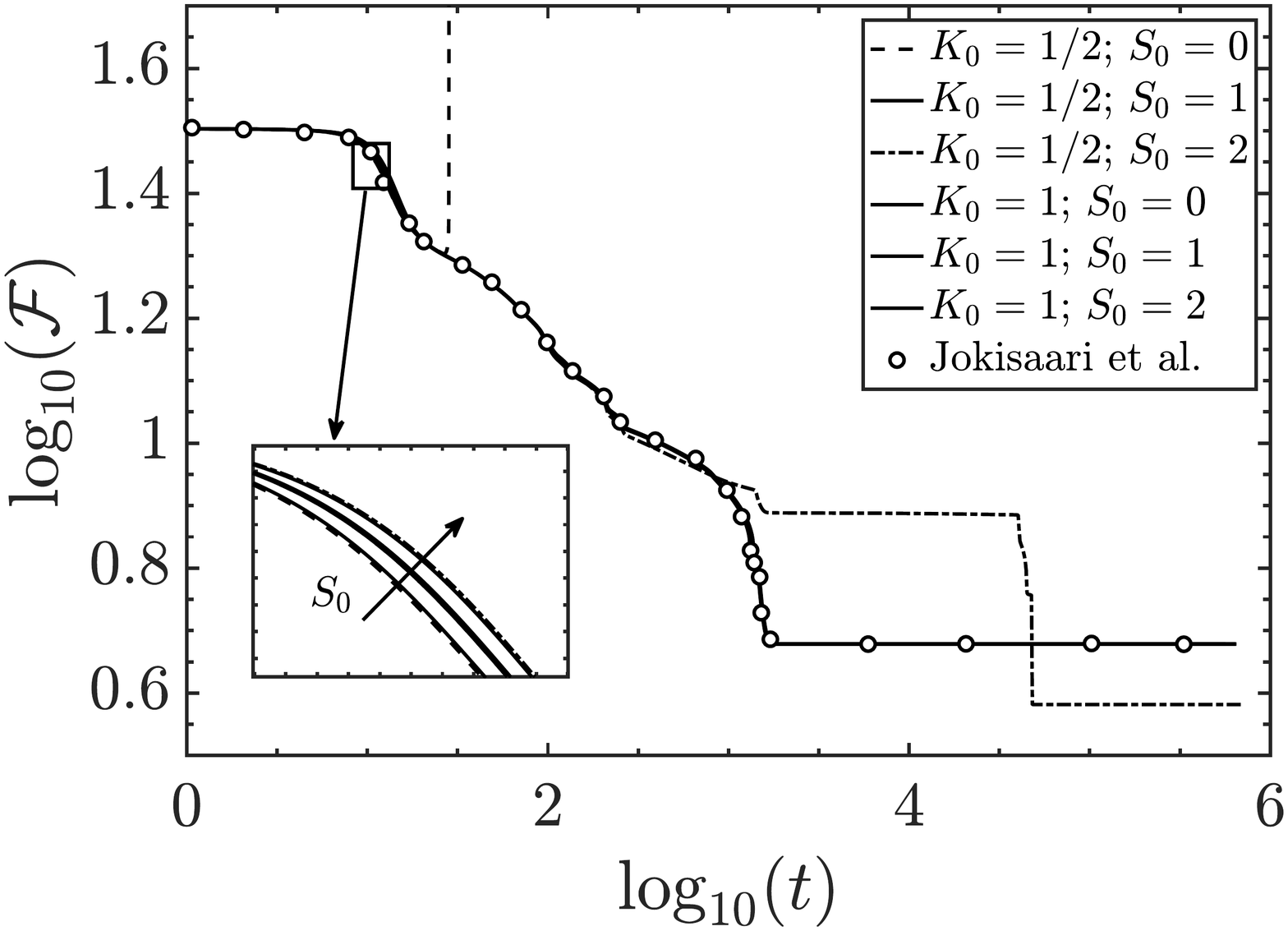}\label{fig:Numerical:Jokisaari:T-Unstructured-Results}}
\caption{Evolution of the free--energy $\mathcal F$ for the ``T'' domain. We consider two $K_0$ cases: Crank--Nicolson ($K_0=1/2$) and
backward Euler ($K_0=1$). We find that the solution depends on $S_0$ only if we use the Crank--Nicolson scheme, 
producing an
unstable scheme if $S_0=0$ (dashed line), and a different solution for $S_0=2$ (dot--dash line in the unstructured mesh). For the rest of the cases (solid lines), 
the solution agrees with \cite{2016:Jokisaari}, represented with circles}
\label{fig:Numerical:Jokisaari:T-Result}
\end{figure}

The evolution of the free energy is shown in Fig. \ref{fig:Numerical:Jokisaari:T-Result} 
for all parameter combinations and the two meshes (curvilinear in Fig. \ref{fig:Numerical:Jokisaari:T-Curved-Results}, and unstructured in Fig. \ref{fig:Numerical:Jokisaari:T-Unstructured-Results}). The circles represent the reference solution 
from \cite{2016:Jokisaari}. On the one hand, 
for the backward Euler scheme ($K_0=1$, solid lines in Fig. \ref{fig:Numerical:Jokisaari:T-Result}), there is 
no high impact of $S_0$ on the free energy. 
The scheme is stable even for the lowest $S_0$ value. Thus, the physical dissipation and the dissipation 
introduced by backward Euler scheme are enough to balance non--linear 
instabilities in this case. On the other hand, for the Crank--Nicolson scheme ($K_0=1/2$), the solution depends highly 
on $S_0$. For $S_0=0$ (represented with a dashed line), the scheme is unstable. The 
physical dissipation is not sufficient to counteract the instabilities that arise from 
the non--linear terms. For $S_0=1$ (solid line), we 
obtain approximately the same solution as the backward Euler scheme in both meshes. However, for 
$S_0=2$ (dot--dash line), and only for the unstructured mesh, the final solution is different from the others, leading to 
lower free--energy values in the steady--state.  Furthermore, in the detailed view (Fig. \ref{fig:Numerical:Jokisaari:T-Result}) we find small differences in the free--energy evolution when varying $S_0$. In a non--intuitive way, the free--energy decreases at a slower rate for higher $S_0$ values. Nonetheless, this result is still in agreement with the theoretical bound \eqref{eq:FullyDiscrete:stability-all-timesteps}.

Overall, we recommend the 
use of the backward Euler scheme to remove the final solution dependency on $S_0$, and to use $S_0=1$ 
to avoid  instabilities from non--linear terms, although they have not appeared 
in this experiment (with $K_0=1$). 

\begin{figure}[h]
	\centering
	\subfigure[Curved mesh]{\includegraphics[height=5cm]{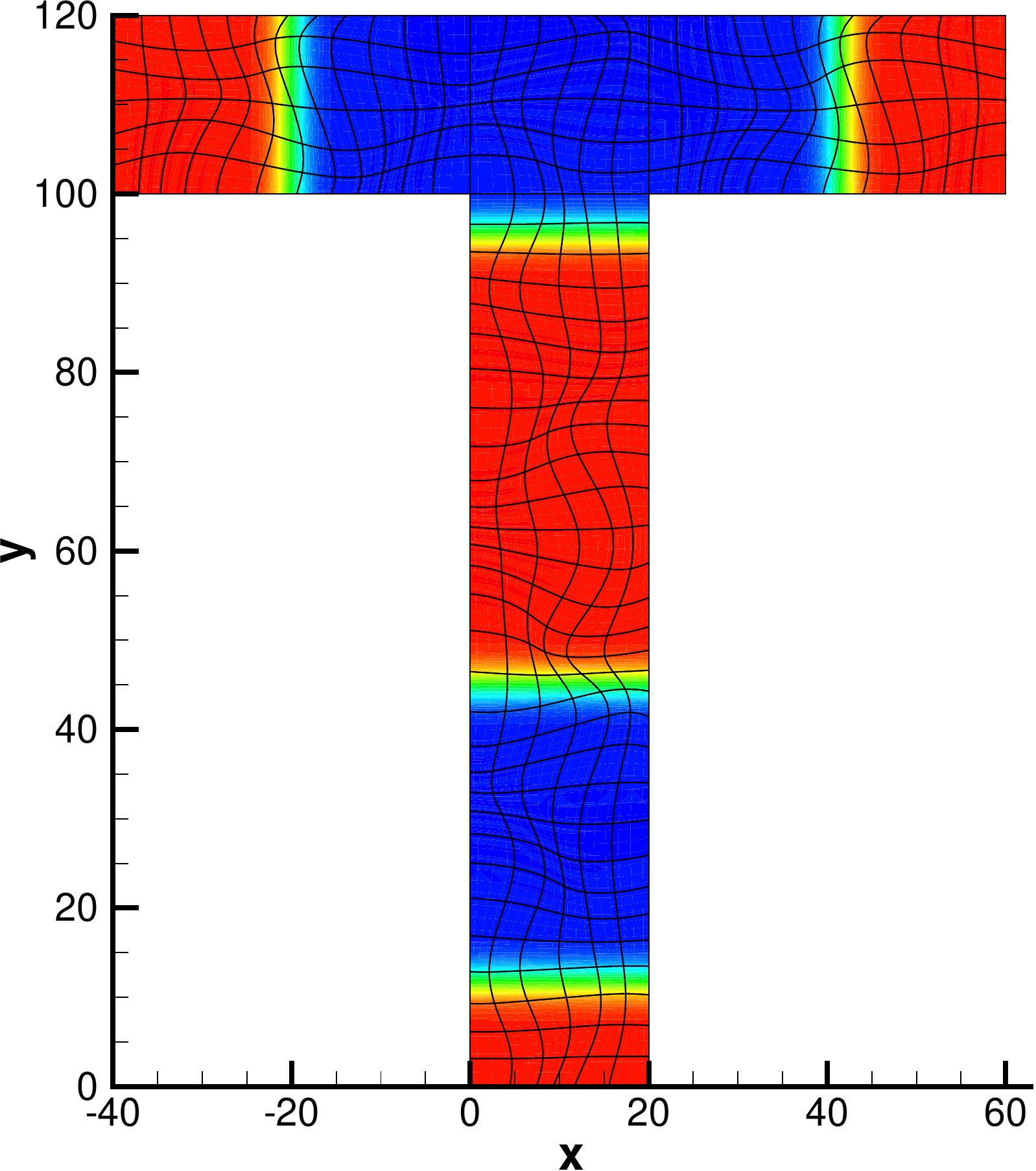}}
	\subfigure[Unstructured mesh]{\includegraphics[height=5cm]{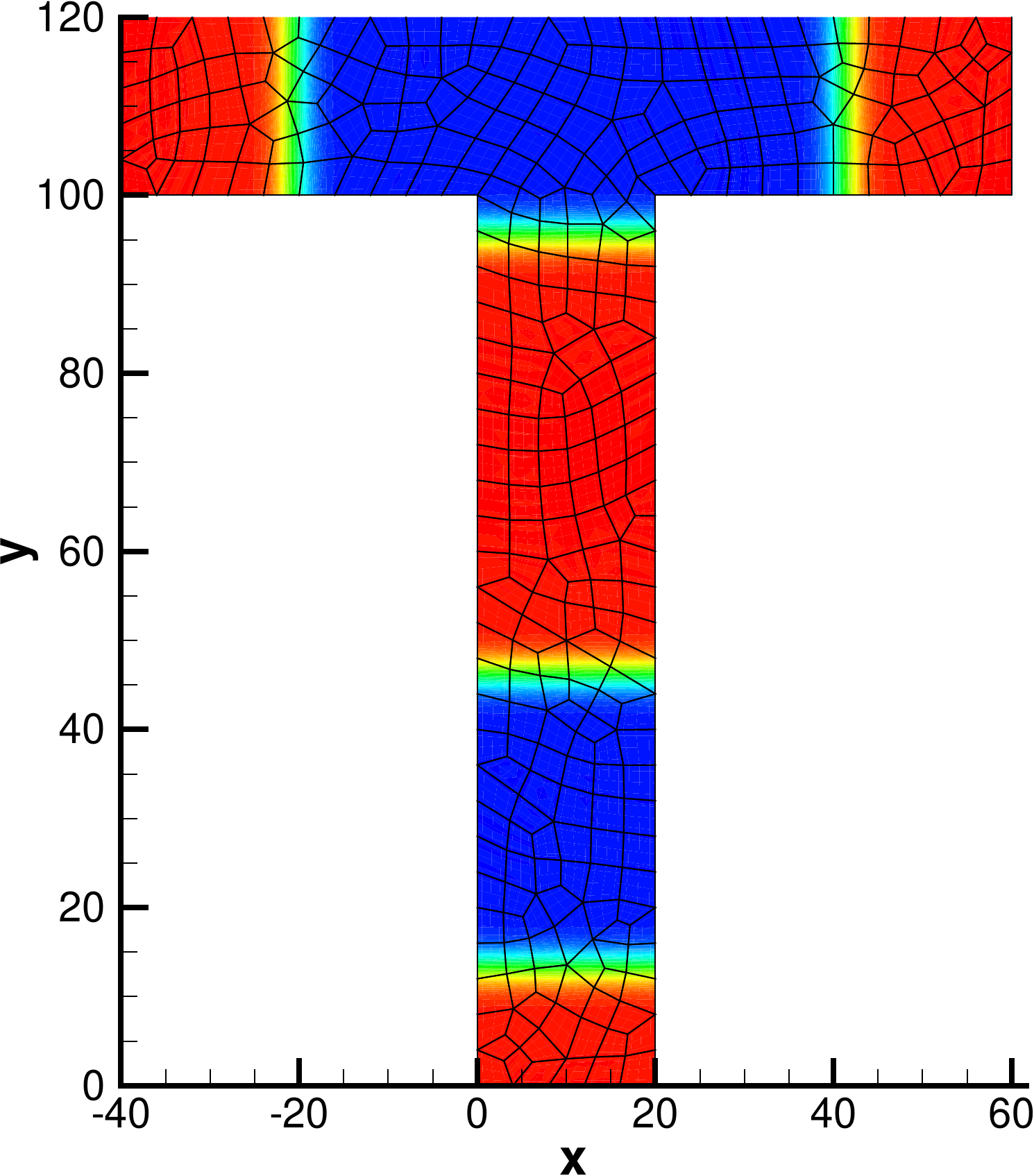}}	
	\caption{Representation of the final equilibrium solution in the curved and unstructured mesh. We find that the impact from the mesh on the solution is minimal, since interfaces are not aligned with the elements faces}
	\label{fig:Numerical:Jokisaari:FinalState}
\end{figure}

For this scenario, we represent the differences in the free--energy evolution for the rest of the meshes considered 
in Fig. \ref{fig:Numerical:Jokisaari:T-meshes}. First, in Fig. \ref{fig:Numerical:Jokisaari:FinalState} we represent 
the final state for the curvilinear and unstructured meshes, where we find that 
the mesh has no visual impact on the solution, since the resulting interfaces are not aligned 
with the mesh.

\begin{figure}[h]
	\centering
	\includegraphics[scale=0.2]{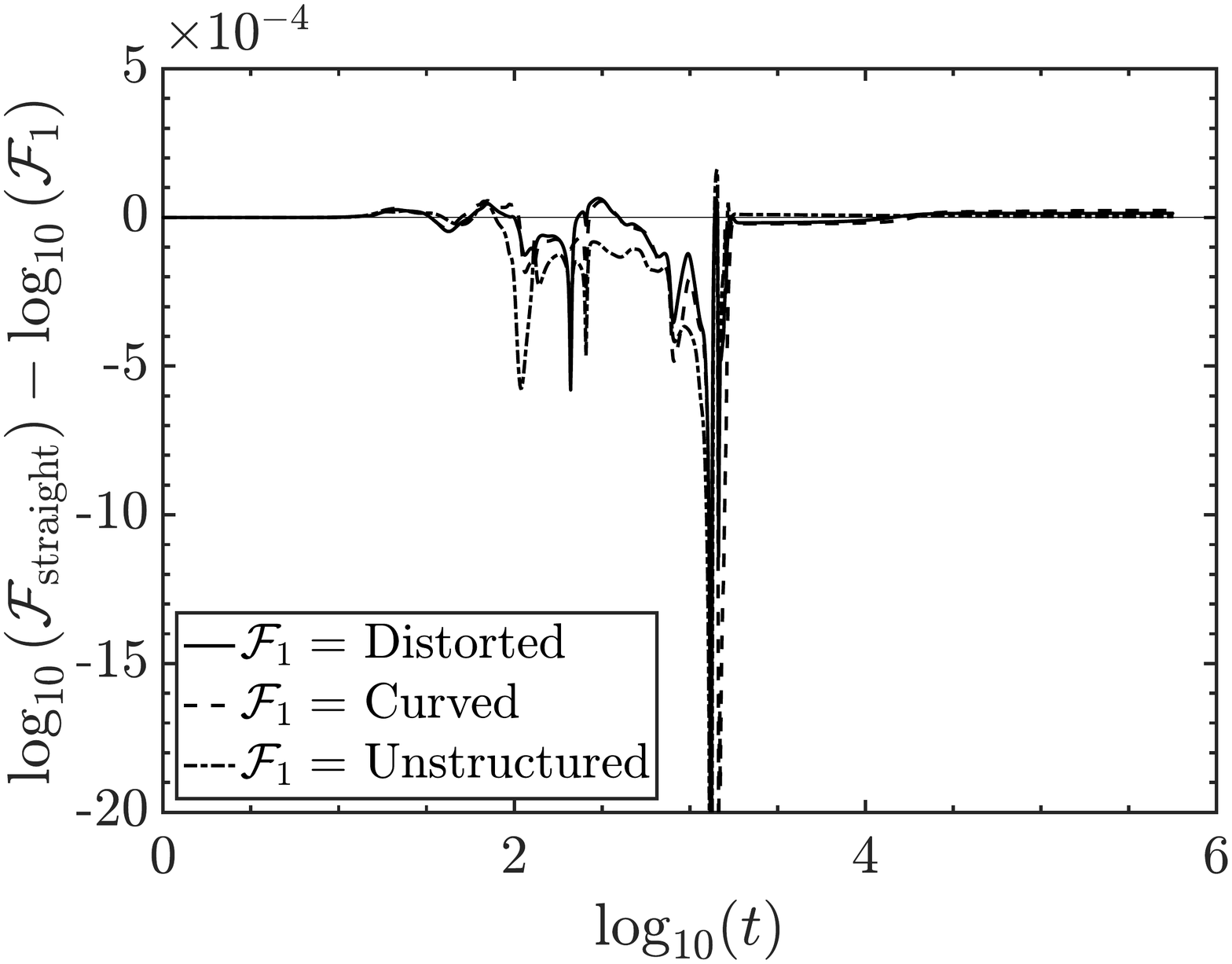}
	\caption{Free--energy evolution comparison of the distorted (solid line), curved (dashed line), and unstructured (dash--dot line) 
	meshes with the Cartesian mesh. We find that, despite the initial and final values being
	 approximately the same, in the evolution the errors are maintained low (of order $10^{-3}$) for the same number
	 of degrees of freedom}
	\label{fig:Numerical:Jokisaari:FreeEnergyComparison}
\end{figure}

To quantify the differences due to the mesh, we consider the Cartesian mesh 
(Fig. \ref{fig:Numerical:Jokisaari:T-Straight}) as the reference solution, and represent the 
free--energy difference of the other meshes with that. The result is
 represented in Fig. \ref{fig:Numerical:Jokisaari:FreeEnergyComparison}. In the $y$--axis 
 we represent the difference $\log_{10}\left( \mathcal F_{\text{cartesian}}\right)-\log_{10}\left(\mathcal F_{1}\right)$, 
 between the Cartesian mesh and $\mathcal F_1$ stands for the free--energy of the distorted 
 (solid line), curved (dashed line) and unstructured (dash--dot line). In any case, the errors (and more substantially, in the initial and final states) 
 remain low, and all meshes considered were stable and accurate once the IMEX parameters 
 were set-up appropriately.

\subsection{Three dimensional spinodal decomposition}

\begin{figure}[h]
  \centering
  \includegraphics[scale=0.2]{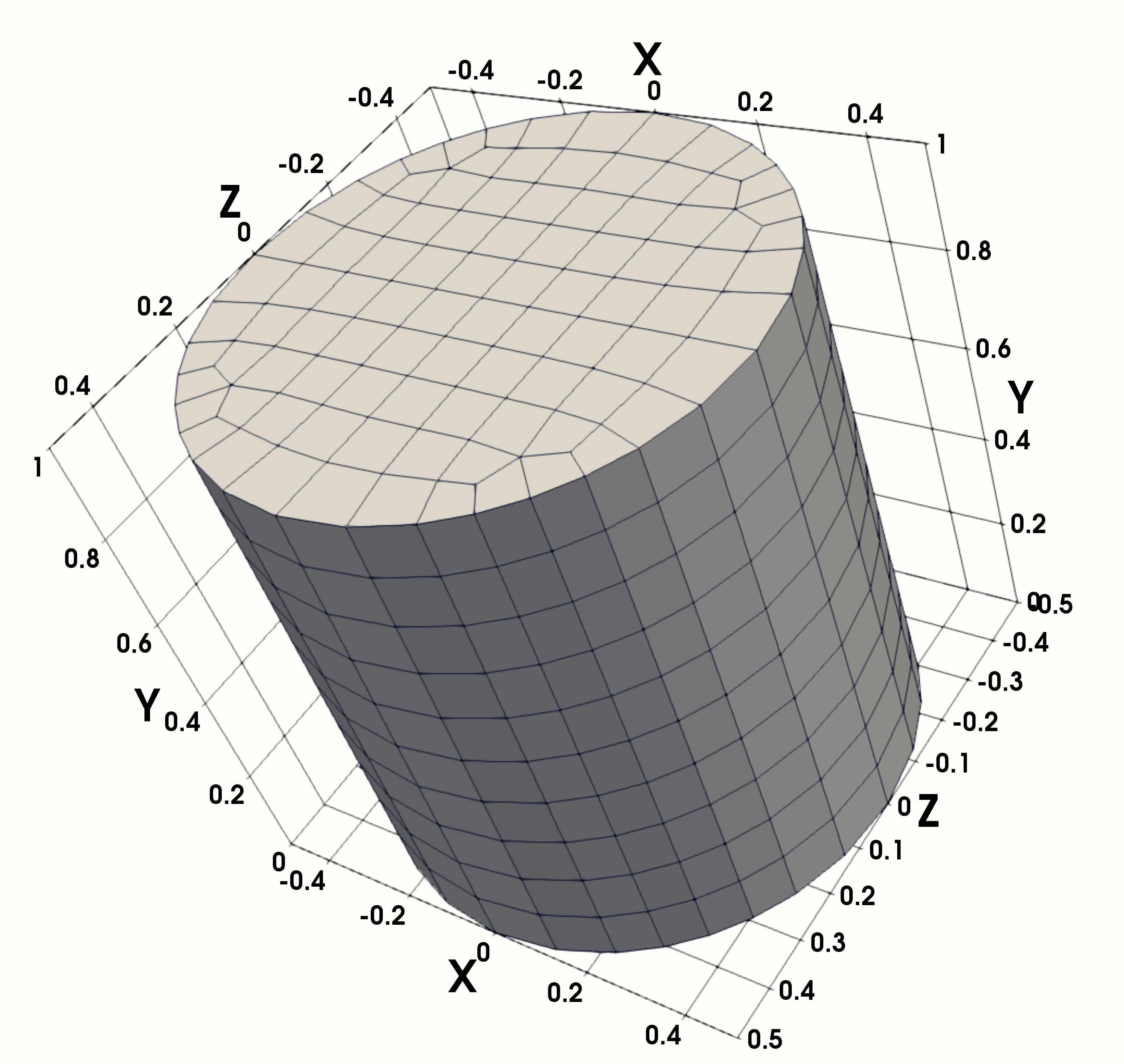}
  \caption{Representation of the mesh used for the three dimensional spinodal decomposition simulation}
  \label{fig:Numerical:3DCyl:Mesh}
\end{figure}

We also show the solver's capability to solve the spinodal decomposition in three 
dimensions. We consider the interior of a cylinder ($L=D=1$), which we divide into $920$ elements. 
A representation of the cylinder and the mesh is provided in Fig. 
\ref{fig:Numerical:3DCyl:Mesh}.
We use order $N=3$ polynomials, and the following values for the 
rest of the parameters:
\begin{equation}
  \varepsilon = 0.1, ~~ M = 1.0, ~~\Delta t = 10^{-4},~~K_0=1,~~S_0=1.0.
\end{equation}
Initially, as in the two dimensional example, the phases are mixed, and we trigger 
the decomposition with the  initial condition,
\begin{equation}
  \begin{split}
\phi_0(x,y,z) = & \phantom{{}+{}} 0.015 \cos(5x-10z) \cos(7x+10zy+1) \\
&+ 0.02\cos(20y^2+15x^2)\sin(5x+2y+3x) \\
&+ 0.02\cos(10\sqrt{y^2+z^2}) \cos(15xy)\sin(20x +10z)\\
&+ 0.01\cos(3x)\cos(3z)\cos(4y).
 \end{split}
\end{equation}

\begin{figure}[h]
  \centering
  \subfigure[$t=0$]{\includegraphics[scale=0.2]{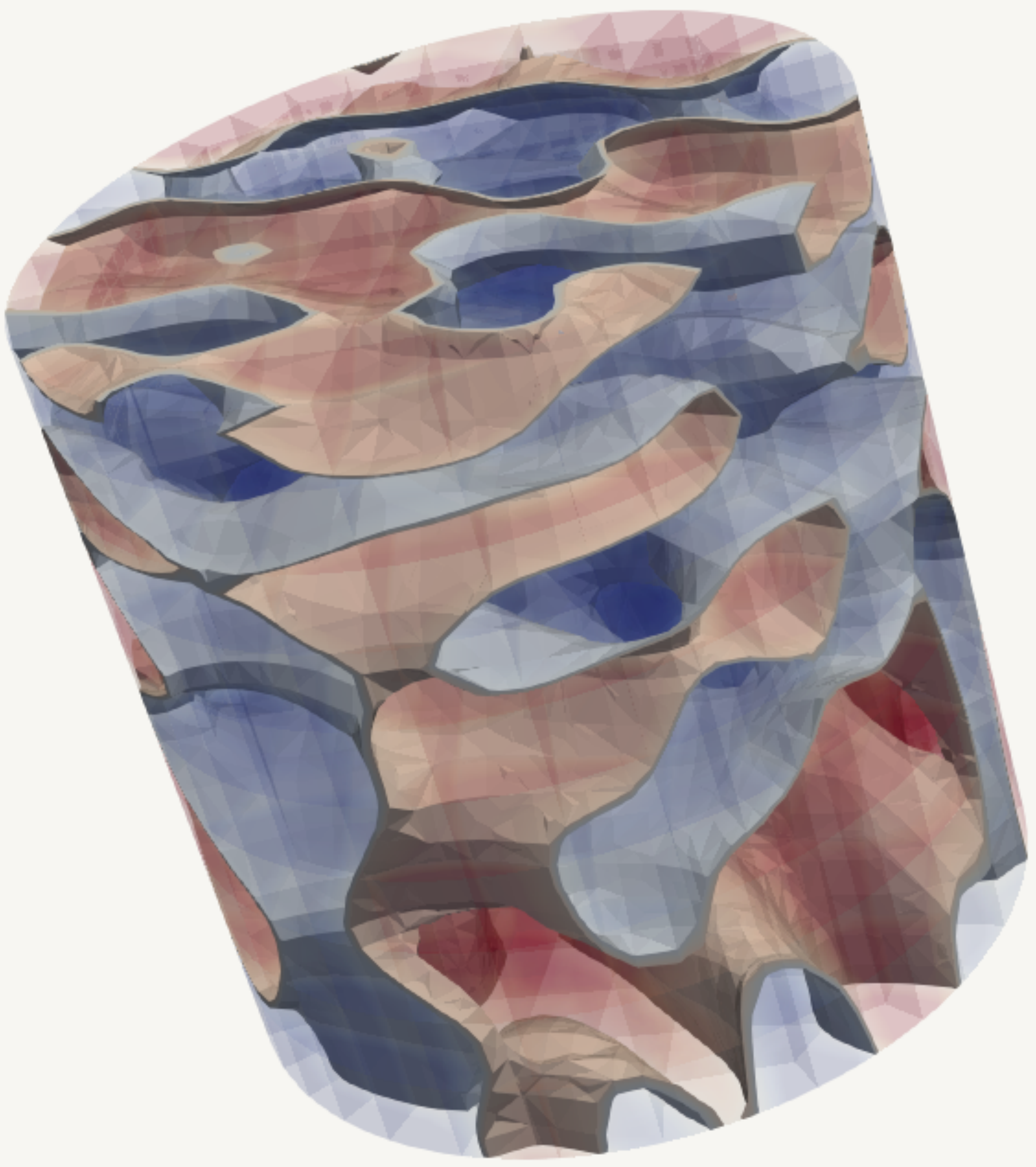}\label{fig:Numerical:3DCyl:IC}}
    \subfigure[$t=10^{-4}$]{\includegraphics[scale=0.2]{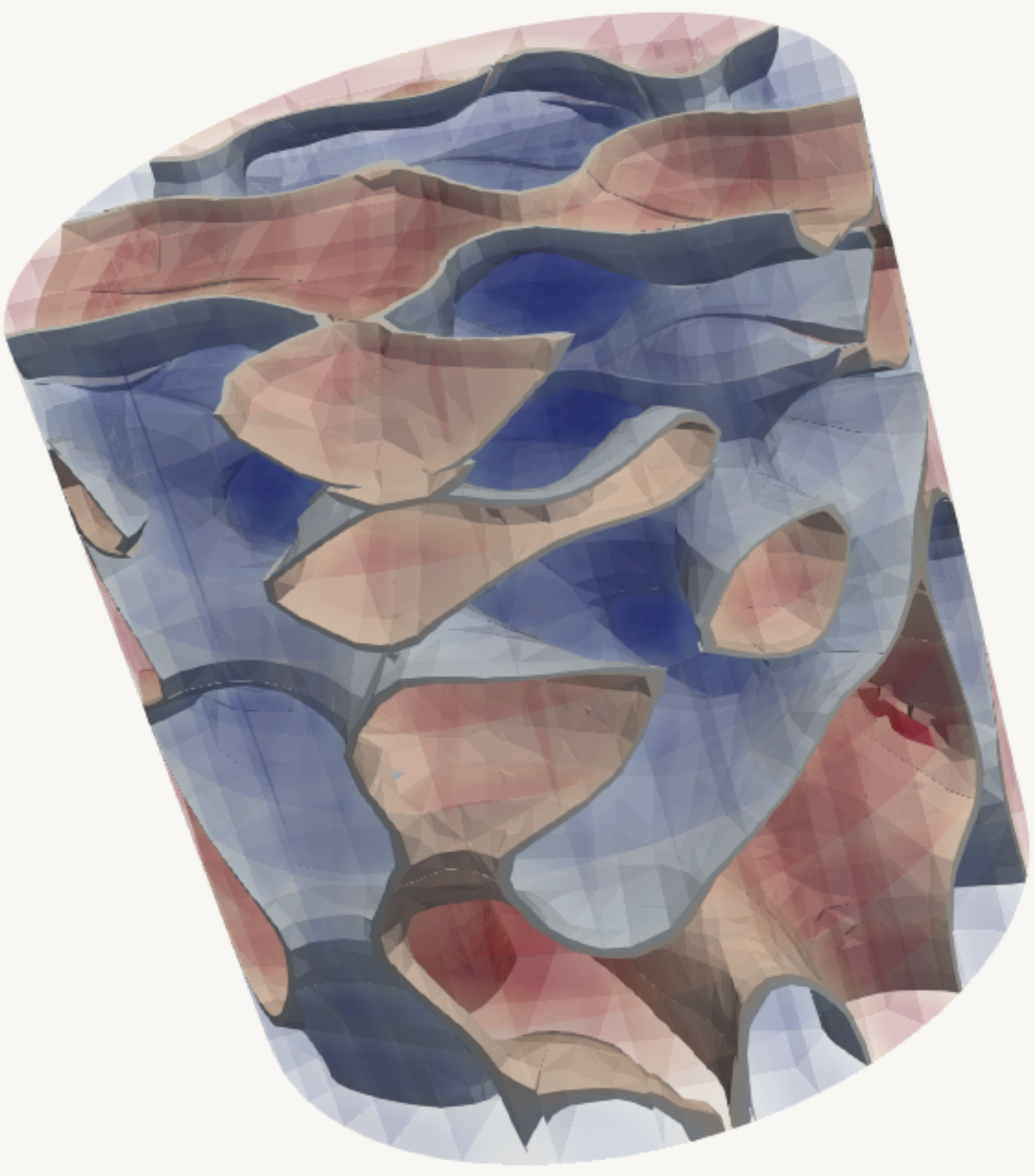}}
      \subfigure[$t=10^{-3}$]{\includegraphics[scale=0.2]{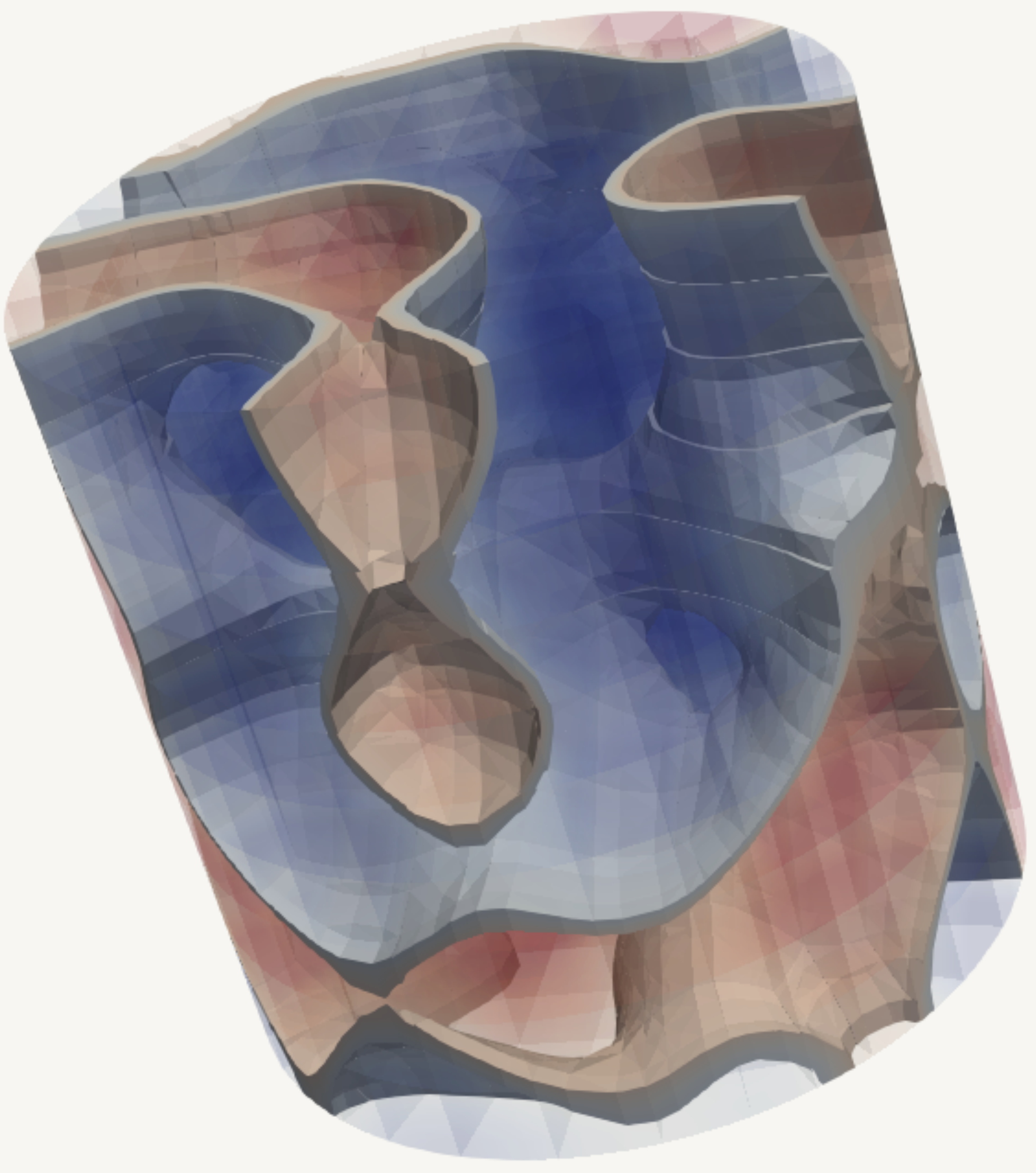}}
        \subfigure[$t=10^{-2}$]{\includegraphics[scale=0.2]{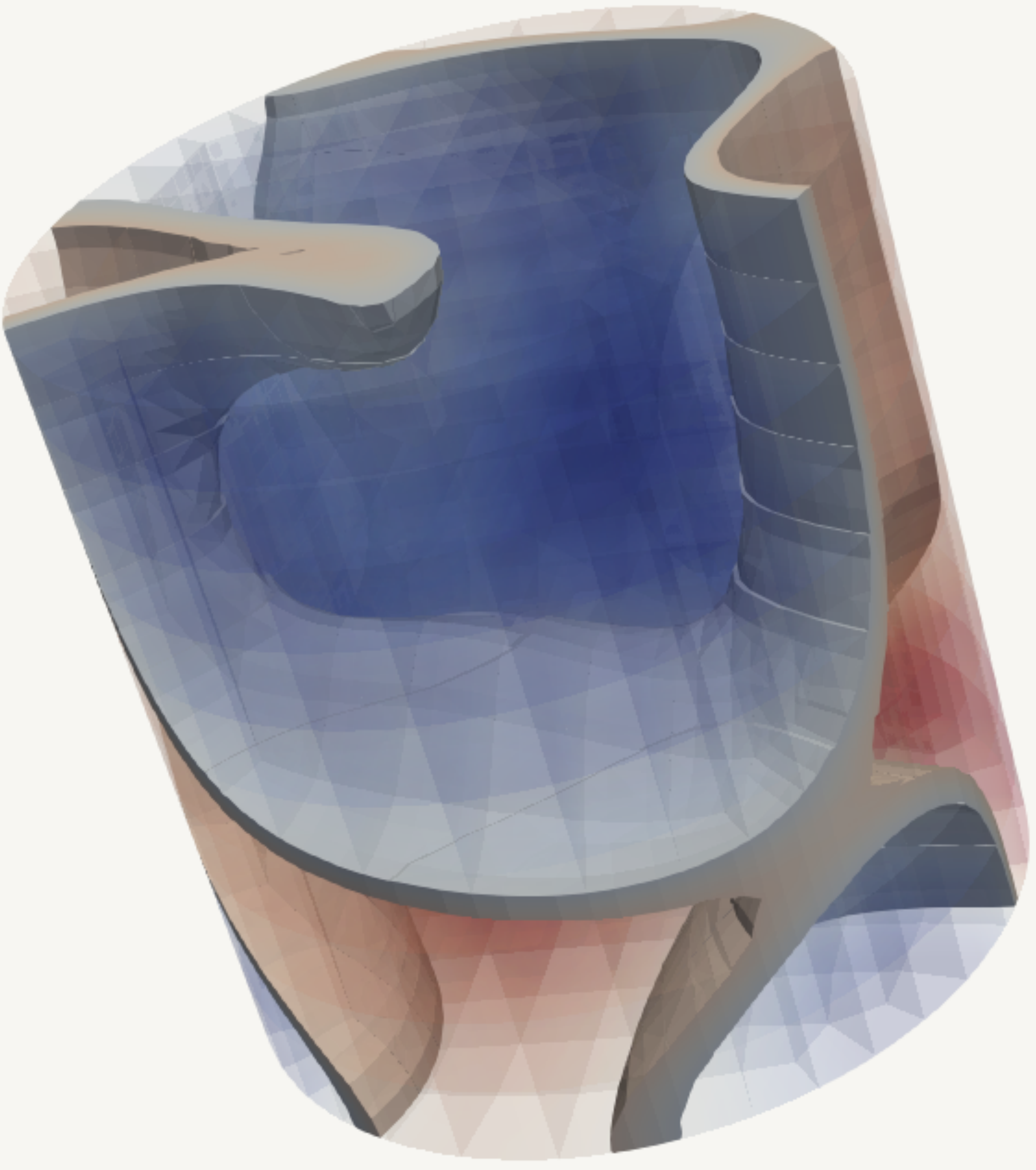}}\\
          \subfigure[$t=0.1$]{\includegraphics[scale=0.15]{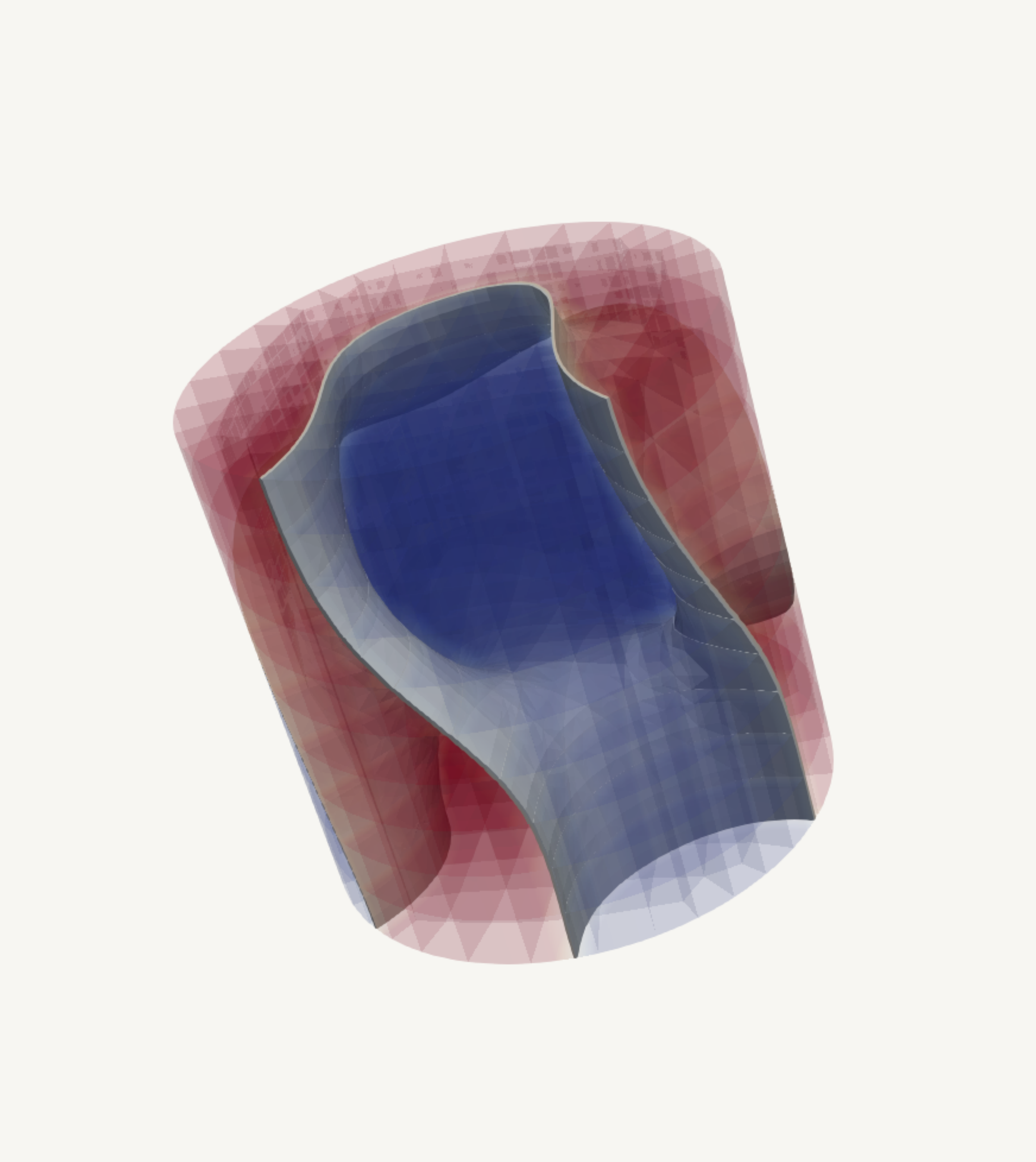}}
            \subfigure[$t=1.0$]{\includegraphics[scale=0.15]{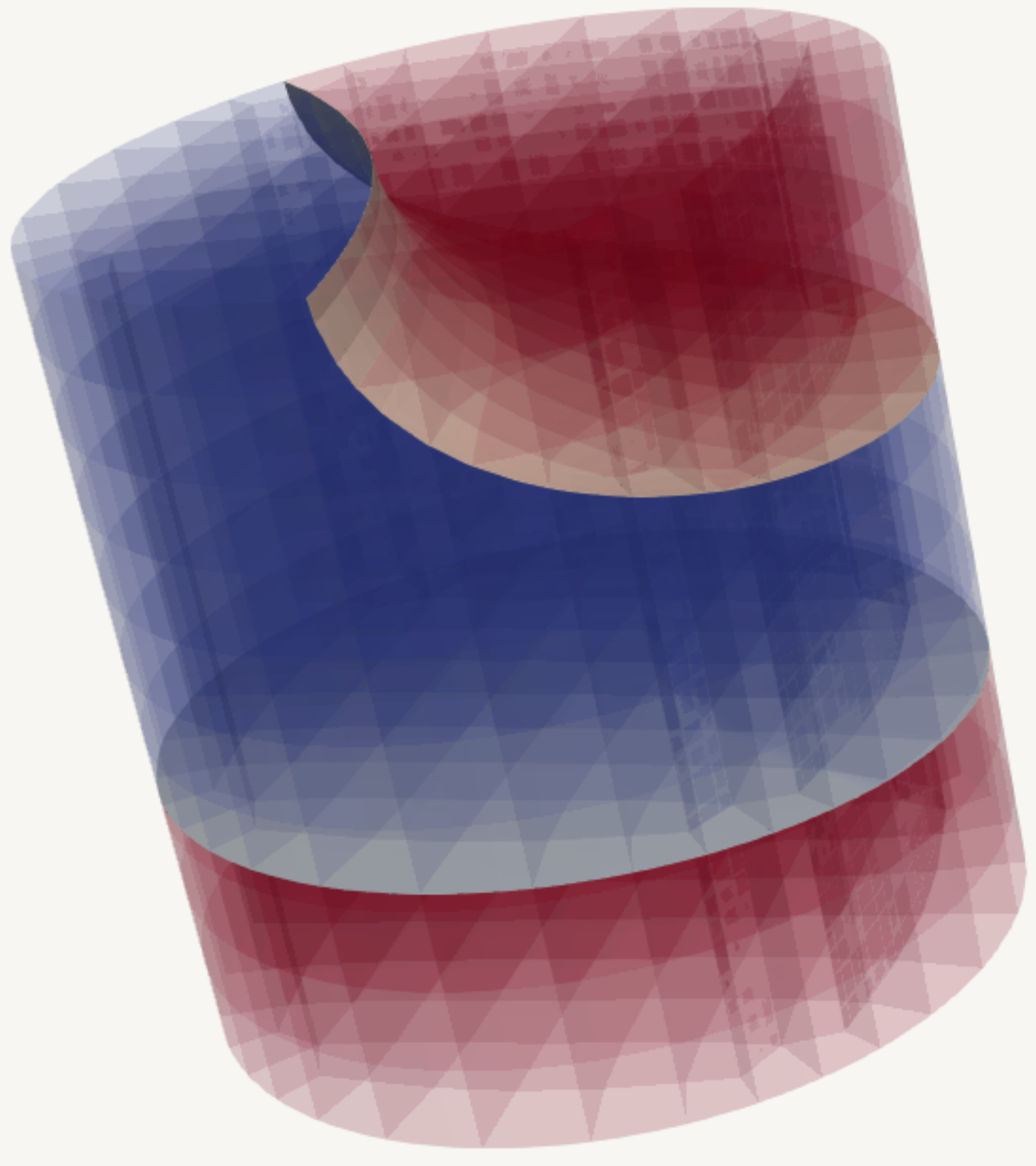}}
              \subfigure[$t=10.0$]{\includegraphics[scale=0.15]{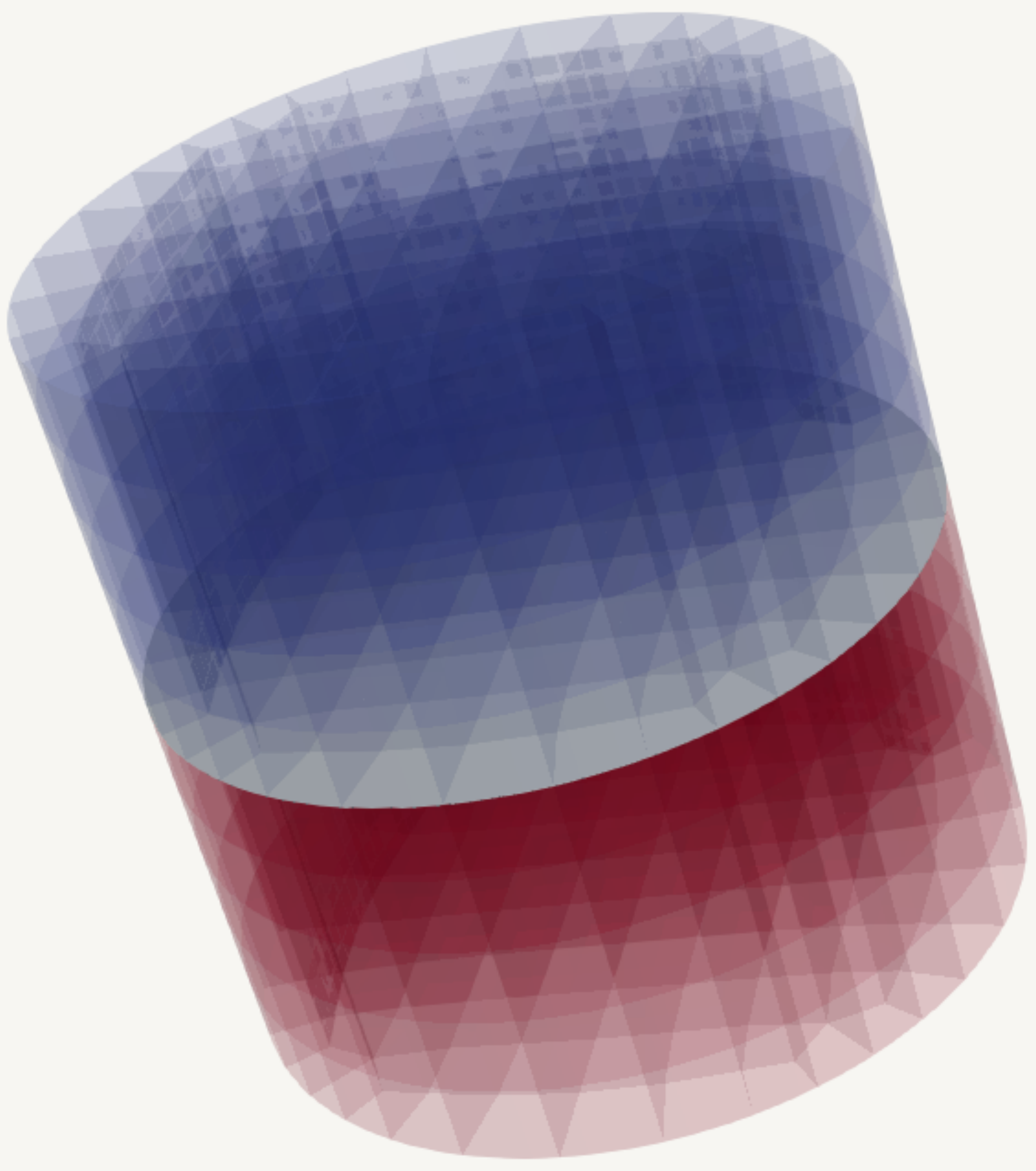}\label{fig:Numerical:3DCyl:SS}}
              \caption{Evolution of the phases with time for the three dimensional spinoidal decomposition. Blue and red contours represent the equilibrium phases $\phi=-1$ and $\phi=1$ respectively}
              \label{fig:Numerical:3DCyl:Evolution}
\end{figure}

\begin{figure}[h]
  \centering
  \includegraphics[scale=0.6]{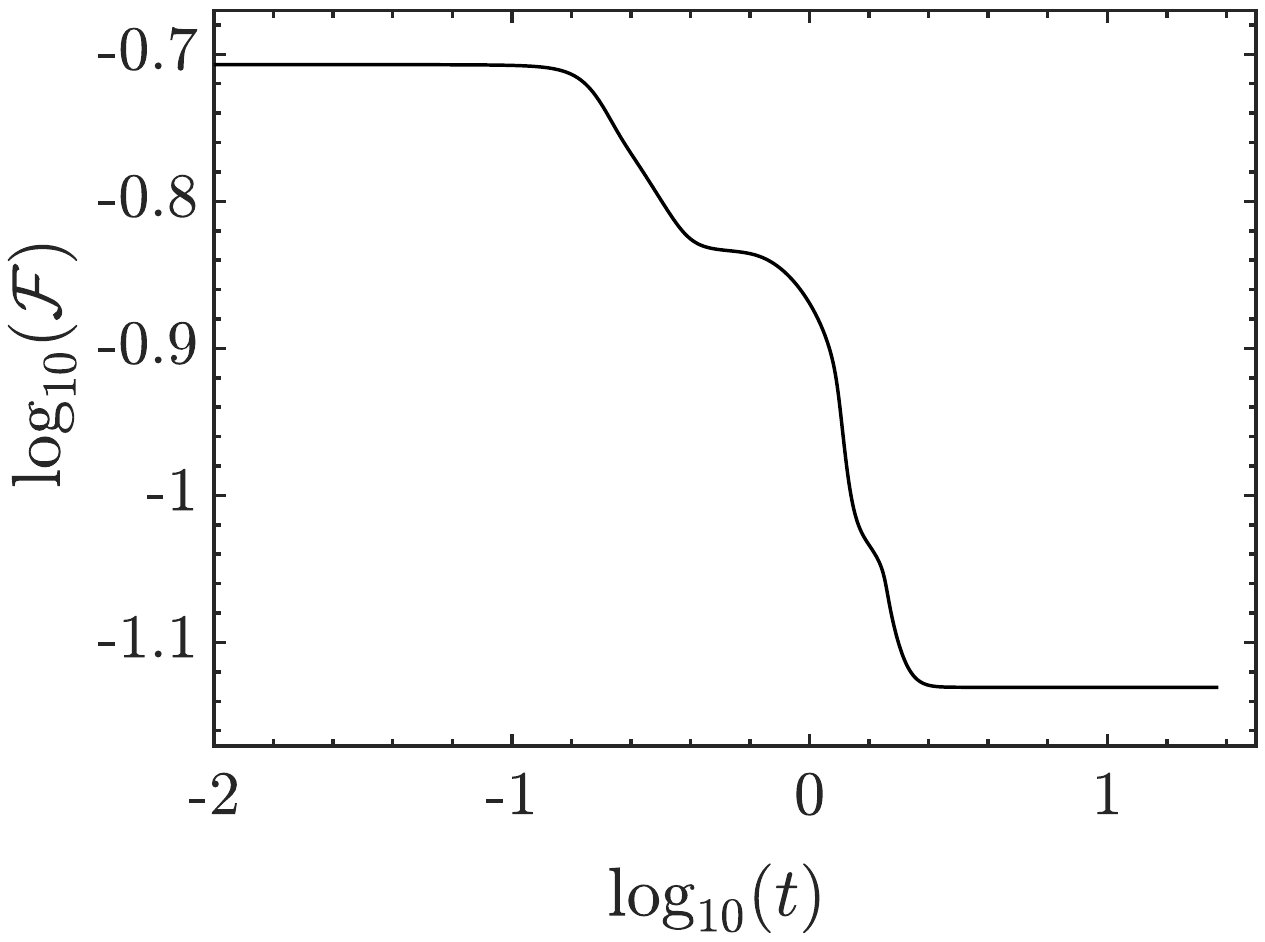}
  \caption{Evolution of the free--energy $\mathcal F$ with time for the three-dimensional spinodal decomposition simulation}
  \label{fig:Numerical:3DCyl:FreeEnergy}
\end{figure}

We run the simulation until the steady--state is reached. We show the evolution 
of the phases in Fig. \ref{fig:Numerical:3DCyl:Evolution}, where we have 
represented the interfaces and coloured the rest with blue ($\phi=-1$) and red 
($\phi=1$). Each plot in Fig. \ref{fig:Numerical:3DCyl:Evolution} represents a 
different time instant, one for each power of $10$. The initial condition corresponds to Fig. \ref{fig:Numerical:3DCyl:IC}, and the final
steady--state is Fig. \ref{fig:Numerical:3DCyl:SS}. Note that the final state is achieved with a flat interface separating both phases. Additionally, in Fig. 
\ref{fig:Numerical:3DCyl:FreeEnergy} we depict the evolution of the free--energy 
$\mathcal F$, showing that like in two dimensional simulations, it decreases monotonically from the initial condition 
to the steady state.
	
\section{Summary}\label{sec:Conclusions}
We have developed a nodal Discontinuous Galerkin (DG) spectral element method (DGSEM) to solve the Cahn--Hilliard equation. We used the  Gauss--Lobatto variant of the DG method to use the summation--by--parts simultaneous--approximation--term (SBP--SAT) property and showed that it is discretely stable. The spatial discretization uses the Bassi--Rebay 1 (BR1) method to couple inter--element fluxes, and the time discretization uses an efficient IMEX scheme. We first show the semi--discrete stability analysis (i.e. continuous in time), and later we show the fully--discrete one using the IMEX scheme. Both analyses show that the discrete free--energy is bounded in time by the initial value, which is in accordance with the continuous energy estimate. Lastly, we perform a convergence study of the scheme, we compare the scheme with previous results from the literature, and perform a simulation in three--dimensional curvilinear geometries, showing that the scheme is stable (i.e. its free--energy decreases) under all conditions tested herein 
as proved.

\acknowledgement{The authors would like to thank Dr. Gustaaf Jacobs of the San Diego State University for his hospitality. This work was supported by a grant from the Simons Foundation (\#426393, David Kopriva).
This work has also been partially supported by Ministerio de Economía y Competitividad under the research grant TRA2015-67679-C2-2-R.
The authors acknowledge the computer resources and technical assistance provided by the Centro de Supercomputaci\'on y Visualizaci\'on de Madrid (CeSViMa).
}

\appendix

\section{Stability analysis of the Crank--Nicolson scheme}
It has been shown in \eqref{eq:FullyDiscrete:IMEX-dissipation} that the implicit Euler scheme used for $\vec{Q}$ in \eqref{eq:FullyDiscrete:weakforms-with-surface-w} adds numerical dissipation that is proportional to the jumps in time $\Delta \vec{Q}$ squared, 
\begin{equation}
-\frac{1}{2}k\left\langle \mathcal J\Delta\vec{Q},\Delta\vec{Q}\right\rangle_{E,N} \leqslant 0.
\end{equation}
This dissipation can be effectively controlled using a linear combination of $\vec{Q}^{n+1}$ and $\vec{Q}^{n}$. To show this, we define an intermediate state $\vec{Q}^{\theta}$,
\begin{equation}
\vec{Q}^{\theta}=K_0\vec{Q}^{n+1} + (1-K_0)\vec{Q}^n,
\label{eq:FullyDiscrete:CN-Q}
\end{equation}
such that $K_0=1$ recovers backward Euler, $K_0=1/2$ is Crank--Nicolson, and $K_0=0$ is forward (explicit) Euler. 

To study the stability, we consider \eqref{eq:FullyDiscrete:weakforms-with-surface-w} and \eqref{eq:FullyDiscrete:Q-diff-weak}, where instead of $\vec{Q}^{n+1}$ we use $\vec{Q}^{\theta}$ defined in \eqref{eq:FullyDiscrete:CN-Q},
\begin{subequations}\label{eq:FullyDiscrete:weakforms-with-surface-CN}
\begin{align}
\left\langle \mathcal J W^{\theta},\varphi_W\right\rangle_{E,N} =& \left\langle  \left(\frac{\diff \Psi}{\diff \Phi}\right)^{n}+S_0\Phi^{n+1}-S_0\Phi^{n} ,\mathcal J\varphi_W\right\rangle_{E,N} \notag \\
&-k\int_{\partial E,N}\varphi_W\tilde{\boldsymbol{Q}}^{\star,\theta}\cdot\hat{\boldsymbol{n}}\diff S_\xi + k\left\langle\tilde{\boldsymbol{Q}}^{\theta},\nabla_\xi\varphi_W\right\rangle_{E,N}, \label{eq:FullyDiscrete:weakforms-with-surface-w-CN}\\
  \left\langle \mathcal J\frac{\Delta\vec{Q}}{\Delta t},\vec{\varphi}_Q\right\rangle_{E,N} =& \int_{\partial E,N}\left(\frac{\Delta\Phi^{\star}}{\Delta t}-\frac{\Delta\Phi}{\Delta t}\right) \tilde{\boldsymbol{\varphi}}_Q\cdot\tilde{\boldsymbol{n}}\diff S_\xi +  \left\langle \frac{\nabla_\xi\left(\Delta\Phi\right)}{\Delta t},\tilde{\boldsymbol{\varphi}}_Q\right\rangle_{E,N}. \label{eq:FullyDiscrete:weakforms-with-surface-q-CN}
\end{align}
\end{subequations}
We replace $\vec{\varphi}_{\vec{Q}} = \vec{Q}^\theta$ in \eqref{eq:FullyDiscrete:weakforms-with-surface-q-CN}, so that the left hand side is
\begin{equation}
\begin{split}
  \left\langle \mathcal J\Delta\vec{Q},\vec{Q}^{\theta}\right\rangle_{E,N} =&\left\langle \mathcal J\left(\vec{Q}^{n+1}-\vec{Q}^{n}\right),K_0\vec{Q}^{n+1} + (1-K_0)\vec{Q}^n,\right\rangle_{E,N}  \\
  =& \frac{1}{2}\left\langle \mathcal J \vec{Q}^{n+1}, \vec{Q}^{n+1}\right\rangle_{E,N}  -\frac{1}{2}\left\langle \mathcal J \vec{Q}^{n}, \vec{Q}^{n}\right\rangle_{E,N} \\
& +\left(K_0-\frac{1}{2}\right) \left\langle \mathcal J\Delta\vec{Q},\Delta\vec{Q},\right\rangle_{E,N}.
  \end{split}
  \end{equation}
When we replace $\varphi_{W} = \frac{\Delta \Phi}{\Delta t}$, we obtain the more general expression of \eqref{eq:FullyDiscrete:W+Q-equation-simpler},
  \begin{equation}
\begin{split}
\left\langle \mathcal J W^{\theta},\frac{\Delta\Phi}{\Delta t}\right\rangle_{E,N} =& \left\langle\left(\frac{\diff \Psi}{\diff \Phi}\right)^{n}+S\Delta\Phi, \mathcal J \frac{\Delta\Phi}{\Delta t}\right\rangle_{E,N}+k  \left\langle \mathcal J\frac{\Delta\vec{Q}}{\Delta t}, \vec{Q}^{\theta}\right\rangle_{E,N}  \\
&-k\text{BT}_{E,N}\left(\frac{\Delta\Phi}{\Delta t}, \vec{Q}^{\theta}\right).
\end{split}
\end{equation}
Following the rest of the steps in the fully--discrete analysis, which remain the same, we arrive to the same bound in \eqref{eq:FullyDiscrete:stability-all-timesteps}, but with new IMEX dissipation, $\text{diss}_{\text{IMEX}}^{\theta,E,N}$ in terms of $K_0$
\begin{equation}
\text{diss}^{\theta,E,N}_{\text{IMEX}} = - \left\langle \mathcal J \Pi, 1\right\rangle_{E,N} -\frac{1}{2}k\left(K_0-\frac{1}{2}\right)\left\langle \mathcal J\Delta\vec{Q},\Delta\vec{Q}\right\rangle_{E,N} \leqslant 0.
\end{equation}
Thus, the amount of dissipation related to $\vec{Q}$ added by the implicit scheme varies linearly with $K_0$, and vanishes with $K_0=1/2$ (Crank--Nicolson scheme). Note that, with this approach, we cannot confirm that the explicit Euler is stable, which does not necessarily mean that it is unstable. The stability of the explicit Euler scheme is not analyzed in detail in this work as we consider it impractical due to the time step limitation imposed by the fourth order spatial derivate of the Cahn-Hilliard equation. In the numerical experiments, we set $K_0$ to both 1 (implicit Euler) and $1/2$ (Crank--Nicolson) to see whether it is enough with the dissipation provided by $- \left\langle \mathcal J \Pi, 1\right\rangle_{E,N}$ and the physical dissipation for the approximation to be stable.

\section{Effect of interface stabilization in $\vec{Q}^\star$}\label{app:StabilizationQ}

We have only considered interface stabilization in the chemical potential 
gradient weak form, $\vec{\mathcal F}^{\star}$. We have not considered interface stabilization in the 
phase field gradient, $\vec{Q}$, since we have found that it pollutes the 
free--energy $\mathcal F$. In this Appendix we show the effect of adding 
interface stabilization also to $\vec{Q}^{\star}$,
\begin{equation}
\vec{Q}^{\star} = \average{Q} - \sigma_{q}\jump{\Phi},
\label{eq:FullyDiscreteIP:flux-Q}
\end{equation}
where $\sigma_{q}$ is a positive penalty parameter. When studying the boundary terms generated by \eqref{eq:FullyDiscreteIP:flux-Q}, 
\begin{equation}
\text{BT}_{E,N}\left(\frac{\Delta\Phi}{\Delta t}, \vec{Q}^{n+1}\right) =- \sigma_{q} \int_{\partial e, N}\jump{\Phi^{n+1}}\jump{\Delta \Phi}\diff S,
\label{eq:FullyDiscreteIP:BTEN}
\end{equation}
we find that \eqref{eq:FullyDiscreteIP:BTEN} cannot be bounded. By rearranging,

\begin{equation}
\begin{split}
-\sigma \int_{\partial e, N}\jump{\Phi^{n+1}}\jump{\Delta \Phi}\diff S =& -\frac{ \sigma_{q}}{2} \int_{\partial e, N}\jump{\Phi^{n+1}}^2\diff S+\frac{ \sigma_{q}}{2} \int_{\partial e, N}\jump{\Phi^{n}}^2\diff S \\
&-\frac{ \sigma_{q}}{2} \int_{\partial e, N}\jump{\Phi^{n+1}-\Phi^{n}}^2\diff S,
\end{split}
\end{equation}
we can see that by modifying the free--energy $\mathcal F$ to
\begin{equation}
\mathcal F^\sigma = \mathcal F + \frac{\sigma_{q}}{2}\sum_{\interiorfaces}\int_{f,N}\jump{\Phi^{n+1}}^2\diff S,
\end{equation}
then the dissipation is bounded and equal to
\begin{equation}
\text{diss}_\sigma = -\frac{\sigma_{q}}{2}\sum_{\interiorfaces} \int_{f, N}\jump{\Phi^{n+1}-\Phi^{n}}^2\diff S.
\end{equation}
Although the interface penalization contribution is always positive, and therefore the free--energy $\mathcal F$ is bounded in time, we cannot confirm that the free--energy is strictly monotonic in time. In most practical cases it is, and other methods that rely on this stabilization (e.g. the interior penalty method \cite{1978:Wheeler}) have been found by us to work well in practice with the Cahn--Hilliard 
equation.

\bibliography{mybibfile}
\end{document}